\documentclass[reqno]{amsart}

\usepackage{amsmath}
\usepackage{amsfonts}
\usepackage{rotating}
\usepackage{bbm}
\usepackage[all,cmtip]{xy}
\usepackage{amscd}
\usepackage{pgfplots}
\usepackage{tikz-cd}              
\usepackage{tikz}
 \usetikzlibrary{arrows.meta}

\usetikzlibrary{decorations.pathreplacing}

\usepackage{graphicx}
\usepackage{subcaption}
 \usetikzlibrary{arrows}
 \tikzset{>=latex}
\usetikzlibrary{arrows.meta}
\pgfplotsset{compat=1.10}
\usetikzlibrary{decorations.markings}

\usepgfplotslibrary{fillbetween}
\usetikzlibrary{patterns}
\usepackage{mathrsfs}
\usepackage{color}
\usepackage{xcolor}
\usepackage{mathtools}
\usepackage{amssymb}
\usepackage[latin1]{inputenc}
\usepackage{graphicx}
\usepackage{dsfont}
 \usepackage{hyperref}
\usepackage{enumerate}
\usepackage{marginnote}

\usepackage{subcaption}
 \usepackage{comment}
\pgfplotsset{compat=1.18}

\DeclareMathOperator{\sech}{sech}

\newtheorem{theorem}{Theorem}[section]

\newtheorem{lemma}[theorem]{Lemma}
\newtheorem{proposition}[theorem]{Proposition}

\theoremstyle{definition}
\newtheorem{definition}[theorem]{Definition}

\theoremstyle{remark}
\newtheorem{remark}[theorem]{Remark}

\newcommand{\C}{{\mathbb{C}}}

\newcommand{\R}{{\mathbb{R}}}

\newcommand{\Z}{{\mathbb{Z}}}
\newcommand{\N}{{\mathbb{N}}}

\renewcommand{\epsilon}{\varepsilon}
\renewcommand{\theta}{\vartheta}

\let\oldaddcontentsline\addcontentsline
\newcommand{\stoptocentries}{\renewcommand{\addcontentsline}[3]{}}
\newcommand{\starttocentries}{\let\addcontentsline\oldaddcontentsline}

\begin{document}

\title[The equal mass restricted three-body problem]{On the Contact Geometry of the Equal Mass Restricted Three-Body Problem}

\author[S. Kim]{Seongchan Kim}
   \address{ Department of Mathematics Education, Kongju National University, Gongju 32588, Republic of Korea}
   \email {seongchankim@kongju.ac.kr}

 \author[K.Ruck ]{Kevin Ruck}
   \address{ Seoul National University, Department of Mathematical Sciences, Research institute in Mathematics, Gwanak-Gu, 
			Seoul 08826, South Korea }
   \email {kevin.ruck@snu.ac.kr}

\setcounter{tocdepth}{3}
 
\maketitle
 \begin{abstract}

   We study the contact-type property of regularized energy hypersurfaces in the planar circular restricted three-body problem with equal masses. Extending the result of  \cite{AFvKP12}, we prove that the regularized bounded component is of contact type for every energy between the first and second critical values, denoted by $c_1$ and $c_2,$ respectively.  As applications, we derive two dynamical consequences. First, combining the contact-type property with   the two-or-infinity dichotomy  recently proved in \cite{CGHHL}, we obtain   the existence of infinitely many geometrically distinct periodic orbits. Second, we study symmetric consecutive collision orbits using equivariant Lagrangian Rabinowitz Floer homology and prove that  there are  infinitely many distinct symmetric consecutive collision orbits unless a  periodic symmetric collision orbit exists. 
  For   energies slightly above $c_2,$ we cap the regularized energy hypersurface beyond the necks at $L_2$ and $L_3$   and prove that the capped hypersurface is  of contact type. Combining equivariant Lagrangian Rabinowitz Floer homology with   local index estimates,    we show that for every $N \in \N$ there exists  $\epsilon_N>0$ such that  each    $c \in  ( H(L_2),  H(L_2)+\varepsilon_N)$ admits   at least $N$ distinct symmetric consecutive collision orbits, unless a   periodic symmetric collision orbit exists. 
   \end{abstract}

     \medskip

     \noindent
{\centerline{\small \textbf{Keywords:} Three-Body Problem $\cdot$ Contact Geometry $\cdot$ Floer Homology  }}

     \smallskip

     \noindent
{\centerline{\small \textbf{Mathematics Subject Classification:} 70F07 $\cdot$ 53D99 $\cdot$ 53D40 $\cdot$ 37J12}}

  \tableofcontents
   


\section{Introduction}

The planar circular restricted three-body problem (PCR3BP) studies the motion of a massless body moving under the gravitational influence of two massive primaries that orbit each other in circular orbits about their common center of mass, all within a single fixed plane.  
In recent decades, modern techniques from contact and symplectic geometry have provided new tools for studying the PCR3BP. The first fundamental result in this direction is due to  Albers, Frauenfelder, van Koert and Paternain \cite{AFvKP12}, who proved that the regularized bounded components   are of contact type for all energies below or  slightly above the first critical value. This result opened the way to applying tools from contact geometry, pseudoholomorphic curve theory and Floer theory to understand  the dynamics of the PCR3BP. This  viewpoint led to further developments in  \cite{AFFHvK12GSS, AFFvK13CZ,BF25WFH,CFvK14Finsler, FvKZ22RKP, frauenfelder2019a,Kang14a, kang2025consecutive, Kim18a, Kim22a, Kim24a, LS25finite, ruck24b}.
We recommend the book by Frauenfelder and van Koert \cite{FvK18book} for a comprehensive introduction to this approach.

The main purpose of this paper is to address the contact-type  property for higher energies. We restrict ourselves to the equal-mass case $\mu=\frac{1}{2}.$  In this case,  the Hamiltonian of the PCR3BP is given by
\[
H(q,p) = \frac{1}{2}\left[ (p_1 - q_2)^2 + (p_2 + q_1)^2  \right] - \frac{\frac{1}{2}}{\lvert q- \frac{1}{2}\rvert} -   \frac{\frac{1}{2}}{\lvert q+ \frac{1}{2}\rvert} - \frac{1}{2}(q_1^2 +q_2^2).
\]
The first critical point is $L_1 = (0,0,0,0),$ with critical value  $c_1:=H(L_1)=-2.$ The second critical value is realized by the two points $L_2/L_3 = (\pm q_c,0,0,\mp q_c),$ where $q_c>1$ is the unique positive solution of
\begin{equation}\label{eq:definexc}
\frac{1}{(x - \frac{1}{2})^2} + \frac{1}{( x + \frac{1}{2})^2} = 2x,  \quad \quad x>\frac{1}{2}.
\end{equation}
Numerically,   $q_c \approx 1.198406$ and $c_2:=H(L_2) =H(L_3)\approx -1.728398.$
The remaining two critical points are $L_4 =(0,   \sqrt{3}/2,   \sqrt{3}/2,0)$ and $ L_5=(0, - \sqrt{3}/2, - \sqrt{3}/2,0)$ and their common critical value is   $c_3:=H(L_4)=H(L_5)=-11/8 .$ Throughout the paper we assume   $c < c_3.$

We mainly consider energies that are higher than $c_1,$ but   below $c_2.$ In this range, the energy level $H^{-1}(c)$ contains two connected components: one is near the two primaries, and the other far away.  Due to collisions at the primaries, the component near the primaries is never compact. To apply the machinery of contact and symplectic geometry, we regularize these singularities via elliptic coordinates.

The main result of this paper is the following.

\begin{theorem}\label{thm:main} For every $H(L_1)<c<H(L_2),$ there exists a Liouville vector field that is positively transverse to the regularized bounded component. Consequently, this component  is of contact type. 
\end{theorem}

The proof relies crucially on  elliptic coordinates. An advantage of these coordinates is that  the regularized effective potential becomes decoupled, namely, it is the sum of two single-variable smooth functions. This separation property  is the main technical tool of the paper. 

\begin{remark}
The use of  elliptic coordinates in the present paper  is inspired by the recent work of Liu and Salom\~ao \cite{LS25finite}, where   elliptic coordinates are used to  verify  strict convexity of the bounded components for energies slightly above  $H(L_1).$ This convexity is an essential ingredient in their  construction of  finite energy foliations.  For sufficiently small mass ratios, they also construct finite energy foliations between the first and second critical values and on some compact region containing the primaries for energies slightly above the second critical value, see \cite{LS26}.
\end{remark}

\begin{remark} The contact-type  property  established in \cite{AFvKP12}  was later extended to the spatial problem by Cho, Jung and Kim \cite{CJK20} using a similar approach. 
However, the argument in the present paper using  elliptic coordinates is essentially planar. The natural three-dimensional analogue of the planar elliptic coordinates is given by the   prolate spheroidal coordinates $(x_1, x_2, \varphi) \in \R \times \R / 2 \pi \Z \times \R / 2 \pi \Z$:
    \[
    q_1 =\frac{1}{2}\cosh x_1 \cos x_2, \quad q_2 = \frac{1}{2}\sinh x_1 \sin x_2 \cos \varphi, \quad q_3 = \frac{1}{2}\sinh x_1 \sin x_2 \sin \varphi
    \]
An issue is that in dimension three the additional angular variable $\varphi$ introduces new artificial singularities, not corresponding to collisions. These occur along $
x_1=0,$    $ 0 < x_2 < \pi,$ 
which is mapped to the line segment joining  the two primaries. 
Thus, the elliptic coordinates argument developed in this paper does not extend directly to the spatial problem.
\end{remark}

 As an application of Theorem \ref{thm:main}, we obtain infinitely many periodic orbits.

\begin{theorem}\label{thm:infinitelymany} For every $H(L_1)<c<H(L_2),$ the regularized bounded component  carries infinitely many geometrically distinct periodic orbits.
\end{theorem}
\begin{proof}
For $c \in ( H(L_1), H(L_2)),$ the Moser regularized bounded component is diffeomorphic to   $ \R P^3 \# \R P^3.$ This is not a lens space and satisfies $H^2(\R P^3 \# \R P^3;\Z) = \Z_2 \oplus \Z_2.$ Thus, the first Chern class of the induced contact structure is torsion.  The claim therefore follows from     \cite[Corollary 1.7]{CGHHL}. 
\end{proof}

The contact-type property established in Theorem~\ref{thm:main} also allows us to study \textit{consecutive collision orbits}. These are trajectories  between two  collisions with the same primary, described after regularization by chords with endpoints on its corresponding collision Lagrangian.   The motivation to study this kind of orbit  comes from the interest in finding suitable orbits to perform a powered flyby maneuver. For such a maneuver it is crucial to fly as close as possible to the primary. Once we find a consecutive collision orbit   we can perturb the orbit slightly to get an orbit that starts and ends in an  almost collision  with the primary, i.e.~an orbit one can use to perform a powered flyby. The advantage of considering consecutive collision orbits is that they can be described as Reeb chords between Legendrian submanifolds, allowing  us to use well-established tools from symplectic geometry.

The singularities at $q = (\pm \frac{1}{2},0)$ are resolved by regularization, in which a collision with any of the two primaries is represented by a Lagrangian submanifold.   There is also a natural anti-symplectic involution,   given by the reflection of the position  across the $q_1$-axis and reflection of the momentum across the $p_2$-axis.   A symmetric consecutive collision orbit is then a consecutive collision orbit   invariant under this involution. Equivalently,   half of such an orbit can be viewed as  a   Hamiltonian chord  that starts on the   collision Lagrangian and ends on the fixed point set of the anti-symplectic involution.
Using    $\Z_2$-equivariant Lagrangian Rabinowitz Floer homology (RFH), we prove the following result. Each assertion holds for collisions with either prescribed primary.

\begin{theorem}\label{thm:ccorbits} The following statements hold:
\begin{itemize}
    \item[$i)$] For all energies $H(L_1)<c<H(L_2) ,$ there exists either a periodic symmetric collision orbit or infinitely many distinct symmetric consecutive collision orbits.
    \item[$ii)$] For every $N \in \N,$ there is an $\epsilon_N>0$ such that, for all energies $c\in (H(L_2),H(L_2)+\epsilon_N),$ we have either a periodic symmetric collision orbit or at least $N$   distinct symmetric consecutive collision orbits.
\end{itemize}
\end{theorem}

Heuristically periodic symmetric collision orbits are highly special: they satisfy several   boundary and symmetry conditions simultaneously. Hence, one expects that  they should not persist under a small perturbation of the energy. A full proof of this expectation remains open. 
However, according to \cite[Proposition~4.3 and Corollary~4.4]{kang2025consecutive},
at a generic energy slightly above the first critical value and for a generic mass ratio, any symmetric periodic collision orbit through the prescribed primary has an even number of such collisions, implying that it
contains at least two geometrically distinct symmetric consecutive collision orbits.

For the first part of Theorem \ref{thm:ccorbits}  we related elliptic and Birkhoff regularization by constructing an explicit symplectomorphism.  An equivariant contact isotopy then identifies the subcritical hypersurfaces, together with the collision and symmetry loci, with the contact model due to Kang and Ruck  \cite{kang2025consecutive}, provided energies are slightly above $c_1.$ This model RFH computation gives the required chords.  Then the symplectomorphism allows us to transfer this computation to our case. Finally, by   homotopy invariance of RFH, we   extend the result to all energies $c\in (c_1,c_2)$.

The second part requires an additional construction.  For $c>c_2,$ the energy level becomes non-compact and therefore, the argument used for   subcritical energies cannot be applied directly.  Instead, we modify the Hamiltonian beyond the neck regions   corresponding to the critical points $L_2$ and $L_3,$  leaving a neighborhood of each neck unchanged. This  yields a compact  energy level $\tilde{\Sigma}_c$  diffeomorphic to $S^1 \times S^2.$ We then construct a Liouville vector field transverse to $\tilde{\Sigma}_c.$  The four neck equilibria will be canceled, so that the resulting slightly perturbed new Liouville vector field expresses $\tilde{\Sigma}_c$  as a graph over a fixed subcritical hypersurface and gives an equivariant RFH continuation isomorphism. Finally, adapting the local index estimates from \cite[Section 7]{LS25finite}, we show that trajectories passing through the neck regions have arbitrarily large positive index as they the energy decreases to $c_2.$ We compare this index with the Floer grading and choose generators in a finite range of degrees. Their Hamiltonian segments are contained in the unmodified region, hence, provide the desired collision orbits. 

\stoptocentries
\subsection*{Acknowledgments}  SK acknowledges   support  by Samsung Science and Technology Foundation under Project Number SSTF-BA2502-02. 
KR was supported by Samsung Science and Technology Foundation under Project Number SSTF-BA2502-02 and the National Research Foundation of Korea under the grants RS-2023-00211186, (MSIT) RS-2023-NR076656 and  NRF-2020R1A5A1016126.

\subsection*{Use of AI Tools} 
 The authors used large language models   to proofread the manuscript, refine the language and exposition, generate figures, and assist with the review.  All mathematical statements, derivations, proofs,  and references were independently checked and approved by the authors, who take full responsibility for the content of the article.

\starttocentries


 \section{The planar circular restricted three-body problem}


\subsection{Elliptic coordinates}
\label{sec:ellipticcoord}

Elliptic coordinates $ (x_1, x_2) \in \R \times \R/2\pi \Z $ are given by
\[
q_1 = \frac{1}{2} \cosh x_1  \cos x_2 , \quad   q_2 = \frac{1}{2} \sinh x_1  \sin  x_2  .
\]
Let $(y_1, y_2)\in \R \times \R$ denote the conjugate momenta determined by    $ydx= pdq$. Explicitly, one has 
\[
\begin{pmatrix} p_1 \\ p_2 \end{pmatrix} = \frac{2}{\cosh^2 x_1 - \cos^2 x_2}\begin{pmatrix} \sinh x_1 \cos x_2 & -\cosh x_1 \sin x_2 \\ \cosh x_1 \sin x_2 & \sinh x_1 \cos x_2 \end{pmatrix} \begin{pmatrix}
  y_1 \\ y_2  
\end{pmatrix}.
\]
For a fixed energy value   $c \in \R,$   define the regularized Hamiltonian  
\[
  K_c(x,y) :=   \frac{1}{4}\left( \cosh^2 x_1   - \cos^2  x_2  \right) \left( H(x,y) - c\right).
\]
A direct computation yields
\begin{equation}\label{eq:regulHam}
    K_c(x,y)  =    \frac{1}{2}\left[ \left( y_1 + A_2(x_2)\right)^2 + \left( y_2 +A_1(x_1)\right)^2     \right] +V_c(x ),
\end{equation}
where 
\[
A_1(x_1)=\frac{1}{8} \sinh(2x_1), \quad  A_2(x_2) =  \frac{1}{8} \sin(2x_2) \]
and 
\begin{equation}\label{eq:Potential}
    V_c(x ) = - \frac{1}{128}\left( \sinh^2(2x_1) + \sin^2(2x_2) \right) - \frac{c}{4}\left( \cosh^2 x_1  - \cos^2 x_2  \right) - \frac{1}{2}\cosh x_1 .
\end{equation}
It is convenient to decompose $V_c$ as  $V_c(x_1, x_2) = V_{c,1}(x_1)+V_{c,2}(x_2),$ where
\begin{equation}\label{eq:decoupled}
\begin{aligned}
    V_{c,1} (x_1)&=   - \frac{1}{128} \sinh^2(2x_1)   - \frac{c}{4}  \cosh^2 x_1    - \frac{1}{2}\cosh x_1 ,\\
    V_{c,2}(x_2)&=    - \frac{1}{128}  \sin^2(2x_2) +\frac{c}{4} \cos^2 x_2  .
\end{aligned}
\end{equation}
The regularized system  has   an additional $\mathbb{Z}_2$-symmetry given by
 \begin{equation}
     \begin{aligned}
         \sigma_{\rm el} \colon \mathbb{Z}_2\times \R \times \R/2\pi \Z\times \R \times \R &\to \R \times \R/2\pi \Z\times \R \times \R; \\ (a,(x_1,x_2,y_1,y_2)) &\mapsto (a x_1,a x_2,a y_1,a y_2).
     \end{aligned}
     \label{eq:ellipticSym}
 \end{equation}
 A direct computation shows that $q(x )$, $p(x,y)$ and the symplectic form $\omega = dy \wedge dx$ are invariant under this action.



\subsection{Hill's regions}\label{sec:Hill}

We first analyze the potential $V_{c,1}.$ We have 
\begin{equation}\label{eq:derofV1}
\begin{aligned}
    V_{c,1}'(x_1) 
    &= -\frac{1}{2} \sinh x_1 ( c \cosh x_1 +1 +\frac{1}{8}\cosh x_1 \cosh 2x_1 ), \\
       V_{c,1}''(x_1)  &= \frac{1}{16} ( -8 c \cosh 2x_1 - 8 \cosh x_1 - \cosh 4x_1).
\end{aligned}
\end{equation}
Thus $x_1=0$ is a critical point for every $c \in \R .$ Moreover, $V_{c,1}''(0)  >0$ for  $c<c_3 ,$ so $x_1=0$ is  a strict local minimum. 

To find the other critical points, set $u   = \cosh x_1 \geq 1.$ Since $\cosh 2x_1 = 2u^2-1, $ every nonzero critical point of $V_{c,1}$ corresponds to a root of the cubic polynomial 
 \begin{equation}\label{eq;polynomialP}
P_c(u)  =  2u^3 + (8c-1)u+8.
\end{equation}
Denote by
\begin{equation}\label{eq:u0}
u_0=u_0(c) = \sqrt{ \frac{1-8c}{6}} > \sqrt{2}
\end{equation}
the unique positive critical point of $P_c,$ so that $P_c$ is strictly  decreasing on $[1,u_0)$ and is strictly  increasing on $(u_0,\infty).$ We then have
\[
P_c(u)  < P_c(1) = 8c+9 <0, \qquad u \in (1,  u_0] 
\]
from which we see that $P_c$ has no root in $[1, u_0].$ Since $P_c \to \infty$ as $u \to \infty,$ it follows that   $P_c$ has a unique root $u_*=u_*(c) \in (u_0,\infty) .$  Note that $P_c$ is negative on $[1,u_*)$ and positive on $(u_*, \infty).$
If we take $x_*=x_*(c)>0$ satisfying $u_* = \cosh x_*,$ then $\pm x_*$ are the two nonzero critical points of $V_{c,1}.$  Since $V_{c,1} \to -\infty$ as $x_1 \to \pm \infty,$ these critical points are local maxima.

We next consider $V_{c,2}.$ A direct computation gives 
\begin{equation}\label{eq:derivofV2}
\begin{aligned}
V_{c,2}'(x_2)&= - \frac{1}{32} \sin 2x_2 ( 8c+ \cos 2x_2), \\   V_{c,2}''(x_2)&= -\frac{1}{16} ( 8c \cos 2x_2 + \cos 4 x_2).
\end{aligned}
\end{equation}
Since $c< c_3,$   the only critical points are $x_2=0, \frac{\pi}{2}, \pi, \frac{3\pi}{2}.$   The  points $x_2=0, \pi$ are local minima, while   $x_2=\frac{\pi}{2}, \frac{3\pi}{2}$ are local maxima.

Combining the critical points of $V_{c,1}$ and $V_{c,2}$, we obtain the following critical points of the potential  $V_c = V_{c,1} +V_{c,2}$: 
\begin{itemize}
    \item[-] $(0,0)$ and $(0,\pi)$ are strict local minima with value $-\frac{1}{2}$;

    \item[-] $(0, \frac{\pi}{2})$ and $(0, \frac{3\pi}{2})$ are saddles with value  $-\frac{c}{4}-\frac{1}{2}$;

    \item[-] $(\pm x_*, 0)$ and $(\pm x_*, \pi)$ are saddles with value  $V_{c,1}(x_*)+\frac{c}{4}$; and

    \item[-] $(\pm x_*, \frac{\pi}{2})$ and $(\pm x_*, \frac{3\pi}{2})$ are strict local  maxima with value $V_{c,1}(x_*).$
    
\end{itemize}

\begin{lemma}\label{lem:orderofvalues}
    For $c\in(H(L_1),H(L_2)+\varepsilon_0)$  with   $\varepsilon_0>0$ sufficiently small,  the four critical values are ordered as
\[
-\frac{1}{2} < - \frac{c}{4}-\frac{1}{2} < V_{c,1}(x_*)+\frac{c}{4}< V_{c,1}(x_*).
\]

\end{lemma}

\begin{proof} The first and last inequalities follow from the fact that $c<0.$ 
 To study  the second inequality, we  introduce  
\begin{equation}\label{eq:B}
B_c(u):=   \frac{1}{32}u^2(1-u^2 )+\frac{  c}{4}(1-u^2 )-\frac{u}{2},
\end{equation}
so that
\[
B_c(u) = V_{c,1}(x_1)+\frac{c}{4}, \qquad u= \cosh x_1.
\]

The function   $B_*(c) := B_c( u_*(c))$ is   smooth. Indeed, $P_c'(u_*(c))>0$ by the fact that $u_*(c) >u_0(c),$ so the implicit function theorem gives smooth dependence of $u_*$ on $c.$ It follows from uniqueness that $B_*$ is smooth.  

We next claim that the function 
\[
\mathfrak{D}(c) := B_c( u_*(c)) + \frac{c}{4}+\frac{1}{2} = B_*(c) + \frac{c}{4} +\frac{1}{2} 
\]
  is strictly decreasing for every $c <c_3.$ To see this, we  compute
\begin{equation}\label{eqdiffofD}
\begin{aligned}
\mathfrak{D}'(c) &= \partial_c B_c( u_*(c)) + \partial_u B_c(u_*(c))u_*'(c)+\frac{1}{4} \\
&= \partial_c B_c( u_*(c)) - \frac{1}{16} P_c( u_*(c))u_*'(c)+\frac{1}{4} \\
&= \frac{1}{4}(2-u_*(c)^2).
\end{aligned}
\end{equation}
Now since $u_*(c) >u_0(c)>\sqrt{2},$ see \eqref{eq:u0}, it follows that $\mathfrak{D}'(c)<0$ for every $c <c_3.$

Recall that   $q_c>1,$ the $q_1$-component of $L_2,$ is the unique solution of  equation  \eqref{eq:definexc}. Denoting by $\hat{u}= 2q_c>2,$ we get from \eqref{eq:definexc} that
\[
\frac{4}{(\hat{u} -1)^2} + \frac{4}{ (\hat{u}+1)^2} = \hat{u},
\]
or equivalently,
\begin{equation}\label{eq:ustaridentty}
\hat{u}^5 - 2\hat{u}^3 - 8\hat{u}^2 +\hat{u} -8=0.
\end{equation}
We now observe that
\begin{equation}\label{eq:lasthl2}
 c_2   = - \frac{\frac{1}{2}}{\lvert q_c- \frac{1}{2}\rvert} -   \frac{\frac{1}{2}}{\lvert q_c+ \frac{1}{2}\rvert} - \frac{1}{2} q_c^2   = - \frac{1}{ \hat{u}-1} - \frac{1}{\hat{u}+1} - \frac{1}{8}\hat{u}^2 = \frac{ -2\hat{u}^3 + \hat{u}-8}{8\hat{u}},
\end{equation}
where in the last identity we used \eqref{eq:ustaridentty}. We thus have
\[  
 8 c_2-1 = \frac{ -2\hat{u}^3 + \hat{u}-8}{\hat{u}} -1= \frac{ -2\hat{u}^3 -8}{\hat{u}} \quad \Rightarrow \quad P_{c_2}(\hat{u})=0.
\]
Since $u_*(c_2)>1$ is the unique root of $P_{c_2} $ in $(1,\infty),$ it follows that $\hat{u}=u_*(c_2).$
Plugging \eqref{eq:lasthl2} into $B_c$ yields
\[
B_{c_2}(\hat{u}) = \frac{\hat{u}^5-2\hat{u}^3-8\hat{u}^2+\hat{u}-8}{32\hat{u}}=0,
\]
where we  used again \eqref{eq:ustaridentty}. It follows that
\begin{align*}
\mathfrak{D}(c_2) = B_{c_2}(  u_*(c_2)) + \frac{1}{4}c_2 + \frac{1}{2} =  \frac{1}{4} c_2 + \frac{1}{2} >0.
\end{align*}
  Since $\mathfrak{D}$ is strictly decreasing, this implies that $\mathfrak{D}(c)>0$ for every $c_1<c<c_2+\varepsilon_0 $ with $\varepsilon_0>0$ sufficiently small. This completes the proof. 
\end{proof}


\begin{remark}\label{remark:valuev1xstar}
The argument in the proof above shows that 
\[
B_*(c)=B_c(u_*(c))=V_{c,1}(\pm x_*)+\frac{c}{4}=V_c(\pm x_*,0)
\]
is positive for $c \in (c_1, c_2)$   and equals zero at $c = c_2.$ Indeed, it is shown that 
\[
\mathfrak{D}(c) = B_*( c ) + \frac{c}{4}+\frac{1}{2}
\]
is strictly decreasing for $c<c_3 .$ We thus have $B_*'(c)< - \frac{1}{4},$ and hence $B_*$ is strictly decreasing as well. Moreover, it was shown that $B_*( c_2)=0.$ In particular, we have $B_*(c)<0$ for $c>c_2$ sufficiently close to $c_2.$

\end{remark}

\begin{remark}
It is easy  to show that   $u_*(c_3)=2,$ so $\mathfrak{D}(c_3)=-3/16<0.$ Since $\mathfrak{D}(c_2)>0$ and  $\mathfrak{D}$ is strictly decreasing for $c < c_3,$ this shows that  there is a unique $\hat{c} \in (c_2,c_3)$ such that $\mathfrak{D}(\hat{c})=0.$ Hence,
$V_{c,1}(x_*)+c/4 <  - c/4-1/2$ for $c \in (\hat{c}, c_3) .$
\end{remark}

Assume first that   $c<c_1=-2,$ so that  $- \frac{c}{4} - \frac{1}{2} >0$.  Hence the Hill's region 
\[
\mathfrak{K}_c:= \{ (x_1, x_2) \mid V_c(x_1, x_2) \leq 0 \}
\]
associated with the    regularized energy level $K_c^{-1}(0)$   consists of four connected components: $\mathfrak{K}_{c}^0,\mathfrak{K}_{c}^\pi $ and $\mathfrak{K}_c^{u,\pm}.$ Here   $\mathfrak{K}_{c}^0$ is the bounded component around $(0,0),$ $\mathfrak{K}_{c}^\pi$ is the bounded component  around $(0,\pi)$ and $\mathfrak{K}_c^{u,\pm}$ are the unbounded components, corresponding to the two ends $x_1 \to \pm \infty$. The two bounded components are symmetric with respect to the line $x_2 = \frac{\pi}{2}.$ Each of them is   diffeomorphic to a closed disk.  The corresponding bounded component of   $K_c^{-1}(0)$ lying over either disk  is diffeomorphic to $S^3.$
At $c=c_1,$ the two bounded components meet precisely at $(0, \frac{\pi}{2})$ and $(0, \frac{3\pi}{2}).$ 
For  $c_1<c<c_2,$ the two bounded components are merged into a single bounded component, denoted by $\mathfrak{K}_{c}^b$, while the two unbounded components remain disjoint. The   bounded component $\Sigma_c$ of $K_c^{-1}(0),$ lying over $\mathfrak{K}_c^b,$ is diffeomorphic to $S^1 \times S^2.$  Its quotient by $\sigma_{\rm el}$ is diffeomorphic to $\R P^3 \# \R P^3.$  At $c = c_2,$  the  component $\mathfrak{K}_{c}^b$ meets the unbounded components precisely at $(\pm x_*, 0)$ and $(\pm x_*, \pi).$ Finally, for $c_2 < c < c_3=- {11}/{8} ,$ the bounded and unbounded components are connected to form a single unbounded component. 
 
\begin{figure}[ht!]
    \centering
        \includegraphics[width=\textwidth]{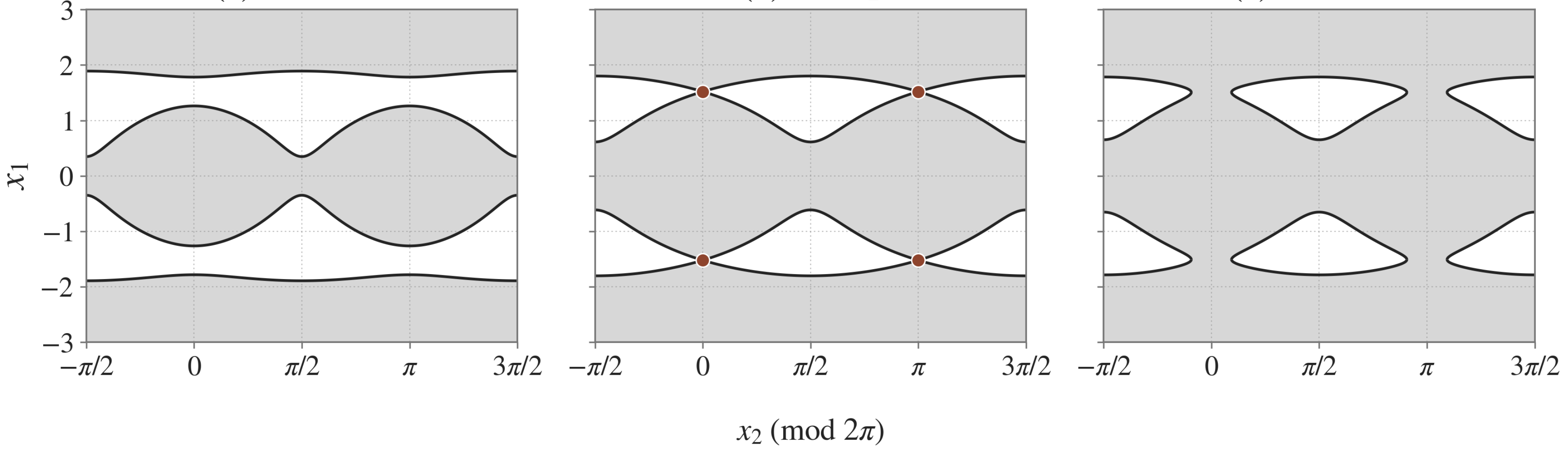}
        \caption{Hill's regions (shaded) in elliptic coordinates for (a) $c_1<c<c_2,$ (b) $c=c_2,$ and (c) $c_2<c<c_3.$ The four dots in (b) indicate the saddle points where the bounded region intersects the two unbounded regions}
        \label{fig:hill}
\end{figure}



\section{Below the second critical value: proof of Theorem \ref{thm:main}}
\label{sec:Proof1.1}
Throughout this section, we assume that $c_1<c<c_2.$
Motivated by the proof of \cite[Proposition 3.7]{dPKSS24}, we define  
\begin{equation}\label{eq:yepsilon}
Y_1 = Y_0 +  X_f 
\end{equation}
 where  
 \[
 Y_0 = y_1 \partial_{y_1} + y_2 \partial_{y_2} , \qquad 
X_f=\frac12\tanh(2x_1)\partial_{x_1}-\sech^2(2x_1)y_1\partial_{y_1}.
\] 
Here, $X_f$
denotes the Hamiltonian vector field generated by the Hamiltonian \[
f = \frac{1}{2} \tanh(2x_1) y_1.\] Both Liouville vector fields $Y_0$ and $Y_1$ are   invariant under   $\sigma_{\rm el},$  see \eqref{eq:ellipticSym}.

Set $p_1 =y_1+A_2(x_2)$ and $p_2 =y_2+A_1(x_1).$ Then we have
\[
K_c(x,y)=\frac{1}{2}(p_1^2+p_2^2)+V_c(x_1,x_2).
\]
On \(K_c^{-1}(0)\), we therefore have the constraint
\begin{equation}\label{eq:constraint}
p_1^2+p_2^2=-2V_c(x_1,x_2)=:R_c^2(x_1,x_2)  \qquad (R_c \geq 0)
\end{equation}
We compute
\begin{equation}\label{eq:Y1Ksecond}
\begin{aligned}
  Y_1(K_c) &= Y_0(K_c) + X_f(K_c)\\
  &= p_2^2+\tanh^2(2x_1)(p_1^2-A_2p_1)+\frac12\tanh(2x_1)V_{c,1}'(x_1) \\
&= R_c^2 - \sech^2(2 x_1)p_1^2 - \tanh^2(2x_1)A_2p_1+\frac12\tanh(2x_1)V_{c,1}'(x_1),
\end{aligned}
\end{equation}
where the last identity is obtained by the constraint $ p_2^2=R_c^2-p_1^2$, see \eqref{eq:constraint}. As a function of $p_1 \in [-R_c,R_c],$ this is a concave quadratic function. Therefore its minimum is attained at $p_1 = \pm R_c,$ and we thus have 
\begin{equation}\label{eq:Y1Kfirst}
\begin{aligned}
Y_1(K_c)&\ge \tanh^2(2x_1)(R_c^2-R_c|A_2|)+\frac12\tanh(2x_1)V_{c,1}'(x_1)\\
&\ge \tanh^2(2x_1)\left(\frac{V_{c,1}'(x_1)}{2\tanh(2x_1)}-\frac{A_2(x_2)^2}{4}
\right).
\end{aligned}
\end{equation}

We claim that 
\[
\Psi_c(x_1, x_2):= \frac{V_{c,1}'(x_1)}{2\tanh(2x_1)}-\frac{A_2(x_2)^2}{4}  
\]
is strictly positive  on the bounded component $\mathfrak{K}_{c}^b$. We   set $u=\cosh x_1\ge1$ and $v=\sin^2x_2 \in [0,1],$
so that 
\[
\frac{A_2(x_2)^2}{4}=\frac{1}{64}v(1-v), \qquad  
\frac{V_{c,1}'(x_1)}{2\tanh(2x_1)} = -\frac{2u^2-1}{64u}P_c (u) =: \frac{1}{64}F_c(u),
\]
where the polynomial $P_c$ is as in \eqref{eq;polynomialP}. At $x_1=0,$ the quotient $V_{c,1}'(x_1) /2\tanh(2x_1)$   is understood by the continuous extension 
\[
\lim_{x_1\to0}\frac{V_{c,1}'(x_1)}{2\tanh(2x_1)}
=
-\frac{8c+9}{64}>0.
\]
Hence, we have
\[
  \Psi_c(x_1, x_2) = \frac{1}{64}\left( F_c(u)-v(1-v) \right).
\]
 The condition $V_c=V_{c,1}+V_{c,2} \leq 0$ is equivalent to 
\begin{equation}\label{eq:Vintouv}
B_c (u)+ \frac{v^2}{32} + av\le0, \qquad  a = a(c)=-\frac{  c}{4}-\frac{1}{32}>0,
\end{equation}
where $B_c$ is as in \eqref{eq:B}.
Since $ v(1-v) \leq v \leq -B_c(u)/a,$  we find 
\[
\Psi_c(x_1, x_2) \geq \frac{1}{64}\left( F_c(u) +\frac{B_c(u)}{a}\right).
\]
Thus,   it now suffices to show that $\mathcal{Q}_c(u):=aF_c(u) + B_c(u)>0.$ 
Recall from  Section \ref{sec:Hill} that $P_c$ has the unique root $u_*=u_*(c)>1 $ and that    it  is negative on $[1,u_*)$ and positive on $(u_*, \infty).$ 
Since $B_c'(u)=-(1/16)P_c(u),$ it follows that $B_c$ is strictly increasing on $[1,u_*)$ and strictly decreasing on $(u_*,\infty).$ Since $B_c(1)=- \frac{1}{2}<0$ and 
\[
B_c(u_*) = V_{c,1}(x_*)+\frac{c}{4} = V_c(\pm x_*,0)>0, 
\]
provided $c<c_2 $ (see Remark \ref{remark:valuev1xstar}), and $B_c (u) \to -\infty$ as $u \to \infty,$  there are exactly two zeros $u_-(c) \in (1, u_*(c))$ and $u_+(c) \in (u_*(c),\infty)$ of $B_c.$   We also have $B_c(u) \leq 0$ if and only if $u \in [1,u_-(c)]\cup [u_+(c),\infty).$ 
    Since $v \geq 0,$ inequality \eqref{eq:Vintouv} implies that $B_c(u) \leq 0$   on $\mathfrak{K}_{c}.$  Hence, no point of $\mathfrak{K}_{c} $ can satisfy $  u_-(c)<u<u_+(c),$  otherwise such a $u$ would have $B_c(u)>0$,   a contradiction.      Recall that the bounded component $\mathfrak{K}_{c}^b$ contains $u=1,$ or equivalently, $x_1=0.$  Therefore, the image of $\mathfrak{K}_c^b$ under the continuous map $(x_1,x_2) \mapsto \cosh x_1$ is a connected subset of $[1,\infty) $ which contains $1.$ Therefore, it does not intersect $[u_+(c),\infty),$ and hence we have 
    \[
    \mathfrak{K}_c^b \subset \{ 1 \leq u \leq u_-(c)\}.
    \]
In particular, $u_{\max}(c) = u_-(c)$   is the maximal $u$-coordinate attained in  $\mathfrak{K}_{c}^b.$  Hence, it is enough to consider $u \in [1, u_{\max}].$

We compute
 \[
 \mathcal{Q}'_c(u)  = \frac{1}{16u^2} \mathcal{D}_c(u),  
\]
where $     \mathcal{D}_c(u)=   (6+64c)u^5 +(128c^2-24c-3)u^3 +64c u^2 +4+32c .$
To show that  $\mathcal{D}$ has a unique zero in $(1,\infty),$   we differentiate it to obtain
\[
\mathcal{D}_c'(u) = u\mathcal{E}_c(u), \quad \mathcal{E}_c(u)  = 5(6+64c)u^3 +3(128c^2-24c-3)u  +128c .
\]
We further compute
\[
\mathcal{E}_c'(u) = 15(6+64c)u^2 + 3(128c^2 - 24c -3)
\]
from which we see that $\mathcal{E}'_c$ is strictly decreasing in $u>0.$ Therefore, for all $u \geq 1$  
\[
\mathcal{E}_c'(u) \leq \mathcal{E}'_c(1) = 384c^2 + 888c +81 <0,  
\]
where the inequality follows from $c \in  (c_1, c_3), $
and thus $\mathcal{E}_c$ is strictly decreasing on $[1,\infty).$ 
We also find
\[
\mathcal{E}_c(1) = 384c^2 + 376c+21 >0, \qquad \mathcal{E}_c(u) \to -\infty \; \text{ as } \; u \to \infty,
\]
and hence there exists a unique $\bar{u} =\bar{u}(c)>1$ such that $\mathcal{E}_c(\bar{u})=0.$ Consequently,  $\mathcal{D}_c'(u) > 0$ if $1 < u < \bar{u}$ and $\mathcal{D}_c'(u)  <0$ and if $u > \bar{u},$
 which implies that $\mathcal{D}_c$ first increases and then decreases on $[1,\infty).$ Note that $\mathcal{D}_c \to -\infty$ as $ u \to \infty.$ 
Since
\[
\mathcal{Q}'_c(1) = \frac{1}{16}\mathcal{D}_c(1) = \frac{128c^2 + 136c+7}{16}>0,
\]
this implies that $\mathcal{Q}_c$ is either strictly increasing on $[1,u_{\max}]$ or is increasing first and then decreasing. In either case, the minimum of $\mathcal{Q}_c$ on $[1,u_{\max}]$ is attained at one of the endpoints.  Note that 
\begin{align*}
\mathcal{Q}_c(1) &= aF_c(1)+B_c(1) = \frac{64c^2+ 80c-7}{32}>0 , \\
\mathcal{Q}_c(u_{\max})&=aF_c(u_{\max}) = -a \frac{2u_{\max}^2-1}{u_{\max}}P_c(u_{\max})>0
\end{align*}
 since    $P_c(u)<0$ for all $1\leq u < u_*$ and $u_{\max} <u_*.$ We then conclude that $\mathcal{Q}_c(u)>0$ for all $u \in [1,u_{\max}],$ proving the claim. 

Consequently, 
\[
Y_1(K_c) \geq0 \quad \text{ on }\;\; \Sigma_c.
\]
Since $\Psi_c>0$   on $\mathfrak{K}_{c}^b,$ the lower bound  in   \eqref{eq:Y1Kfirst} is positive whenever $ x_1 \neq 0.$ Therefore, $Y_1(K_c)$ can vanish only when $x_1=0.$  Plugging this into \eqref{eq:Y1Ksecond} yields $Y_1(K_c)=p_2^2.$ Thus, $Y_1(K_c)=0$ only along $\Sigma_{c,0} \subset\Sigma_c,$ where
 \[
 \Sigma_{c,0} = \{ x_1=p_2=0 \} \cap \Sigma _c.
 \]

 We now claim that $Y_0(K_c) >0$ along $\Sigma_{c,0}.$  Indeed, on $\Sigma_{c,0}$ we have
\[
Y_0(K_c) = p_1^2 - A_2p_1  .
\]
Since
\[
p_1^2 = -2V_c(0,x_2) = \frac{1}{64} \sin^2(2x_2) + \frac{c}{2} \sin^2 x_2+1 >0,
\]
 we find that $p_1 \neq0$ and 
 \[
 p_1^2 - A_2^2 = 1+\frac{c}{2} \sin^2 x_2 \geq 1 +\frac{c}{2}>0 .
 \]
  We then conclude that
  \[
  Y_0(K_c) \geq \lvert p_1 \rvert ( \lvert p_1 \rvert  -\lvert A_2 \rvert) >0  \quad \text{ on } \Sigma_{c,0},
  \]
proving the claim.

For $\theta \in (0,1),$ we define the $\sigma_{\rm el}$-invariant Liouville vector field
\[
Y_\theta:= (1-\theta)Y_0 + \theta Y_1 = Y_0 + \theta X_f
\]
and claim that it is positively transverse to $\Sigma_c$ for $\theta$ sufficiently close to $1.$  
Indeed, by continuity the previous claim shows that there is an open neighborhood $\mathcal{U}_c \subset \Sigma_c$ of $\Sigma_{c,0}$ and a positive constant $m_c>0$ such that
\[
Y_0(K_c) \geq m_c \qquad \text{ on } \mathcal{U}_c.
\]
Since $Y_1(K_c) \geq0,$  we then obtain
\[
Y_\theta (K_c)=(1-\theta)Y_0(K_c) +\theta Y_1(K_c)  \geq (1-\theta)m_c >0  \qquad \text{ on } \mathcal{U}_c.
\]
Since $\Sigma_c \setminus \mathcal{U}_c$ is compact and disjoint from $\Sigma_{c,0},$ where $Y_1(K_c)$ vanishes, there is a positive constant $b_c$ such that 
\[
Y_1(K_c) \geq b_c \quad \text{ on } \Sigma_c \setminus \mathcal{U}_c.
\]
Taking $\theta \in (0,1)$ sufficiently close to $1,$ we find
\[
 \theta b_c - (1-\theta) \cdot \max_{\Sigma_c \setminus \mathcal{U}_c}  \lvert Y_0(K_c)\rvert        >0
\]
from which we obtain
\[
Y_\theta (K_c)  \geq  \theta b_c - (1-\theta) \cdot \max_{\Sigma_c \setminus \mathcal{U}_c}  \lvert Y_0(K_c)\rvert        >0  \quad \text{ on } \Sigma_c \setminus \mathcal{U}_c.
\]
Therefore, $Y_\theta$ is positively transverse to $\Sigma_c.$ Since $Y_\theta$ is $\sigma_{\rm el}$-invariant, it descends to a Liouville vector field positively transverse to the regularized bounded component.  This completes the proof of Theorem \ref{thm:main}.


\begin{remark}\label{remark:slightlyaboveL2}
Although the Liouville vector field $Y_\theta$ constructed above is sufficient for the subcritical case, namely, for energies below the second critical value, in   the supercritical case,  it would be convenient to have a transverse Liouville vector field which agrees precisely with $Y_1$ near the neck regions. 
Hence we refine the above construction by interpolating $Y_\theta$ near $\Sigma_{c,0}$ and $Y_1$ away from $\Sigma_{c,0}.$

 Let $x_{c_2} = x_*(c_2 )>0,$  so that the four saddle points of the regularized potential $V_c$ at $c=c_2$ are $(\pm x_{c_2}, 0)$ and $(\pm x_{c_2}, \pi).$  The two  critical points with $x_2=0$ correspond to $L_2$, while the two critical points with $x_2=\pi$ correspond to $L_3.$  Choose   pairwise disjoint  neighborhoods $\mathcal{N}^\pm_0$ and $\mathcal{N}^\pm_\pi$ in the $(x_1, x_2)$-plane of these saddle points, where  $\mathcal{N}_0^\pm$ correspond to the neck region around $L_2$ and $\mathcal{N}_\pi^\pm$ correspond to the neck region around $L_3.$ Define  
  $\mathcal{N} :=\mathcal N_0^+\cup \mathcal N_0^-\cup \mathcal N_\pi^+\cup \mathcal N_\pi^- .$ 
 We choose these neighborhoods so small that $pr^{-1}(\mathcal{N})$ is disjoint from $\Sigma_{c_2,0}, $  where $pr$ denotes the footpoint projection. 
 Denote by $\Sigma_{c_2}$ the compact connected component of   $K_{c_2}^{-1}(0) \setminus  pr^{-1}(\mathcal{N})$ that contains $\Sigma_{c_2,0}.$

The argument used in the proof of the above theorem still applies at $c=c_2$ and shows that 
 $Y_1(K_{c_2})\geq 0 $ on the part of $K_{c_2}^{-1}(0)$ lying over the bounded Hill's region, and the zero set consists of $\Sigma_{c_2,0}$ together with the   four critical points. Indeed, since $B_{c_2}(u_*(c_2))=0,$ and $u_*(c_2)$ is the unique maximum of $B_{c_2},$  we see that $u_*(c_2)=u_\pm(c_2)=u_{\max}(c_2)$ is the unique root of $B_{c_2}.$
 Since $P_{c_2}(u_*(c_2))=0,$ this shows that $\mathcal{Q}_{c_2}(u_{\max})=0. $  The other arguments hold for $c<c_3,$ so this proves non-negativity of   $Y_1(K_{c_2}). $ 
To find the zero set, suppose first that $x_1=0,$ so that $Y_1(K_{c_2})=p_2^2.$ Therefore, as before, $Y_1(K_{c_2})$ vanishes along $\Sigma_{c_2,0}.$ If $x_1 \neq0,$ then $Y_1(K_{c_2})$ can vanish only when $\mathcal{Q}_{c_2}=0,$ namely, at $u=u_*(c_2).$ Together with $B_{c_2}(u_*(c_2))=0,$ inequality \eqref{eq:Vintouv} implies that $v=0.$ Since $K_{c_2}=0,$ this shows $p_1=p_2=0,$    providing the four critical points.  Since the   critical points are not contained in $\Sigma_{c_2},$ we conclude that $Y_1(K_{c_2}) \geq 0$ on   $\Sigma_{c_2},$ and the zero set is precisely $\Sigma_{c_2,0}.$ We also note that $Y_0(K_{c_2})>0$ on $\Sigma_{c_2,0}.$

Let $s =  x_1^2 + p_2^2 $ so that $s=0$ on $\Sigma_{c_2,0}.$ For $0<s_-<s_+ $   sufficiently small,  we choose small $\sigma_{\rm el}$-invariant open  neighborhoods of $\Sigma_{c_2,0}$
\begin{equation}\label{eq:VandUI}
\mathcal{V}  = \{ s < s_- \} \subset \mathcal{U}  = \{ s < s_+ \} 
\end{equation}
 in the phase space $T^*(\R \times S^1).$ Choose $s_+$ small enough, so that $Y_0(K_{c_2})>0$ on $\Sigma_{c_2} \cap \overline{\mathcal{U} }$ and, moreover, $\overline{\mathcal{U}}\cap pr^{-1}(\mathcal{N}) = \emptyset.$
 Set
 \[
 \mathcal{A}_{c_2} := \Sigma_{c_2} \cap   \overline{\mathcal{U}  \setminus \mathcal{V} }   = \Sigma_{c_2} \cap \{  s_- \leq s \leq s_+\} ,
\]
which is   disjoint from $\Sigma_{c_2,0}$ and the four critical points. Consequently, by compactness we find a positive constant $k_*>0$ such that
\begin{equation}\label{eqlontheannul}
Y_0(K_{c_2}) \geq k_*, \quad Y_1(K_{c_2}) \geq k_* \quad \text{ on }  \mathcal{A}_{c_2}.
\end{equation}

Choose a smooth non-decreasing function $\varphi \colon [0,\infty)  \to [0,1]$ such that $\varphi(t)=0$ for $t \leq s_-,$ $\varphi(t)=1$ for $t \geq s_+$ and $0 \leq \varphi'(t) \leq 2/(s_+-s_-)$ for all $  t \geq 0.$   Define the $\sigma_{\rm el}$-invariant Liouville vector field 
\begin{equation}\label{eq:Liouvilevectorifle}
\hat{Y}:= Y_0 + X_{  \varphi_\theta f}, \quad  \; \varphi_\theta=\theta+(1-\theta)\varphi.
\end{equation} See Figure \ref{fig:zerointerpolation}.
Note that $\hat{Y}=Y_\theta$ on $\mathcal{V} $ and $\hat{Y}=Y_1$ outside $\mathcal{U} .$  In particular, since $\overline{\mathcal{U}}\cap pr^{-1}(\mathcal{N}) = \emptyset,$ we have $\hat{Y}= Y_1$ in a small neighborhood of the neck regions. 
 The previous discussion shows that $\hat{Y}(K_{c_2})>0$  on $\Sigma_{c_2} \setminus \mathcal{A}_{c_2}.$

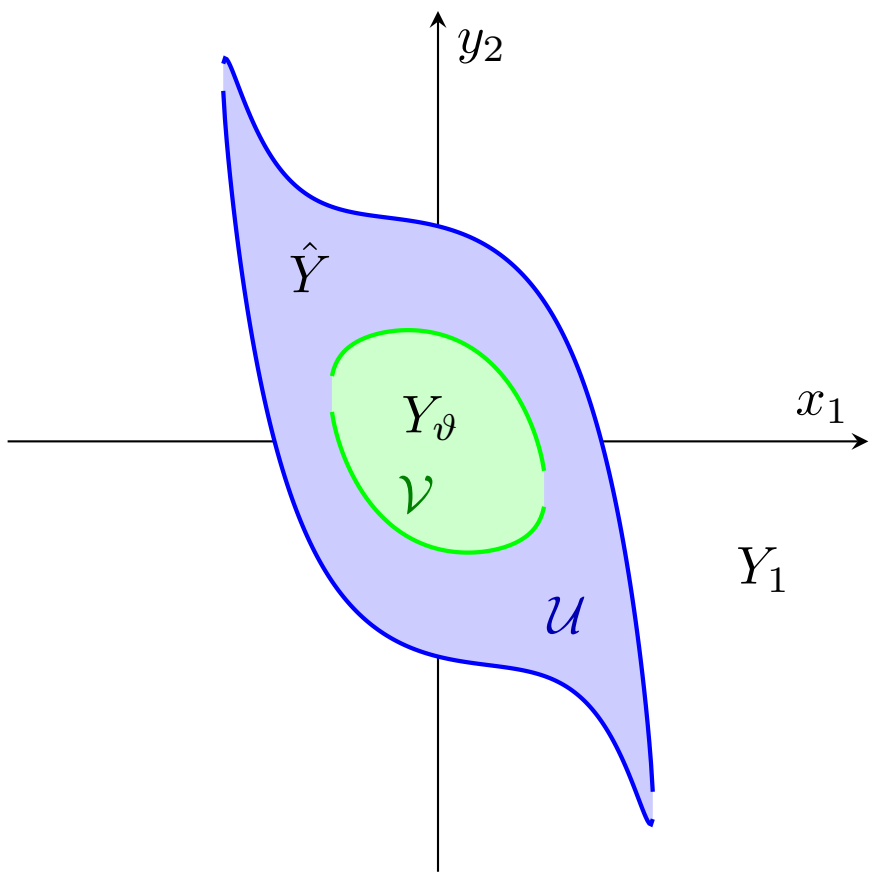
\begin{figure}
    \centering
\begin{tikzpicture}
\begin{axis}[
    axis lines=middle,
    xlabel={$x_1$},
    ylabel={$y_2$},
    xmin=-4, xmax=4,
    ymin=-4, ymax=4,
    samples=500,
    domain=-1:1,
    axis equal image,
    ticks=none,
]

\addplot[
    name path=upper,
    draw=none,
    domain=-4:4,
] {-(exp(2*x)-exp(-2*x))/16 + sqrt(4 - x^2)};

\addplot[
    name path=lower,
    draw=none,
    domain=-4:4,
] {-(exp(2*x)-exp(-2*x))/16 - sqrt(4 - x^2)};

\addplot[
    fill=blue!20,
    draw=none,
] fill between[
    of=upper and lower
];

\addplot[
    blue,
    thick,
    domain=-4:4,
] {-(exp(2*x)-exp(-2*x))/16 + sqrt(4 - x^2)};

\addplot[
    blue,
    thick,
    domain=-4:4,
] {-(exp(2*x)-exp(-2*x))/16 - sqrt(4 - x^2)};
\addplot[
    name path=upper,
    draw=none,
    domain=-4:4,
] {-(exp(2*x)-exp(-2*x))/16 + sqrt(1 - x^2)};

\addplot[
    name path=lower,
    draw=none,
    domain=-4:4,
] {-(exp(2*x)-exp(-2*x))/16 - sqrt(1 - x^2)};

\addplot[
    fill=green!20,
    draw=none,
] fill between[
    of=upper and lower
];

\addplot[
    green,
    thick,
    domain=-4:4,
] {-(exp(2*x)-exp(-2*x))/16 + sqrt(1 - x^2)};

\addplot[
    green,
    thick,
    domain=-4:4,
] {-(exp(2*x)-exp(-2*x))/16 - sqrt(1 - x^2)};
\end{axis}
\draw[blue!70!black] (3.7,1.7) node {$\mathcal{U} $};
\draw[green!50!black] (2.7,2.5) node {$\mathcal{V} $};
\draw (2,4) node {$\hat{Y}$};
\draw (2.8,3) node {$Y_\theta$};
\draw (5,2) node {$Y_1$};
\end{tikzpicture}
\caption{An  illustration of the  transition regions $\mathcal{U} $ and $\mathcal{V} $ in the $(x_1, y_2)$-plane }
\label{fig:zerointerpolation}
\end{figure}

   To show $\hat{Y}(K_{c_2})>0$ on $\Sigma_{c_2},$ it suffices to consider the annulus $\mathcal{A}_{c_2}.$ 
   We compute
\begin{align*}
\hat{Y}(K_{c_2}) &= (1-\varphi_\theta)Y_0(K_{c_2}) + \varphi_\theta Y_1(K_{c_2}) + (1-\theta)f X_{\varphi}(K_{c_2}) \\& \geq k_* + (1-\theta)f X_{\varphi}(K_{c_2}). 
\end{align*}
See \eqref{eqlontheannul}.
It remains to estimate the last term in the right-hand side above. 
Since $\mathcal{A}_{c_2}$ is compact, we find $C_*:=\max_{\mathcal{A}_{c_2}} \lvert f X_{\varphi}(K_{c_2}) \rvert <\infty.$
We now choose $\theta<1$ sufficiently close to $1,$ so that $(1-\theta)C_* \leq  {k_*}/{2}.$  
It follows that
\[
\hat{Y}(K_{c_2}) \geq k_*  - (1-\theta)C_* \geq \frac{1}{2}k_*>0 \quad \text{ on }  \mathcal{A}_{c_2}.
\]
Consequently, the Liouville vector field $\hat{Y}$ is positively transverse to $\Sigma_{c_2}.$

We now continue the discussion above for energies slightly above the second critical value. Choose $\varepsilon_1 \in (0, \varepsilon_0)$  small enough, where $\varepsilon_0$ is as in Lemma \ref{lem:orderofvalues},  so that   for every  $c \in (c_2, c_2+\varepsilon_1)$ the complement  $K_c^{-1}(0)\setminus pr^{-1}(\mathcal{N})$ contains  the compact connected component, denoted by 
$\hat{\Sigma}_c, $ converging to  $\Sigma_{c_2}$ as $c \to c_2^+.$ 

We claim that $\hat{Y}(K_c)>0$ on $\hat{\Sigma}_c$ for every $c \in (c_2, c_2+\varepsilon_1),$ after shrinking $\varepsilon_1>0$ if necessary.  Assume by contradiction that we find sequences $c_n>c_2$ such that $c_n \to c_2$ as $n \to \infty$ and $z_n \in \hat{\Sigma}_{c_n}$ such that $\hat{Y}(K_{c_n})(z_n) \leq 0.$ After passing to a subsequence we may assume that $z_n \to z \in \Sigma_{c_2}.$ We then have
\[
\hat{Y}(K_{c_2})(z)  = \lim_{n \to \infty} \hat{Y}(K_{c_n})(z_n) \leq0,
\]
a contradiction to the fact that $\hat{Y}(K_{c_2})>0$ on $\Sigma_{c_2}.$  Consequently, by decreasing $\varepsilon_1>0$ further if necessary, we find 
\[
\hat{Y}(K_c)>0 \quad \text{ on } \hat{\Sigma}_c 
\]
for every $c \in (c_2, c_2+\varepsilon_1).$

Fix $d \in (c_1, c_2).$ For $s_+>0$ sufficiently small, $Y_0(K_d)>0$ on $\Sigma_d \cap \mathcal{V},$ so $\hat{Y}(K_d)= Y_\theta(K_d)>0$ there. Since $Y_1(K_d)>0$ on the compact set $\Sigma_d \setminus \mathcal{V},$ the Liouville vector field $\hat{Y}$ with the same parameter is positively transverse to $\Sigma_d,$ provided $\theta<1$ is sufficiently close to $1.$

 \end{remark}



\section{Slightly above the second critical value}\label{sec:4}

We   consider energies  slightly above the second critical value $c_2=H(L_2)$, where the critical values are still ordered as in Lemma \ref{lem:orderofvalues}.
We shall modify the regularized Hamiltonian $K_c$ outside a compact set containing the neck regions 
in such a way that   the resulting perturbed energy level is   compact  and   diffeomorphic to $S^1 \times S^2.$  We combine the techniques from \cite{AFvKP12} and \cite{dPKSS24}.

\subsection{In the neck regions}

Let $c >c_2$ be sufficiently close to $c_2$ as in Remark \ref{remark:slightlyaboveL2} so that the Liouville vector field $\hat{Y}$ is positively transverse to $\hat{\Sigma}_c.$

Put
\[
z_* = z_*(c) = ( x_*(c), 0,0,-A_1(x_*(c)))  
\]
which is a critical point of $K_c.$ 
 For convenience, we introduce  ($c$-dependent) local coordinates $\zeta = (\xi_1, \xi_2, \eta_1, \eta_2)$ centered at $z_*$, defined by 
\[
\xi_1 = x_1 - x_*(c), \quad \xi_2 = x_2, \quad \eta_1 = y_1, \quad \eta_2 = y_2 + A_1(x_*(c)).
\]
Since $x_2$ is $2\pi$-periodic, the coordinate $\xi_2$ is understood only locally   near $x_2=0.$  For $r>0,$ denote by 
\[
\mathcal{B}_r(z_*(c)):= \{  \lvert \zeta \rvert <r\}
\]
the coordinate $r$-ball.

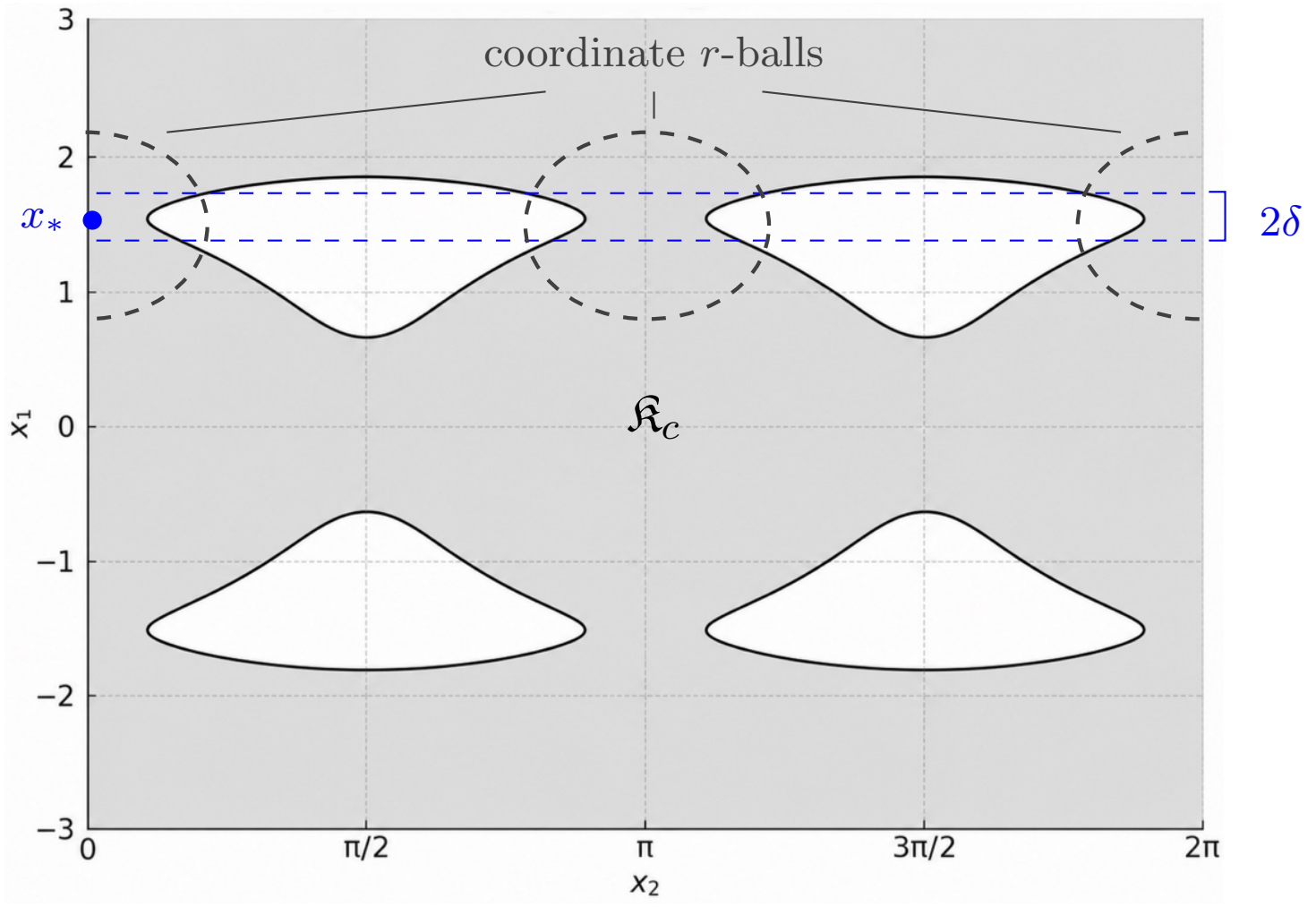
\begin{figure}[ht]
\centering
\begin{tikzpicture} 
  \node[anchor=south west,inner sep=0] at (0,0) {\includegraphics[width=9cm]{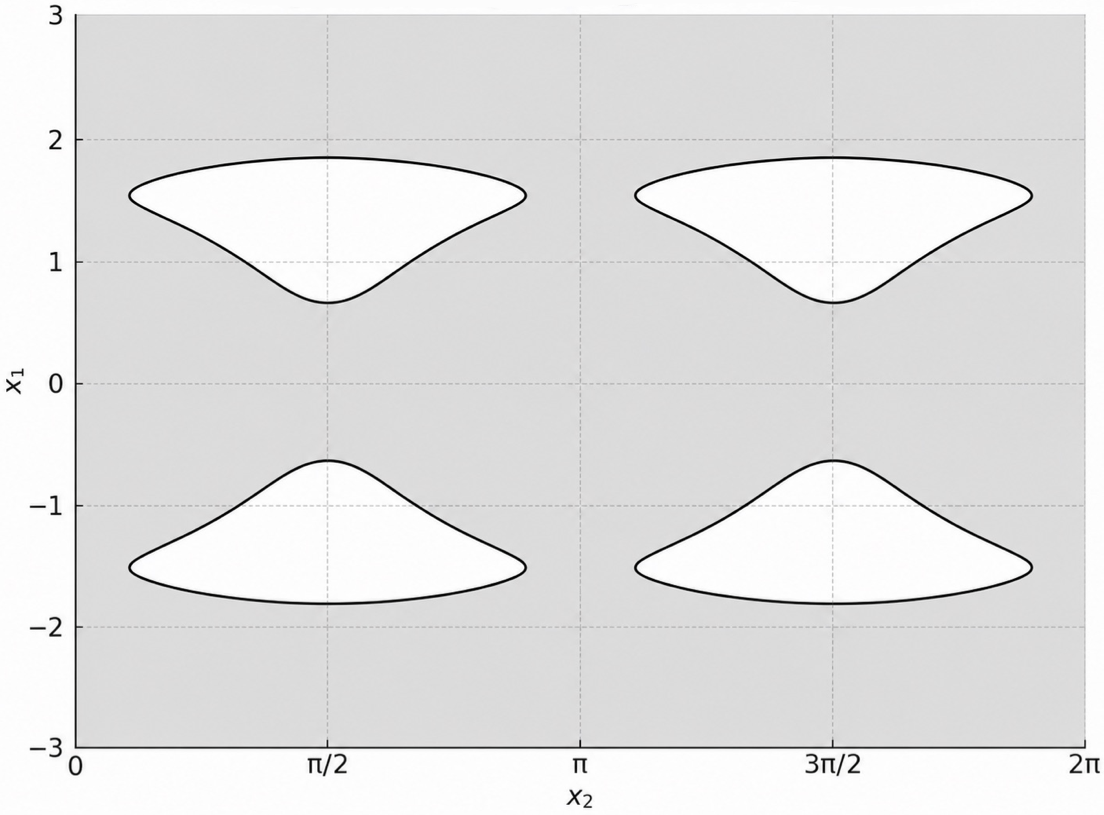}};

  \draw[blue,dashed] (0.68,5.25) -- (8.83,5.25)
  (0.68,4.9 ) -- (8.83,4.9 );
  \node[left,blue] at (0.58,5.05) {\small $x_* $};
\fill[blue] (0.65,5.05) circle (2 pt);
  \begin{scope}[xshift=-1.1cm, yshift=-0.6cm]
  \draw[blue]
  (10,5.50) -- (10.12,5.50) -- (10.12,5.86) -- (10,5.86);
  \end{scope}
  \node[right,blue] at (9.15,5.05) {\small $2\delta$};
\begin{scope}[xshift=0.1cm]
   \draw[darkgray,dashed,thick]    (0.5,5.70) arc (90:-90:0.90 and 0.69);
   \end{scope}
   \begin{scope}[yshift=-0.2cm]
  \draw[darkgray,dashed,thick]    (4.75,5.21) ellipse (0.90 and 0.69);
  \end{scope}
  \begin{scope}[xshift=0.3cm]
  \draw[darkgray,dashed,thick]   (8.53,4.32) arc (-90:-270:0.90 and 0.69);
  \end{scope}
  
  \node[darkgray] at (4.8,6.30) {\small coordinate $r$-balls};
  \draw[darkgray] (4,6 ) -- (1.20,5.7);
  \draw[darkgray] (4.8, 6) -- (4.8,5.8);
    \draw[darkgray] (5.6, 6) -- (8.25,5.7);
  \node at (4.8,3.6) {\large $\mathfrak{K}_c$};
\end{tikzpicture}
\caption{An illustration of Lemma \ref{lem:locallemman} in the position plane}
\label{Fig:Lemma41}
 \end{figure}

\begin{lemma}\label{lem:locallemman}
    Fix $r>0$ sufficiently small. After decreasing $\varepsilon_1>0$ as in Remark \ref{remark:slightlyaboveL2} if necessary, we find $\delta>0$ such that for every $c \in (c_2, c_2+\varepsilon_1),$ every $(x,y) \in K_c^{-1}(0)$ with $\lvert x_1 - x_*(c)\rvert <\delta$ lies either in $\mathcal{B}_r(z_*(c))$ or in the corresponding coordinate $r$-ball centered at $(x_*(c), \pi,0,-A_1(x_*(c))).$  See Figure \ref{Fig:Lemma41}.
\end{lemma}
\begin{proof}
    Assume by contradiction that such $\varepsilon_1>0$ and $\delta>0$ do not exist. Then  we find sequences $c_n \in (c_2,c_2+1/n)$ and $z_n=(x_{1,n}, x_{2,n}, y_{1,n}, y_{2,n}) \in K_{c_n}^{-1}(0),$ lying outside both coordinate balls, satisfying  $\lvert x_{1,n} - x_*(c_n) \rvert<1/n.$    We see, in particular, that $c_n \to c_2$ and $ x_{1,n} \to x_*(c_2). $
Recall from Remark \ref{remark:valuev1xstar} that $B_*(c_2)=0,$ or equivalently, $V_{c_2,1}(x_*(c_2))=-c_2/4.$ We therefore have
\begin{align*}
V_{c_2}(x_*(c_2),x_2) &= V_{c_2,1}(x_*(c_2))+V_{c_2,2}(x_2)\\
&= -\frac{1}{4} \left( c_2 + \frac{1}{8} \cos^2 x_2\right) \sin^2 x_2 \\
&\geq 0,
\end{align*}
where the inequality follows from $c_2 <c_3=-11/8.$ Note that the equality holds if and only if $x_2=0$ or $x_2=\pi.$ Since $V_{c_n}( x_{1,n}, x_{2,n}) \leq 0$ for all $n,$ it follows that, after passing to a subsequence, $x_{2,n} \to 0$ or $x_{2,n} \to \pi,$ and consequently $V_{c_n} (x_{1,n}, x_{2,n}) \to 0.$ Since $p_{1,n}^2 +p_{2,n}^2 = -2V_{c_n}(x_{1,n}, x_{2,n})$ on $K_{c_n}^{-1}(0),$ where  $p_{1,n}=y_{1,n}+A_2(x_{2,n})$ and $p_{2,n}=y_{2,n}+A_1(x_{1,n}),$  we also have $p_{1,n} \to 0$ and $p_{2,n} \to 0.$ 

We first consider the case $x_{2,n}\to 0,$ so that $A_2(x_{2,n}) \to 0. $ Hence, $\eta_{1,n} = y_{1,n} \to 0.$ Since $A_1(x_{1,n}) - A_1(x_*(c_n)) \to 0,$ we also have $\eta_{2,n} = y_{2,n} +A_1(x_*(c_n))\to 0.$ Thus, in the local coordinates described above, 
\[
(x_{1,n}-x_*(c_n), x_{2,n}, \eta_{1,n}, \eta_{2,n}) \to 0,
\]
implying that   $z_n \in \mathcal{B}_r(z_*(c_n)) $ for all $n$ large enough, a contradiction. The case $x_{2,n} \to \pi$ can be treated in a similar way.
This completes the proof.
\end{proof}

It follows that it suffices to carry out a local construction in $\mathcal{B}_r(z_*(c))$ for $r>0$ sufficiently small. The corresponding constructions near the other critical points are  obtained  by symmetry.

Following \cite[Section 7]{AFvKP12}, we consider  the Weinstein-like vector field 
\[
W_{a,b} = a \xi_1 \partial_{\xi_1 }+b \xi_2 \partial_{\xi_2} + (1-a) \eta_1 \partial_{\eta_1} + (1-b)\eta_2 \partial_{\eta_2}.
\]
 Given $a, b \in \R,$ let 
\[
g_c(\xi) = g_{a,b,c}(\xi) = \frac{1}{4} \delta_c ( a+1-b) \xi_1 \xi_2, 
\]
where $\delta_c := \cosh(2x_*(c))= 2u_*(c)^2-1>1$
(recall  that   $u_*=\cosh x_*>1 $ is the unique root of   $P_c $ in $(1,\infty),$ see   Section \ref{sec:Hill}),
and define the Liouville vector field 
\[
X_c= X_{a,b,c}   :=W_{a,b} + X_{g_c }  .
\]
Explicitly, we have
\begin{equation}\label{eq:defofxfg}
X_c =  a \xi_1 \partial_{\xi_1 }+b \xi_2 \partial_{\xi_2} + \left( (1-a) \eta_1  +\sigma_c \xi_2  \right)\partial_{\eta_1} + \left( (1-b)\eta_2   +\sigma_c \xi_1  \right)\partial_{\eta_2} ,
\end{equation}
where $\sigma_c:= -\frac{1}{4} \delta_c ( a+1-b). $

\begin{lemma}\label{lem:Weinstein} 
After decreasing $\varepsilon_1>0$ if necessary, we find constants 
\[
a<0, \quad 0<b<1, \quad \mathfrak{w}>0, \quad r_0>0,
\]
independent of $c,$ such that for every $c \in (c_2, c_2+\varepsilon_1),$ we have
\[
X_c(K_c)(z_*(c)+\zeta) \geq \mathfrak{w}\lvert \zeta \rvert^2  \quad \text{ on }    \mathcal{B}_{r_0}(z_*(c)).
\]
Thus, the Liouville vector field $X_c$ is positively transverse  to $K_c^{-1}(0)$ near $z_*(c).$

\end{lemma}
\begin{proof}
In a sufficiently small neighborhood of  $z_*(c)$ we have the Taylor expansion 
\begin{equation}\label{eq:Taylorofkc}
\begin{aligned}
K_c( z_*+\zeta) &= \frac{1}{2}\left( ( \eta_1 + A_2'(0)\xi_2)^2 + ( \eta_2 + A_1'(x_*)\xi_1)^2 \right) + V_{c,1}(x_*)\\
& \quad +\frac{1}{2}V_{c,1}''(x_*)\xi_1^2 + V_{c,2}(0)+\frac{1}{2}V_{c,2}''(0) \xi_2^2 +\mathcal{R}_c(\zeta) \\
&= K_c(z_*)+  Q_c(\zeta)+\mathcal{R}_c(\zeta) ,
\end{aligned}
\end{equation}
where the quadratic term $Q_c$ is given by
\begin{equation}\label{eq:qudarictopart}
    Q_c(\zeta)  =  \frac{1}{2}\left( ( \eta_1 + \frac{1}{4}\xi_2)^2 + ( \eta_2 +  \frac{1}{4}\delta_c\xi_1)^2  +   \lambda_1(c) \xi_1^2  + \lambda_2(c) \xi_2^2    \right)
\end{equation}
with  $\lambda_1(c) = V_{c,1}''(x_*)<0,$   $\lambda_2(c) = V_{c,2}''(0)>0,$ and 
$\mathcal{R}_c(\zeta) = O(\lvert \zeta \rvert ^3)$ is the remainder term of order at least three. We write $Q_c = Q_{c,1} + Q_{c,2},$ where
 \begin{align*}
     Q_{c,1} (\xi_1, \eta_2)&=   \frac{1}{2} \left(  ( \eta_2 +  \frac{1}{4}\delta_c\xi_1)^2  +   \lambda_1(c) \xi_1^2      \right) \\
      Q_{c,2} (\xi_2, \eta_1) &=  \frac{1}{2}\left( ( \eta_1 + \frac{1}{4}\xi_2)^2   + \lambda_2(c) \xi_2^2    \right)
 \end{align*}
A direct computation shows that
\begin{align*}
X_c(Q_{c,1}) &= \left( a \lambda_1(c) - \frac{\delta_c^2}{16}(1-b) \right) \xi_1^2 + (1-b)\eta_2^2 \\
X_c(Q_{c,2}) &= \left(  b \left( \lambda_2(c) + \frac{1}{16}  \right) + \frac{1}{4} \sigma_c \right) \xi_2^2 + (1-a) \eta_1^2 + \left( \frac{1}{4} ( b+1-a) + \sigma_c \right) \xi_2  \eta_1
\end{align*}
Take  $\tau \in (0,1)$ and set 
\[
a =  \frac{\delta_{c_2}^2}{8\lambda_1(c_2)}\tau <0, \qquad b = 1-\tau>0,
\]
which are independent of $c.$
We compute
\[
a\lambda_1(c_2)-\frac{\delta_{c_2}^2}{16}(1-b) = \frac{\delta_{c_2}^2}{16}\tau>0,
\]
showing that   $X_{c_2}(Q_{c_2,1})$ is positive-definite for every $\tau \in (0,1).$  On the other hand, we take the limit $\tau \to 0^+,$ so that  $a \to 0, b \to 1$ and $\sigma_c \to 0.$ In the limit, we obtain the positive-definite quadratic form 
\[
X_{c_2}(Q_{c_2,2}) =  \lambda_2(c_2) \xi_2^2 + \left( \eta_1+\frac{1}{4}\xi_2\right)^2 = 2Q_{c_2,2} .
\]
Fix $\tau \in (0,1)$  so small   that both $X_{c_2}(Q_{c_2,1})$ and $X_{c_2}(Q_{c_2,2})$ are positive-definite. By continuity, after decreasing $\varepsilon_1>0$ if necessary,   we find $\mathfrak{w}>0$ such that
\[
X_c(Q_c)(\zeta) \geq 2\mathfrak{w}\lvert \zeta \rvert^2, \qquad c \in [c_2, c_2+\varepsilon_1].
\]

The vector field $X_c$ is linear in $\zeta,$ and its coefficients  are uniformly bounded for $c \in [c_2, c_2+\varepsilon_1].$   Hence,  we find $C_0>0,$ independent of $c,$  such that  
\[
\lvert X_c(\mathcal{R}_c)\rvert \leq C_0  \lvert \zeta \rvert^3.
\]
We now choose $r_0>0,$ independent of $c,$ such that $C_0   r_0 \leq \mathfrak{w},$ and hence  
\[
\lvert X_c(\mathcal{R}_c)\rvert \leq \mathfrak{w} \lvert \zeta \rvert^2   \quad \text{ on } \mathcal{B}_{r_0}(z_*(c)),
\]
uniformly for $c \in (c_2, c_2+\varepsilon_1).$ 
Consequently,
\[
X_c(K_c) = X_c(Q_c)+X_c(\mathcal{R}_c) \geq \mathfrak{w}\lvert \zeta \rvert^2 \quad \text{ on }  \mathcal{B}_{r_0}(z_*(c))  .
\]

We define $k(c):= K_c( z_*(c))=B_*(c).$ The proof of Lemma \ref{lem:orderofvalues} shows that it is a smooth function of $c,$ provided $c<c_3 .$  We differentiate
\[
k'(c) = \partial_c K_c ( z_*(c)) = - \frac{1}{4}( \cosh^2 x_*(c) -1)<0,
\]
where in the first identity we have used the fact that $z_*$ is a critical point.  Since $k( c_2)=0$ (cf. Remark \ref{remark:valuev1xstar}), it follows that $K_c( z_*(c))  <0$ 
for all $c >c_2$ sufficiently close to $c_2. $ Hence, $z_*(c) \notin K^{-1}_c(0).$ Therefore, every point in $K^{-1}_c(0)\cap \mathcal{B}_{r_0}(z_*(c))$ satisfies $\zeta \neq 0,$ and hence the preceding estimates complete the proof. 
\end{proof}

\begin{remark} One may wonder whether   the  Weinstein-like vector field $W_{a,b}$ is already sufficient. However, a direct computation shows that $W_{a,b}(Q_{c,1})$ is not positive-definite  for any choice of $a$ and $b.$   
\end{remark}

In the following lemma we express $\lambda_1$ and $\lambda_2$ given in \eqref{eq:qudarictopart} in terms of $u_* = \cosh x_*.$ Recall that $u_*$ is the unique zero of the polynomial $P_c$ in $(1,\infty),$ see \eqref{eq;polynomialP}.

\begin{lemma} \label{lem:lambda12u}For every $c < c_3,$ we have
\[
\lambda_1(c) = \frac{ (1-u_*(c)^2)(u_*(c)^3-2)}{4u_*(c)}, \qquad \lambda_2 (c)= \frac{ u_*(c)^3 - u_*(c)+4}{8u_*(c)}.
\]
 \end{lemma}
\begin{proof} Dividing $P_c(u_*) =  2u_*^3 + (8c-1)u_*+8=0$ by $ u_* >0$ gives 
\[
8c  = 1-2u_*^2 - \frac{8}{u_*}.
\]
Plugging this into   \eqref{eq:derofV1} and \eqref{eq:derivofV2} and using  $\cosh 2x_* = 2u_*^2-1$ and $\cosh 4x_* =  8u_*^4 - 8u_*^2+1$ provides the desired identities. 
\end{proof}


 \subsection{First interpolation}\label{sec:firstin}

Let the  constants $\varepsilon_1, a,b,\mathfrak{w}$ and $r_0$ be as in Lemma \ref{lem:Weinstein}, so that the perturbed Weinstein-like vector field $X_c$ is positively transverse to $K_c^{-1}(0)$ in $\mathcal{B}_{r_0}(z_*(c))$ for every $c \in (c_2, c_2+\varepsilon_1).$
We now interpolate the two Liouville vector fields $\hat{Y}$ and $X_c$ near  $z_*(c) .$   As before, we work with the $c$-dependent local coordinates $\zeta = (\xi, \eta)$ centered at $z_*(c).$

 We first consider the following linear function, separating the bounded part from the unbounded part of the neck region:
\begin{equation}\label{eq;alphad}
   \alpha_{ c}(\zeta)  = (\mu_{ c} -\frac{1}{4}) \xi_1 - \eta_2, \qquad \mu_{ c} =  \frac{ 4 \lambda_2 (c)    }{\delta_c - 1 }>0.
\end{equation}
Recall that $\delta_c = \cosh ( 2x_*) = 2u_*(c)^2-1 >1$ and $\lambda_2(c) = V_{c,2}''(0) >0.$

\begin{remark}[{cf. \cite[Remark 7.2]{AFvKP12}}] \label{rmk:separation} We explain the role of  the function $\alpha_c.$ At $c=c_2,$ the local hypersurface $\alpha_{c_2}^{-1}(0)$ separates the   two local branches of the neck region. More precisely, for $r>0$ sufficiently small,  
\begin{equation}\label{eq:atcritialenergy}
K^{-1}_{c_2}(0) \cap   \alpha_{c_2}^{-1}(0)  \cap  \mathcal{B}_r( z_{c_2})=  \{ z_{c_2}\},
\end{equation}
where $z_{c_2}=z_*(c_2) .$
Moreover, the local branch of $K_{c_2}^{-1}(0),$ obtained as the limit of the bounded components for $c<c_2,$  lies in $\{ \alpha_{c_2} <0\}$ away from $z_{c_2}.$

To see this, we first suppose that   $Q_{c_2}(\zeta)=0.$ Since $Q_{c_2,2}$ is positive-definite,     $Q_{c_2,1}  \leq0,$ and hence 
\begin{equation}\label{eq:xieat2remark}
\lvert \eta_2 + \frac{1}{4}\delta_{c_2} \xi_1 \rvert \leq \sqrt{-\lambda_1(c_2)} \lvert \xi_1 \rvert.
\end{equation}
 Therefore, we find $-1 \leq \gamma \leq 1,$  depending on $\xi_1$ and  $\eta_2$, such that $\eta_2 + \frac{1}{4}\delta_{c_2} \xi_1 = \gamma \sqrt{-\lambda_1(c_2)}\xi_1,$ and hence $\eta_2 = ( - \frac{1}{4}\delta_{c_2} + \gamma \sqrt{-\lambda_1(c_2)} )\xi_1.$ We then obtain
 \[
 \alpha_{c_2}(\zeta) = \left(  \mu_{ c_2} - \frac{1}{4}  + \frac{1}{4}\delta_{c_2} - \gamma \sqrt{-\lambda_1(c_2)} \right)\xi_1 .
 \]
In view of Lemma \ref{lem:lambda12u} we find that
\[
 ( \mu_{c_2} - \frac{1}{4}  + \frac{1}{4}\delta_{c_2} )^2 +\lambda_1(c_2) = \frac{ 8 u_{c_2}^7 + u_{c_2}^6 - 8u_{c_2}^5 - 2u_{c_2}^4 + u_{c_2}^2 + 16}{ 16u_{c_2}^2 ( u_{c_2}^2 -1)^2},  
\]
where $u_{c_2}=u_*(c_2).$
Using \eqref{eq:ustaridentty}, the numerator reduces to $16 ( 4u_{c_2}^4 + u_{c_2}^3 + 8u_{c_2}^2 +5),$ which is positive. 
This implies that
\[
 \mu_{ c_2} - \frac{1}{4}  + \frac{1}{4}\delta_{c_2} >  \sqrt{-\lambda_1(c_2)}  ,
\] 
and hence 
\begin{equation}\label{eq:signalphaxi}
{\rm sign}(\alpha_{c_2}(\zeta)) = {\rm sign}(\xi_1)  \quad \text{ on }  \{ Q_{c_2}(\zeta)=0\} \setminus \{0\}. 
\end{equation}  In particular, $\{ Q_{c_2}=0\} \cap \{\alpha_{c_2}=0\} = \{ 0 \}.$ 
We now assume by contradiction that  \eqref{eq:atcritialenergy} fails, so  there is a sequence   $\zeta_n  \neq 0 $ converging to $0$ as $n \to \infty,$ satisfying $z_{c_2} + \zeta_n \in K_{c_2}^{-1}(0)$ and $\alpha_{c_2}(\zeta_n)=0.$
 Since $K_{c_2}(z_{c_2})=0,$ the Taylor expansion \eqref{eq:Taylorofkc} implies that $Q_{c_2}(\zeta_n) + \mathcal{R}_{c_2}(\zeta_n)=0.$
We set $\bar{\zeta}_n:= \zeta_n/ \lvert \zeta_n \rvert,$ so that $\lvert \bar{\zeta}_n \rvert =1$ for all $n.$ After passing to a subsequence, we may assume that $\bar{\zeta}_n \to \bar{\zeta}$ as $n \to \infty.$ We then have $\alpha_{c_2}(\bar{\zeta})=0.$ Moreover, since $Q_{c_2}$ is a quadratic form, we have 
\[
Q_{c_2}(\bar{\zeta}_n) + \frac{\mathcal{R}_{c_2}(\zeta_n)}{\lvert \zeta_n\rvert^2}=0.
\]
Passing to the limit $n \to \infty,$ we obtain $Q_{c_2}(\bar{\zeta})=0$ because $\mathcal{R}_{c_2}$ is of order at least three. The above argument then shows that $\bar{\zeta}=0,$   contradicting   $\lvert \bar{\zeta}\rvert=1.$ This establishes \eqref{eq:atcritialenergy}. Moreover, for $c<c_2,$ the bounded component approaches the $L_2$-neck region from the side $x_1< x_*(c).$   In the local coordinates this corresponds to the branch with $\xi_1<0.$ Thus, by \eqref{eq:signalphaxi}, the limiting bounded branch has negative $\alpha_{c_2}$-sign near $z_{c_2}.$ Since $\alpha_{c_2}$ has a constant sign on each connected local branch,   the limiting local branch of $K_{c_2}^{-1}(0),$ with $z_{c_2}$ removed, is contained in  $\{\alpha_{c_2}<0\}.$
\end{remark}

Let $\varepsilon_1>0$ and $r_0>0$ be as in Lemma  \ref{lem:Weinstein} and set
\[
 {\varepsilon}_{\rm Tay}:= \min \left\{ \varepsilon_1, \frac{c_3-c_2}{100} \right\}>0.
\]
After decreasing $r_0>0$ if necessary, we assume that \eqref{eq:atcritialenergy} holds for $r=r_0 $ and that the four coordinate $r_0$-balls corresponding to the four critical points are pairwise disjoint.
Since the Taylor expansion smoothly depends on $c,$ we find $r_{\rm Tay}>0$ and $C_{\rm Tay}>0$ such that the higher order term $\mathcal{R}_c$ satisfies
\begin{equation}\label{eq:CTay}
\lvert \mathcal{R}_c(\zeta)\rvert \leq C_{\rm Tay} \lvert \zeta \rvert^3, \qquad \lvert d \mathcal{R}_c(\zeta)\rvert \leq C_{\rm Tay} \lvert \zeta \rvert^2 
\end{equation}
for every $c \in [c_2, c_2+ {\varepsilon}_{\rm Tay}]$ and every $\lvert \zeta \rvert \leq r_{\rm Tay}.$ Recall from Remark \ref{remark:slightlyaboveL2} that $\hat{Y} = Y_1$ outside $\mathcal{U} = \{ s< s_+\},$ where $s = x_1^2 +p_2^2.$ Since $z_*(c_2) \notin \overline{\mathcal{U}} \cup \Sigma_{c_2,0},$ we find $\hat{\varepsilon}>0$ and $\hat{r}>0$ such that 
\[
\overline{\mathcal{B}}_{\hat{r}}(z_*(c_2)) \cap \Sigma_{c_2,0} = \emptyset, \qquad \overline{\mathcal{B}}_{\hat{r}}(z_*(c)) \subset \{ s > s_+\}, \;\;\forall c \in [c_2, c_2+\hat{\varepsilon}].
\]
We then set
\[
\varepsilon_*:= \frac{1}{2}\min \{ \varepsilon_{\rm Tay}, \hat{\varepsilon}\}>0, \qquad r_*:= \frac{1}{2} \min \left\{ r_{\rm Tay}, r_0, \hat{r} , \frac{\pi}{4}\right\}>0 
\]
and $I_*:=[c_2, c_2+\varepsilon_*].$ It follows that 
\[
\hat{Y}=Y_1  \quad \text{ on } \mathcal{B}_{r_*}(z_*(c)), \quad c \in I_*.
\]
For $c \in I_* \setminus \{ c_2\},$     the Liouville vector field $X_c$ is positively transverse to $K_c^{-1}(0)$ in $\mathcal{B}_{r_*}(z_*(c)).$  At $c=c_2,$   the limiting local branch of $K_{c_2}^{-1}(0),$ coming from   the bounded components for $c<c_2,$ is contained in  $\{\alpha_{c_2}<0\}$ away from $z_*(c_2).$  The corresponding statements near $x_2=\pi$ and near the other two critical points with $x_1 = -x_*(c)$ follow from the symmetries.

 For a sufficiently small parameter $\rho>0,$ take  a smooth non-decreasing function  $\chi \colon \R \to [0,1]$   satisfying $\chi(s) = 0$ for all $s \leq -2\rho,$  $\chi(s)=1$ for all $s \geq -\rho$ and  $0 \leq \chi'(s) \leq 2 /\rho $ for all $s.$ Define
\[
\ell_c(\zeta):= \chi (\alpha_c(\zeta)).
\]
Note that $\ell_c(\zeta)=0$ for $\alpha_c(\zeta) \leq -2\rho$ and $\ell_c(\zeta)=1$ for $\alpha_c(\zeta) \geq -   \rho.$  The transition region is given by
\begin{equation}\label{eq:transitionregion}
\mathcal{S}_{c,\rho } := K_c^{-1}(0) \cap \mathcal{B}_{r_*}(z_*(c)) \cap \{ - 2 \rho \leq \alpha_c \leq - \rho \}.
\end{equation}

\begin{figure}
    \centering
    \begin{tikzpicture}
        \draw (5,1.5) to[out=225,in=350,looseness=0.5] (2,1);
        \draw (5,-1.5) to[out=135,in=10,looseness=0.5] (2,-1);
        \draw[blue] (2,1) to[out=170,in=320] (1,1.5);
        \draw[blue] (2,-1) to[out=190,in=40] (1,-1.5);
        \draw (1,1.5) to[out=140,in=315] (0,2.5);
        \draw (1,-1.5) to[out=220,in=45] (0,-2.5);
        \draw[blue] (2,1) -- (2,-1);
        \draw[blue] (1,1.5) -- (1,-1.5);
        \fill[blue!15] (2,1) to[out=170,in=320] (1,1.5) -- (1,-1.5) to[out=40,in=190] (2,-1) -- (2,1);
        \fill[green!15] (5,1.5) to[out=225,in=350,looseness=0.5] (2,1) -- (2,-1) to[out=10,in=135,looseness=0.5] (5,-1.5) -- (5,1.5);
        \fill[lightgray!20] (1,1.5) -- (1,-1.5) to[out=220,in=45] (0,-2.5) -- (0,2.5) to[out=315,in=140] (1,1.5);
        \draw[green!50!black] (3,-1.5) node {$\alpha_c>-\rho$};
        \draw[blue!80!black] (2.5,1.8) node {$-2\rho<\alpha_c<-\rho$};
        \draw[gray!80!black] (0.5,-3) node {$\alpha_c<-2\rho$};
        \clip (2.5,0) circle (1.3cm);
        \fill[orange!15] (1,-1.5) to[out=40,in=190] (2,-1) -- (2,1) to[out=170,in=320] (1,1.5) -- (1,-1.5) ;
        \draw[blue] (2,1) to[out=170,in=320] (1,1.5);
        \draw[blue] (2,-1) to[out=190,in=40] (1,-1.5);
        \draw[blue] (2,1) -- (2,-1);
        \draw[blue] (1,1.5) -- (1,-1.5);
        \draw[orange!70!black] (1.6,0) node {$\mathcal{S}_{c,\rho}$};
    \end{tikzpicture}
    \caption{An illustration of the regions of  the energy hypersurface determined by $\alpha_c$ and   the transition region $\mathcal{S}_{c,\rho}$}
    \label{fig:firstinterpolation}
\end{figure}
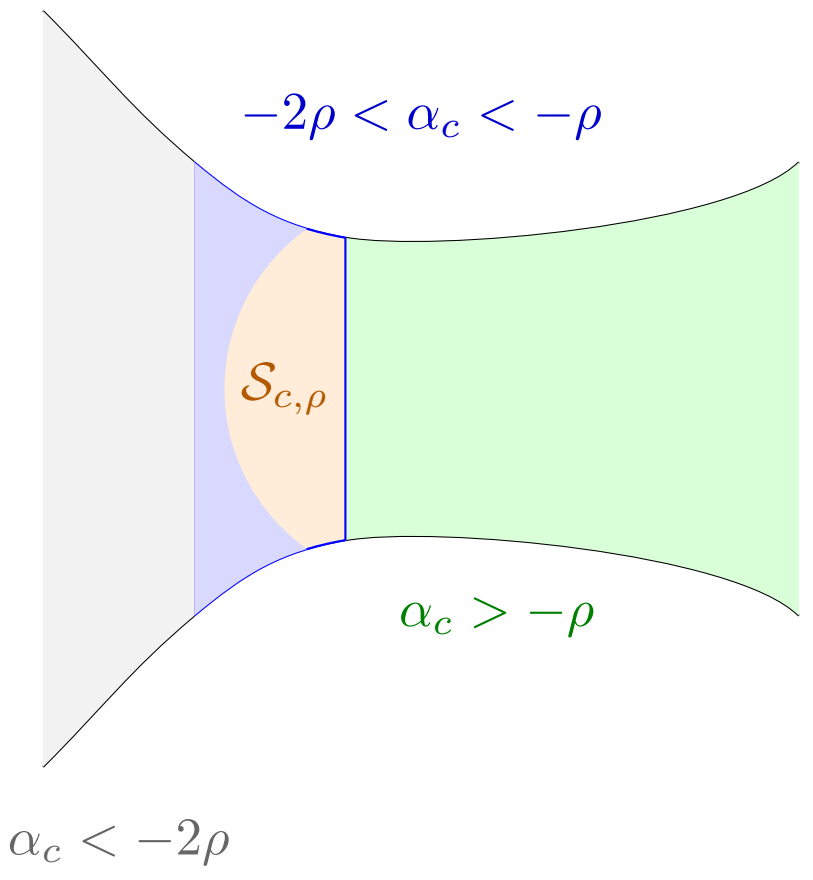

Since $\hat{Y}=Y_1$ in the coordinate ball, we have
\[
\hat{Y} = \frac{1}{2} \tanh( 2\xi_1 + 2x_*(c)) \partial_{\xi_1} + \tanh^2( 2\xi_1 + 2x_*(c)) \eta_1 \partial_{\eta_1} + ( \eta_2 - A_1(x_*(c)))\partial_{\eta_2}.
\]
We then find $X_c - \hat{Y}=X_{E_c},$ where the local Hamiltonian $E_c$ is given by
\begin{equation}\label{eq:firstinterlE}
E_c(\zeta) = \left( a \xi_1 - \frac{1}{2} \tanh( 2\xi_1+ 2x_*(c)) \right)\eta_1 + b \xi_2 \eta_2 - \sigma_c \xi_1 \xi_2 - A_1(x_*(c))\xi_2.
\end{equation}
Recall that $\sigma_c= -\frac{1}{4} \delta_c ( a+1-b). $
We define the  local   Liouville vector field by
\[
\tilde{Y}_{\rm loc} := \hat{Y}+
X_{\ell_c E_c}.
\]
Note that    $\tilde{Y}_{\rm loc} = \hat{Y} $ if $\alpha_c \leq -2\rho$ and $\tilde{Y}_{\rm loc} =X_c$ if $\alpha_c \geq -\rho.$ It follows from  Lemma \ref{lem:Weinstein}  that $\tilde{Y}_{\rm loc}$ is positively transverse to $K_c^{-1}(0)$ in $\mathcal{B}_{r_*}( z_*(c))$ where $\alpha_c \geq -\rho.$ We consider the region $\{ \alpha_c \leq -2\rho\}$ below.

 We compute  
 \begin{align*}
 \tilde{Y}_{\rm loc}(K_c) &= (1-\ell_c) Y_1 (K_c) + \ell_c X_c(K_c) + E_c X_{\ell_c} (K_c) \\
 &= (1-\ell_c) Y_1 (K_c) + \ell_c X_c(K_c) + \chi'(\alpha_c) E_c X_{\alpha_c}(K_c)   .
 \end{align*}
In order to get estimates for each term in the right-hand side above, we need the following assertion, whose proof is included in the appendix.

\begin{lemma} \label{lem:transtionregion}
There exist constants
\[
   \rho_0>0,\;\;m_0>0,
\;\; M_0>0,\;\; \kappa>0
\]
such that  if  $0<\rho<\rho_0$ and $c_2 < c < c_2+\rho^3,$ then    every  $ \zeta \in \mathcal{S}_{c,\rho}$  
satisfies
\[
m_0 \rho \leq \lvert \zeta \rvert \leq M_0 \rho, \qquad \xi_1 \leq - \kappa \rho.
\]
\end{lemma}

We now decrease $\rho_0>0$ obtained in   Lemma \ref{lem:transtionregion} so that  $ \rho_0< \min \left\{1, \varepsilon_*^{1/3}\right\}.$

We first study the term $Y_1(K_c).$ We write $Y_1 = Y_{\rm const} + Y_{\rm err},$
where
\begin{align*}
Y_{\rm const} = \frac{1}{2} \tanh(2x_*(c)) \partial_{\xi_1} - \frac{1}{8} \sinh(2x_*(c)) \partial_{\eta_2} 
\end{align*}
is the constant part of $Y_1.$  
We compute
\begin{equation}\label{eq:TTT}
    Y_{\rm const}(Q_c)  = T(c)\xi_1, \quad T(c):=\frac{\lambda_1(c)}{2} \tanh(2x_*(c)) <0 .
    \end{equation}
Since $Q_c$ is a quadratic form whose coefficients smoothly depend on $c,$ by compactness of $I_*,$ there exists $C>0$ such that   
\[
 \lvert \partial_{\zeta_j} Q_c(\zeta) \rvert \leq C \lvert \zeta \rvert , \;\;\;\; c \in I_* , \; \lvert \zeta \rvert \leq r_*,
 \]
 where  $\partial_{\zeta_j}$ denotes the partial derivative along the coordinate  directions.
Hence,
 \begin{align*}
     \lvert Y_{\rm err}(Q_c) \rvert & \leq  \frac{1}{2} \lvert \tanh(2\xi_1 + 2x_*(c))- \tanh(2x_*(c)) \rvert   \lvert \partial_{\xi_1}Q_c\rvert  \\   &  \quad  + \tanh^2( 2\xi_1 + 2x_*(c))\lvert \eta_1 \partial_{\eta_1}Q_c\rvert  + \lvert \eta_2 \partial_{\eta_2}Q_c\rvert\\  & \leq C'\lvert \zeta \rvert^2
 \end{align*}
for some $C'>0,$
where in the last step we used  the mean value theorem to   find $\theta \in (0,1)$ such that
\[
\lvert \tanh(2\xi_1 + 2x_*(c)) - \tanh (2x_*(c)) \rvert \leq \lvert 2 \sech^2( 2 \theta \xi_1 + 2x_*(c))  \xi _1 \rvert \leq 2 \lvert \xi_1 \rvert.
\]
Recall from  \eqref{eq:CTay} that 
\[
\lvert \mathcal{R}_c(\zeta) \rvert \leq C_{\rm Tay} \lvert \zeta\rvert^3, \quad \lvert \partial_{\zeta_j} \mathcal{R}_c(\zeta)\rvert \leq C_{\rm Tay}\lvert \zeta\rvert^2 , \qquad c \in I_*,\;\; \lvert \zeta \rvert \leq r_*.
\]
 Then arguing as above, we find a constant $C''>0$ such that
\[
\lvert Y_1(\mathcal{R}_c) \rvert \leq C''  {\lvert \zeta\rvert^2.} 
\]
Consequently, we obtain
\begin{equation}\label{eq:firstesimatedsdfs}
\lvert Y_1(K_c) - T(c)\xi_1\rvert \leq \lvert Y_{\rm err}(Q_c) \rvert + \lvert Y_1(\mathcal{R}_c) \rvert \leq (C'+C'') \lvert \zeta\rvert^2    .
\end{equation}
 Since $T$   depends smoothly on $c,$ applying the mean value theorem, we find a constant $D>0$ such that
\[
\lvert T(c) - T(c_2) \rvert \leq D \lvert c - c_2 \rvert, \qquad c \in I_*.
\]
 This together with \eqref{eq:firstesimatedsdfs} yields
 \[
 \lvert Y_1(K_c) - T(c_2)\xi_1\rvert  \leq  D' ( \lvert c - c_2 \rvert \lvert \zeta \rvert + \lvert \zeta \rvert^2 )
 \]
for some $D'>0.$ Since $T(c_2)<0,$ in view of Lemma \ref{lem:transtionregion}  we now  have  for $c\in(c_2,c_2+\rho^3)$
\begin{equation}\label{eq:Y1Kc}
    Y_1(K_c)    \geq \left( - T(c_2) \kappa    -  D' M_0( \rho^3   + M_0  \rho  ) \right)\rho.
\end{equation}

We now move to the term $X_c(K_c).$ By   Lemma \ref{lem:Weinstein},   we have
\[
X_c(K_c) \geq \mathfrak{w} \lvert \zeta\rvert^2
\]
on $ \mathcal{B}_{r_*}(z_*(c)).$  Then by Lemma \ref{lem:transtionregion} we have
\[
X_c(K_c) \geq \mathfrak{w}  m_0^2 \rho^2\qquad \text{ on } \; \mathcal{S}_{c,\rho}.
\]

Finally, to study the term $E_cX_{\alpha_c}(K_c),$ we first    introduce the linear term of $E_c$
\begin{equation}\label{eq:Elin}
E_{\rm lin} =  -A_1(x_*(c)) \xi_2 -\frac{1}{2} \tanh(2x_*(c)) \eta_1 
\end{equation}
and set $E_{\rm err} = E_c- E_{\rm lin},$ so that $E _c= E_{\rm lin}+E_{\rm err}.$ 
We estimate
\begin{equation*}\label{eq:esimate1}
\begin{aligned}
\lvert E_{\rm err} \rvert &\leq \lvert  ( a \xi_1 - \frac{1}{2} \tanh( 2\xi_1+ 2x_*(c) ) + \frac{1}{2} \tanh(2x_*(c)) )   \eta_1 \rvert  +   b   \lvert \xi_2   \eta_2 \rvert +\lvert \sigma_c\rvert\lvert \xi_1  \xi_2\lvert \\
& \leq ( \lvert a \rvert + 1 ) \lvert \xi_1  \eta_1\rvert +   b   \lvert \xi_2   \eta_2 \rvert +\lvert \sigma_c\rvert\lvert \xi_1  \xi_2\lvert \\
& \leq C_E \lvert \zeta \rvert^2, \qquad C_E:=2( \lvert a \rvert +1 +b + \sup_{c \in I_*}\lvert \sigma_c \rvert)>0,
\end{aligned}
\end{equation*}
where in the second inequality we used the mean value theorem. 
We also compute
\begin{align*}
X_{\alpha_c} (Q_c) &= - \partial_{\xi_2}Q_c - ( \mu_c - \frac{1}{4} )\partial_{\eta_1}Q_c\\
&= -  \mu_c \eta_1 - ( \lambda_2(c) + \frac{1}{4}\mu_c) \xi_2\\
&= \nu_c E_{\rm lin},
\end{align*}
where
\begin{equation}\label{eq:nuc}
\nu_c:=\frac{8 \lambda_2(c)}{ \sinh(2x_*(c)) - \tanh(2x_*(c)) } = \frac{2\mu_c}{\tanh(2x_*(c))}  >0.
\end{equation}
Since $\nu_c$ depends continuously on $c,$   we find $\nu_\pm >0$ such that $\nu_- \leq \nu_c \leq \nu_+$ for all $c \in I_*.$
We also find a constant $M>0$   such that 
\begin{equation}\label{eq:choseofM}
\lvert \mu_c - \frac{1}{4} \rvert \leq M , \quad c \in   I_*.
\end{equation}
It follows from \eqref{eq:CTay}  that
\[
   \lvert X_{\alpha_c}( \mathcal{R}_c)(\zeta) \rvert  \leq  \lvert \partial_{\xi_2}\mathcal{R}_c(\zeta)\rvert + \lvert \mu_c - \frac{1}{4} \rvert \lvert \partial_{\eta_1}\mathcal{R}_c(\zeta)\rvert  \leq (M+1) C_{\rm Tay}\lvert \zeta\rvert^2 .
\]
Consequently, we obtain
\begin{align*}
E_c X_{\alpha_c}(K_c)   &= ( E_{\rm lin}+ E_{\rm err})( X_{\alpha_c}(Q_c) +  X_{\alpha_c}(\mathcal{R}_c)) \\
&= ( E_{\rm lin}+ E_{\rm err})( \nu_c E_{\rm lin} +  X_{\alpha_c}(\mathcal{R}_c))\\
& = \nu_c E_{\rm lin}^2 + \nu_c E_{\rm lin }E_{\rm err} + E_{\rm lin}X_{\alpha_c}(\mathcal{R}_c) + E_{\rm err}X_{\alpha_c}(\mathcal{R}_c)\\
& \geq \nu_c E_{\rm lin}^2 - ( \nu_+C_E+ (M+1)C_{\rm Tay})\lvert E_{\rm lin}\lvert \lvert \zeta \rvert^2 - (M+1)C_EC_{\rm Tay}\lvert \zeta \rvert^4\\
& \geq \nu_c E_{\rm lin}^2 - C\lvert E_{\rm lin}\lvert \rvert \zeta \lvert^2  - C\lvert \zeta \rvert^4 ,
\end{align*}
where $C:=2 \max \{ \nu_+ C_E, (M+1)C_{\rm Tay}, (M+1)C_EC_{\rm Tay}\}>0.$ Applying Young's inequality $ab \leq  {a^2}/({2\delta})+  {\delta b^2}/{2}$  
with $a = \lvert E_{\rm lin}\rvert, b = C\lvert \zeta \rvert^2$ and $\delta = 1/\nu_c ,$ we obtain
\[
C \lvert E_{\rm lin}\rvert \lvert \zeta\rvert^2 \leq \frac{\nu_c}{2}E_{\rm lin}^2 + \frac{C^2}{2\nu_c} \lvert \zeta \rvert^4\leq \frac{\nu_c}{2}E_{\rm lin}^2 + \frac{C^2}{2\nu_-} \lvert \zeta \rvert^4,
\]
which yields 
\[
  E_cX_{\alpha_c} (K_c) \geq   \frac{\nu_c}{2}E_{\rm lin}^2 - \left( \frac{C^2}{2\nu_-} +C\right) \lvert \zeta \rvert^4     \geq -  C_0 \lvert \zeta \rvert^4
\]
with $C_0 := C^2 / 2\nu_- + C>0.$ Since $0 \leq \chi' \leq 2 /\rho,$ we get
\[
 \chi'(\alpha_c) E_cX_{\alpha_c} (K_c) \geq    -  \frac{2C_0}{\rho} \lvert \zeta \rvert^4
\]
An application of Lemma \ref{lem:transtionregion} yields
\[
\chi'(\alpha_c) E_c X_{\alpha_c}(K_c) \geq - 2 C_0 M_0^4 \rho^3 =: - C_1 \rho^3.
\]

Set $C_2 := \min \{ - \frac{1}{2} T(c_2) \kappa, \mathfrak{w} m_0^2\}>0.$ 
Taking smaller  $\rho_0>0$ if necessary, we may assume that  
\[
D' M_0( \rho_0^3   + M_0  \rho_0 )    \leq - \frac{1}{2}T(c_2)\kappa, \quad C_1 \rho_0 \leq \frac{1}{2}C_2.
\]
Then for every $0<\rho<\rho_0$ and every $c_2<c<c_2+\rho^3,$ by \eqref{eq:Y1Kc} we have
\begin{equation}\label{eq:T1forkcrho}
Y_1(K_c) \geq - \frac{1}{2}T(c_2) \kappa \rho >0 \qquad \text{ on } \; \mathcal{S}_{c,\rho},
\end{equation}
and therefore
\begin{align*}
\tilde{Y}_{\rm loc}(K_c) &\geq - \frac{1}{2}(1-\ell_c)T(c_2) \kappa \rho  + \ell_c\mathfrak{w} m_0^2 \rho^2  - C_1 \rho^3\\
& \geq\left(  - \frac{1}{2}(1-\ell_c)T(c_2) \kappa   + \ell_c\mathfrak{w} m_0^2 \right) \rho^2  - C_1 \rho^3\\
 &\geq  C_2 \rho^2  - C_1 \rho^3 \\
 &\geq \frac{1}{2}C_2 \rho^2
\end{align*}
on $\mathcal{S}_{c,\rho},$ where in the second and last steps we used  that $0<\rho<\rho_0<1.$

It remains to consider the region where $\alpha_c \leq -2\rho.$ 
Recall that in this region we have $\tilde{Y}_{\rm loc}=\hat{Y}=Y_1.$
Since the above estimates hold  for a smaller $r_*>0, $ after decreasing it if necessary, we assume that $(M+1) r_* < 2 \rho_0,$ where $M>0$ is as in \eqref{eq:choseofM}. We then have $\lvert \alpha_c (\zeta)\rvert \leq M \lvert \xi_1 \rvert + \lvert \eta_2 \rvert \leq  (M+1)\lvert \zeta \rvert.$
Then for a given    $\zeta \in K_c^{-1}(0)  \cap \mathcal{B}_{r_*}( z_*(c))$ with $\alpha_c(\zeta) \leq-2 \rho, $ we have
\[
\frac{\lvert \alpha_c(\zeta)\rvert}{2} \leq \frac{M+1}{2} \lvert \zeta \rvert \leq \frac{M+1}{2}r_* < \rho_0.
\]
We now set $s:=  - \alpha_c(\zeta)/2>0.$  
Then applying \eqref{eq:T1forkcrho} with $\rho>0$ replaced by $s>0$  we find
\[
\tilde{Y}_{\rm loc}(K_c)(\zeta)=Y_1(K_c)(\zeta) \geq - \frac{1}{2}T(c_2) \kappa s >0 . 
\]

We   have proved

\begin{lemma}\label{lem:firstinterpolatoinsuccess} After decreasing $\rho_0>0$ and $r_*>0$ if necessary,  for every $0<\rho<\rho_0$ and $c_2<c<c_2+\rho^3,$ the local Liouville vector field $\tilde{Y}_{\rm loc}$ is positively transverse to $K_c^{-1}(0)$ in the coordinate ball $\mathcal{B}_{r_*}(z_*(c)).$

\end{lemma}

Before proceeding further, we obtain the following estimates that will be used later in the second interpolation step.

Set
\begin{align*}
\delta_* =\sup_{c \in I_*}  \lvert \delta_c\rvert  >0, \quad   \lambda_{*,1} =\sup_{c \in I_*}   \lvert \lambda_1(c)\rvert  >0, \quad   \lambda_{*,2}:=\inf_{c \in I_*} \lvert \lambda_2(c)\rvert  >0.
\end{align*}
Then we find
\begin{align*}
\lvert \zeta \rvert^2 &=\xi_1 ^2 + \xi_2^2 +\eta_1^2 +\eta_2^2\\
&= \xi_1^2 + \xi_2^2 + \left( \eta_1 +\frac{1}{4}\xi_2 - \frac{1}{4}\xi_2 \right)^2 + \left( \eta_2 +\frac{\delta_c}{4}\xi_1 -\frac{\delta_c}{4}\xi_1\right)^2\\
&\leq  \xi_1^2 + \xi_2^2 + 2\left( \eta_1 +\frac{1}{4}\xi_2\right)^2 + \frac{1}{8}\xi_2 ^2 +   2\left(\eta_2 +\frac{\delta_c}{4}\xi_1\right)^2 +\frac{\delta_c^2}{8}\xi_1 ^2\\
&\leq C_0 \left(  \xi_1^2 + \xi_2^2 +  \left( \eta_1 +\frac{1}{4}\xi_2\right)^2 + \left( \eta_2 +\frac{\delta_c}{4}\xi_1 \right)^2    \right),  
\end{align*}
where $C_0:=\max \{ 1+ \delta_*^2/8, 2\}.$
We further estimate with $C_1:=\max\{1, 1/\lambda_{*,2}\}$
\begin{align*}
& \xi_2^2 +  \left( \eta_1 +\frac{1}{4}\xi_2\right)^2 + \left( \eta_2 +\frac{\delta_c}{4}\xi_1 \right)^2   \\ \leq & \; C_1 \left(  \lambda_2(c)\xi_2^2 +  \left( \eta_1 +\frac{1}{4}\xi_2\right)^2 + \left( \eta_2 +\frac{\delta_c}{4}\xi_1 \right)^2 \right)\\
    =& \;C_1 \left( 2Q_c(\zeta) - \lambda_1(c) \xi_1^2\right)\\
 \leq & \;  C_1 \left( 2 \lvert Q_c(\zeta)\rvert +\lambda_{*,1}\xi_1^2 \right).
\end{align*}
Therefore, with $\mathcal{C}_Q:=C_0 \cdot \max\{2C_1, 1+C_1 \lambda_{*,1}\}$ we have
\begin{equation}\label{eq:boundeforzetaxc}
\lvert\zeta\rvert^2 \leq \mathcal{C}_Q( \lvert Q_c(\zeta)\rvert +\xi_1^2) , \qquad c \in I_*.
\end{equation}

Set 
\begin{align*}
m&:= \min \left\{ \lvert Q_{c_2}(w)\rvert + \alpha_{c_2}(w)^2\ :\  \lvert w \rvert =1\right\} >0,\\
d_0& := \frac{1}{4}\left( \mu_{c_2} + \frac{1}{4}(\delta_{c_2}-1) - \sqrt{-\lambda_1(c_2)}\right)>0.
\end{align*}
 After decreasing  $r_*>0$ if necessary, we may assume that 
\begin{equation}\label{eq:SDFK}
C_{\rm Tay} r _* \leq \frac{m}{4}, \quad  \mathcal{C}_Q C_{\rm Tay} r_* \leq \frac{1}{2}    \quad \text{ and }  \quad 2 \sqrt{2 r_* C_{\rm Tay}} \mathcal{C}_Q^{3/4} \leq d_0,
\end{equation}
where $C_{\rm Tay}>0$ is as in \eqref{eq:CTay}.

Recall that  $\{Q_{c_2} =0\} \cap \{ \alpha_{c_2}=0\} = \{0\}, $ see Remark \ref{rmk:separation}. 
Choose $0<\rho_1<\rho_0,$ where $\rho_0$ is as in Lemma \ref{lem:firstinterpolatoinsuccess}. After  decreasing $\rho_1>0$ if necessary,  continuity  {in $c$} and homogeneity {in $\zeta$} of $\lvert Q_c(\zeta)\rvert + \alpha_c(\zeta)^2$ give  
\begin{equation}\label{eq:Qminunmd}
\lvert Q_c(\zeta)\rvert+\alpha_c(\zeta)^2 \geq \frac{m}{2}\lvert \zeta\rvert^2 , \quad  c_2 \leq c \leq c_2+\rho_1^3.
\end{equation}
Also,  again by continuity we may assume that
\begin{equation}\label{eq:cdsd02}
\mu_{c} + \frac{1}{4}(\delta_{c}-1) - \sqrt{-\lambda_1(c)} \geq 2d_0, \qquad   c_2 \leq c \leq c_2+\rho_1^3.
\end{equation}

Recall that $r_*>0$ is chosen sufficiently small, so that the four    (closed) coordinate $r_*$-balls centered at the four critical points are pairwise disjoint for every $c \in I_*.$  
Applying Lemma \ref{lem:locallemman} with this radius and decreasing $\rho_1>0$ further if necessary, we obtain   $\delta_0 \in ( 0, \min_{c \in I_*} x_*(c))$ such that for every $0<\rho<\rho_1$ and every $c_2 <c<c_2+\rho^3,$ each point in $K_c^{-1}(0)$ with $\lvert x_1 - x_*(c)\rvert <\delta_0$ lies either in $\mathcal{B}_{r_*}(z_*(c))$ or in the corresponding coordinate $r_*$-ball centered at $(x_*(c), \pi, 0 ,-A_1(x_*(c))).$

\begin{figure}[ht]
\centering
\begin{tikzpicture} 
  \node[anchor=south west,inner sep=0] at (0,0) {\includegraphics[width=9cm]{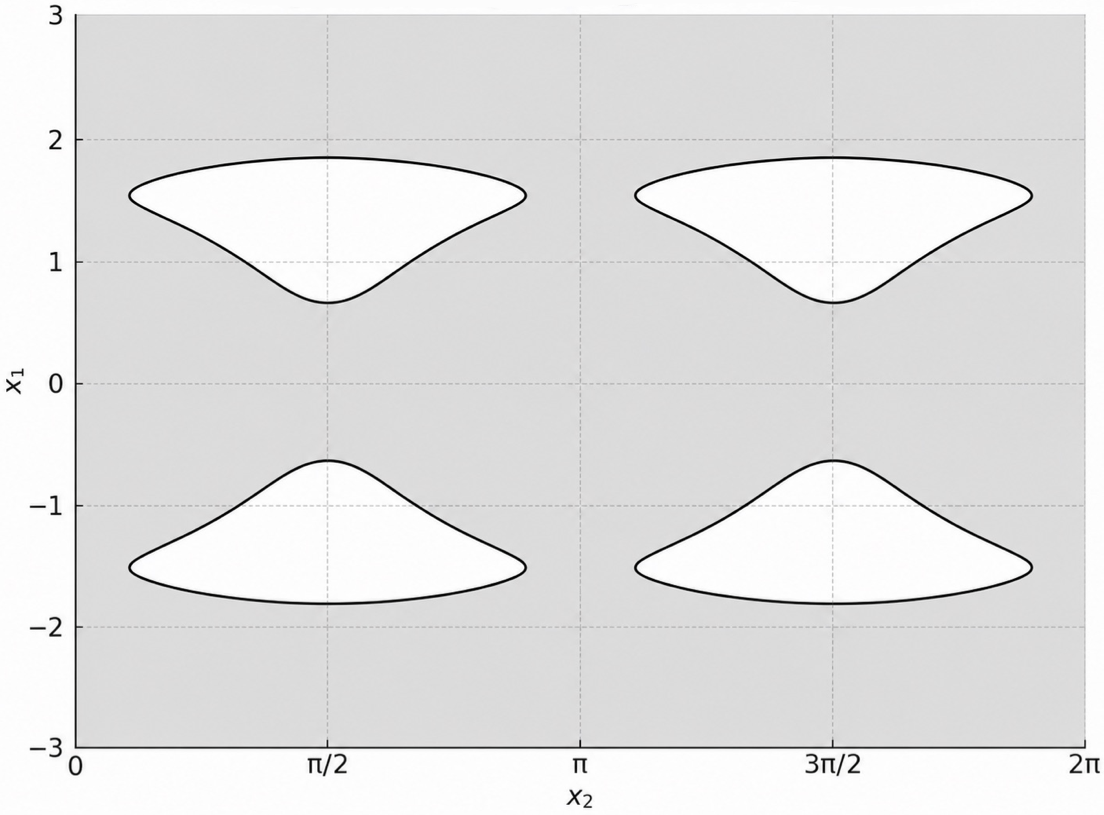}};

  \draw[red!30!black,dashed]   (0.68,4.8 ) -- (8.83,4.8 );
  \node[left,blue] at (0.58,5.05) {\small $x_* $};
\fill[blue] (0.65,5.05) circle (2 pt);

 \draw[red!30!black,dashed]   (0.68,2.35 ) -- (8.83,2.35 );
  \node[left,blue] at (0.58,2.1) {\small $-x_* $};
\fill[blue] (0.65,2.1) circle (2 pt);

\begin{scope}[xshift=7.15cm, yshift=2.35cm]
 \draw[red!30!black,decorate,decoration={brace,amplitude=5pt}]
  (2,2.4) -- (2,0)
  node[midway,right=8pt] {$\mathcal{K}_c$};
\end{scope}  
 \end{tikzpicture}
\caption{The center part $\mathcal{K}_c,$ projected to the position plane}
\label{Fig:Lemma48}
 \end{figure}

 For every $c \in I_*,$ we consider the non-empty closed subset
 \begin{equation}\label{eq:centralpart}
 \mathcal{K}_c:= \{ z \in K_c^{-1}(0) \mid \lvert x_1 \rvert \leq x_*(c)-\delta_0\}.
 \end{equation}
Note that $\mathcal{K}_c$ contains none of the four critical points. See Figure \ref{Fig:Lemma48}.

\begin{lemma}\label{lem:additionalestimateforhaty} There exist sufficiently small $0<\varepsilon_{\ell}  <\varepsilon_*$ and   $\sigma_0>0,$ independent of $c$ and $\rho,$ such that 
\[
\hat{Y}(K_c) \geq \sigma_0  
\]
on $\mathcal{K}_c$ for every $ c \in I_{\ell}:=[c_2, c_2+\varepsilon_{\ell}],$
where $\hat{Y}$ is as in Remark \ref{remark:slightlyaboveL2}.

\end{lemma}
\begin{proof}
Let $\hat{\mathcal{K}} \subset T^*(\R \times S^1)$ be a compact subset such that $\mathcal{K}_c \subset \hat{\mathcal{K}}$ for all $c \in I_*.$  The construction in   Remark \ref{remark:slightlyaboveL2} provides $\hat{Y}(K_{c_2})>0$ along $\mathcal{K}_{c_2}.$ Set
\[
\sigma_0:= \frac{1}{2} \min_{\mathcal{K}_{c_2}}\hat{Y}(K_{c_2})>0.
\]
Assume by contradiction that no such $\varepsilon_{\ell}$ exists. Then  we find   sequences $c_n \in I_*$ converging to $c_2$ and $z_n \in \mathcal{K}_{c_n}$ such that $\hat{Y}(K_{c_n}) < \sigma_0$ at $z = z_n.$ 
By compactness, after passing to a subsequence, we have $z_n \to \bar{z} \in \hat{\mathcal{K}} .$ Then continuity of $K_c$ and $x_*(c)$ in $c$ gives $\bar{z} \in \mathcal{K}_{c_2}.$
By continuity we derive a contradiction
\[
2 \sigma_0 \leq \hat{Y}({K}_{c_2})(\bar{z}) = \lim_{n \to \infty} \hat{Y}(K_{c_n})(z_n) \leq \sigma_0,
\]
which completes the proof.
\end{proof}


\subsection{Capping the energy level}
\label{sec:capping}

We keep $r_*, \delta_0, \rho_1, \varepsilon_{\ell} $ and $ \sigma_0$ as in the previous section with $I_{\ell}=[c_2, c_2+\varepsilon_{\ell}] \subset I_*$ and $0<\rho_1<\rho_0<\min\{ 1, \varepsilon_*^{1/3}\}.$

Recall that  $x_*=x_*(c)$ is a strict local maximum of $V_{c,1}$ and satisfies 
  \begin{equation}\label{eq:twoinde} 
  0< V_{c,1}(x_*(c))< -\frac{c}{4}
\end{equation}
for $c>c_2$ sufficiently close to $c_2.$ See Remark \ref{remark:valuev1xstar}. After decreasing $\rho_1>0$ further if necessary, we may assume that \eqref{eq:twoinde} holds for every $c_2<c<c_2+\rho_1^3.$ Let $0<\rho<\rho_1.$ Choose $\varrho_0=\varrho_0(\rho)>0$ sufficiently small so that 
\begin{equation}\label{eq:varrho0intro}
\varrho_0 < \min\{\rho, \delta_0 /4, \varepsilon_{\ell}^{1/3}\}
\end{equation}
and 
\begin{equation}\label{eq:positiefVatx}
0< V_{c,1}(x_*(c)+s) < -\frac{c}{4},\qquad 
V_{c,1}'(x_*(c)+s) <0
\end{equation}
for every $c \in (c_2, c_2+ \rho^3)$ and $s \in (0, 4\varrho_0].$

Fix $\varrho \in (0,\varrho_0)$ and $c_2 < c \leq  c_2 +  \varrho^3 .$ Because of \eqref{eq:varrho0intro}, we have 
\[
c_2 < c \leq c_2 +\varrho^3  < c_2 + \min \{ \rho^3, \varepsilon_{\ell}\},
\]
and hence the assertions of Lemma \ref{lem:firstinterpolatoinsuccess}  and Lemma \ref{lem:additionalestimateforhaty} hold. Moreover, the strip $\{ x_*(c) < x_1 < x_*(c)+4\varrho\}$ intersects $K_c^{-1}(0)$ only in the two corresponding coordinate $r_*$-balls. By symmetry, the analogous strip near $x_1=-x_*(c)$ intersects $K_c^{-1}(0)$ near the two other necks. In other words, the two strips intersect $K_c^{-1}(0)$ only in small neighborhoods of the four neck regions.

By symmetry, as above, it suffices to describe a perturbation of the potential near $(x_*(c),0).$ 
We shall modify the potential by adding a smooth function depending only on  $\xi_1=x_1 - x_*(c).$ This function vanishes for $\xi_1 \leq 2\varrho,$  creates a new local minimum at $\xi_1=3\varrho,$ and  makes the   modified potential diverge  as $\xi_1 \to 4\varrho^-.$

The proof of the following lemma, motivated by the argument in \cite[Section 3.1]{dPKSS24}, will be given in the appendix.

\begin{lemma}\label{lem:cappingfunction}
After decreasing    $\varrho_0>0$ if necessary, there exist  constants \[
a_0>0, \quad   0<\mathfrak{c}_0<\frac{a_0}{4},  \quad \mathfrak{c}_1>0,  \]
independent of $\rho,c$ and $\varrho,$
such that the following holds:  Let $\varrho \in (0,\varrho_0)$   and $c_2<c \leq    c_2 +   \varrho^3  .$ Then   
\begin{equation}\label{eq:lowerbuonda0}
- V_{c,1}'( x_*(c)+s)   \geq a_0s, \quad \; \forall   s \in [0,4 \varrho_0],
\end{equation}
and  {there exists} a non-decreasing smooth function $\psi=\psi_{c,\varrho} \colon (-\infty, 4 \varrho) \to [0, +\infty)$ satisfying
\begin{enumerate}
    \item[i)] $\psi \equiv 0$  on $(-\infty, 2\varrho]$;

    \item[ii)] $ (s-3\varrho)(V_{c,1}'( x_*(c)+s )+\psi'(s)) \geq \mathfrak{c}_0(s-3\varrho)^2, \quad \forall   s \in  [\varrho, 4\varrho) $;

    \item[iii)] $\psi(s) \to +\infty$ as $s \uparrow 4\varrho$; 

\item[iv)] $0< d_{c,\varrho}:= -c/4 - V_{c,1}(x_*(c)+3\varrho) - \psi(3\varrho) \leq \mathfrak{c}_1 \varrho^2 .$

\end{enumerate}
In particular,
\begin{align*}
&\begin{cases} \psi'(s) < -V_{c,1}'(x_*(c)+s) & 2\varrho < s < 3\varrho \\
\psi'(3\varrho) = -V_{c,1}'(x_*(c)+3\varrho)  &   \\
\psi'(s) > -V_{c,1}'(x_*(c)+s) & 3\varrho<s<4\varrho\end{cases}
\end{align*}
 and $\psi''(3\varrho) \geq -V_{c,1}''(x_*(c)+3\varrho) + \mathfrak{c}_0. $ 
    
\end{lemma}

\begin{figure}[t]
\centering
\includegraphics[width=0.9\textwidth]{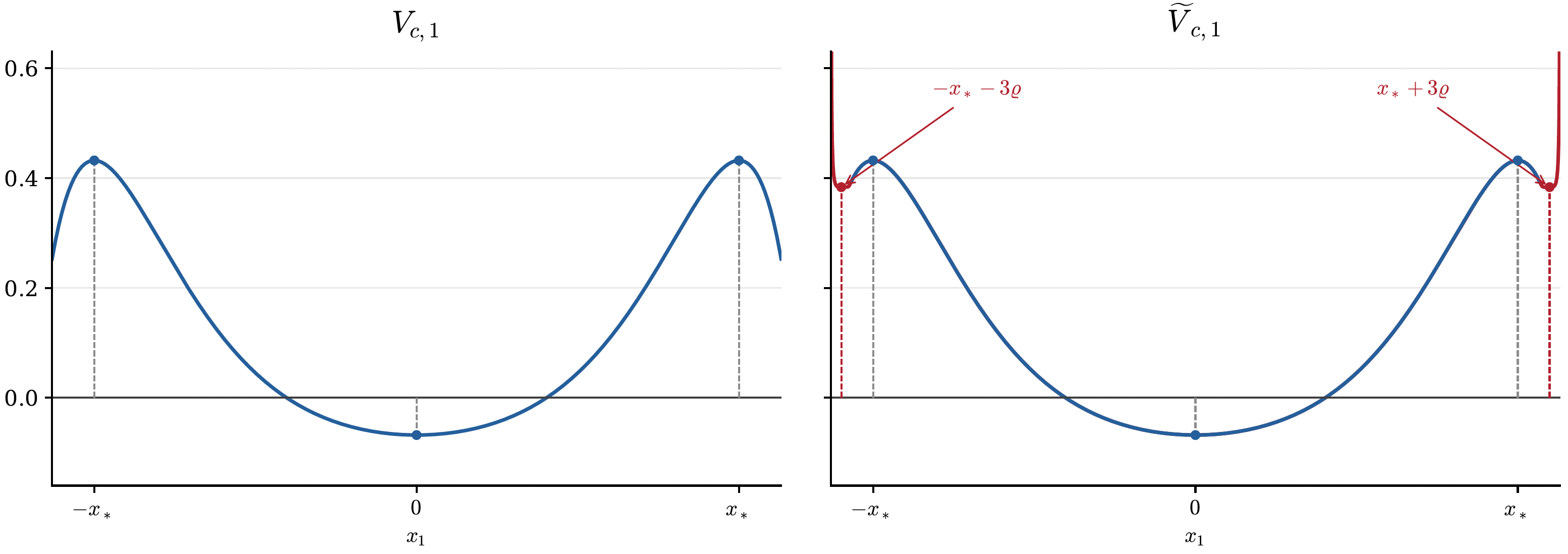}
\caption{An illustration of the potential $V_{c,1}$ and a smooth cap $\tilde{V}_{c,1}.$ The blue curves denote the unmodified parts and the red denote the capped parts}
\label{fig_pert}
\end{figure}

Let $\psi$ be as in the above statement. 
We define 
\[
\bar{V}_{c,\varrho, 1}(x_1):=V_{c,1}(x_1) + \psi( x_1 - x_*(c)), \quad x_1 < x_*(c)+4\varrho.
\]
We understand this perturbation only in a small neighborhood of $x_*(c).$
Note that $\bar{V}_{c,\varrho,1} = V_{c,1}$ for $x_1 \leq x_*(c)+2\varrho.$  Therefore, the critical points $x_1 = 0, x_*(c)$ of $V_{c,1}$ remain critical points of $\bar{V}_{c,\varrho,1}$ with the same critical values.
We compute
\[
\bar{V}_{c,\varrho,1}'(x_1) = V_{c,1}'(x_1) + \psi'( x_1 - x_*(c)).
\]
By construction of $\psi,$ in the transition region $2\varrho<x_1 - x_*(c)<4\varrho$ there is exactly one critical point $x_1 = x_*(c)+3\varrho$ of $\bar{V}_{c,\varrho,1},$ which is a strict local minimum.
Its critical value satisfies
\begin{equation}\label{eq:valueoftildeC}
0 < \bar{V}_{c,\varrho,1}(x_*(c)+3\varrho) = V_{c,1}(x_*(c)+3\varrho) + \psi(3\varrho) < - \frac{c}{4}
\end{equation}
We perform the symmetric construction near $-x_*(c).$ This yields the   perturbed function, denoted by $\tilde{V}_{c,\varrho,1}(x_1).$  See Figure \ref{fig_pert}. 
Note that $\tilde{V}_{c,\varrho,1}(x_1) \to +\infty$ as $x_1 \to x_*(c)+4\varrho$ or $x_1 \to -x_*(c)-4\varrho.$
We now define  
 \[
 \tilde{V}_{c,\varrho}(x_1, x_2):= \tilde{V}_{c,\varrho,1}(x_1) + V_{c,2}(x_2) 
 \]
and  
\[
\tilde{K}_{c,\varrho} ( x,y) := \frac{1}{2}\left[ \left( y_1 + A_2(x_2)\right)^2 + \left( y_2 +A_1(x_1)\right)^2     \right] + \tilde{V}_{c,\varrho}(x_1, x_2).
\]
For brevity, in the following we suppress the subscript $\varrho$ from $\tilde{V}_{c,\varrho,1}, \tilde{V}_{c,\varrho}$ and $\tilde{K}_{c,\varrho}.$
Since $\tilde{V}_{c} $ has no critical point with value $0,$ we see that $0$ is a regular value of $\tilde{K}_{c}.$ It is straightforward to see that the perturbed energy level $\tilde{\Sigma}_{c}:=\tilde{K}_c^{-1}(0)$ is compact and diffeomorphic to $S^1 \times S^2.$  
The capping terms vanish on $\{ \lvert x_1 \rvert < x_*(c)+2\varrho\},$ which contains the central part $\mathcal{K}_c$   (see \eqref{eq:centralpart}) for all $c_2 <c \leq c_2+\varrho^3.$ Therefore, we have $\tilde{K}_c=K_c$ and $d\tilde{K}_c = dK_c$ on a small neighborhood of $\mathcal{K}_c.$
Hence, by Lemma \ref{lem:additionalestimateforhaty} we then have
\[
\hat{Y}(\tilde{K}_c)=\hat{Y}(K_c) \geq \sigma_0 \qquad \text{ on } \; \mathcal{K}_c.
\]


\subsection{Second interpolation}\label{sec:secondinterpo}

We keep $r_* , \delta_0, \rho_1, \varepsilon_{\ell}$ and $\sigma_0$ as in Section \ref{sec:firstin}.  Fix $0<\rho<\rho_1.$ We then choose a sufficiently small $\varrho_0 >0,$ so that the capping function $\psi=\psi_{c,\varrho}$ described in  Lemma \ref{lem:cappingfunction} exists associated with the constants $a_0, \mathfrak{c}_0$ and $\mathfrak{c}_1.$ In particular $\varrho_0$ satisfies \eqref{eq:varrho0intro}, so $\varrho_0^3 < \varepsilon_{\ell} < \varepsilon_*.$  
Hence, for every $0<\varrho<\varrho_0$ and every $c_2 < c \leq c_2 +\varrho^3,$ the conclusions in the previous sections hold. In particular, every point of $K_c^{-1}(0)$ with $0<x_1 - x_*(c)<4\varrho$ lies in one of the two corresponding coordinate $r_*$-balls with $x_1 = x_*(c).$

We continue to  work in the local coordinates $\zeta = (\xi, \eta)$ centered at $z_*(c).$  The     new critical point  $(x_*(c)+3\varrho, 0,0,-A_1(x_*(c)+3\varrho) )$ is now represented by 
\[
\zeta_{c,\varrho}=(3\varrho, 0,0,\Delta_{c,\varrho}),\qquad \Delta_{c,\varrho}:=A_1(x_*(c)) - A_1( x_*(c)+3\varrho).
\]
We also write
\begin{equation}\label{eq:P1P2}
P_1= \eta_1 + A_2(\xi_2), \quad P_2 = \eta_2 - A_1(x_*(c))+A_1( x_*(c)+\xi_1)
\end{equation}
so that 
\[
\tilde{K}_c = \frac{1}{2}(P_1^2 +P_2^2) + \tilde{V}_{c,1}(x_*(c)+\xi_1) + V_{c,2}(\xi_2).
\]

\begin{lemma}\label{lem:perturbedsall} 
After decreasing $\varrho_0>0$ if necessary, the following holds for every $0<\varrho<\varrho_0$ and every $c_2<c \leq c_2+\varrho^3$: 
every  $ (x,y) \in \tilde{\Sigma}_c  \cap \{ \varrho \leq x_1 - x_*(c) \leq 4 \varrho \}$ satisfies either $\lvert x_2\rvert <\pi/4$ or $\lvert x_2 - \pi \rvert<\pi/4$ (understood modulo $2\pi$). Moreover, we find a constant $\mathfrak{C}>0,$ independent of $\rho, c$ and $\varrho,$ such that every such point 
satisfies   
\[
\lvert \zeta \rvert \leq \mathfrak{C}\cdot\varrho
\]
in the corresponding  local coordinates.

\end{lemma}
\begin{proof} Pick $z= (x,y) \in \tilde{\Sigma}_c  \cap \{ \varrho \leq x_1 - x_*(c) \leq 4 \varrho \}.$ In the local coordinates $\zeta = (\xi,\eta),$ this implies $\varrho \leq \xi_1 < 4 \varrho.$  
By property ii) in Lemma \ref{lem:cappingfunction}   the function $s \mapsto \tilde{V}_{c,1}(x_*(c)+s)$ is strictly decreasing on $[\varrho, 3\varrho)$ and strictly increasing on $(3\varrho, 4\varrho).$ Since   $x_*(c)+3\varrho$ is  a strict local minimum,  we have
\[
\tilde{V}_{c,1}(x_*(c)+\xi_1)  \geq \tilde{V}_{c,1}(x_*(c)+3\varrho) 
\]
and therefore we obtain using   $\tilde{K}_c=0$ that 
\begin{equation}\label{eq:sdlifhiesti}
V_{c,2}(x_2) \leq - \tilde{V}_{c,1}(x_*(c)+3\varrho)  =   d_{c,\varrho}+ \frac{c}{4} = V_{c,2}(0)+d_{c,\varrho} = V_{c,2}(\pi)+d_{c,\varrho}.
\end{equation}
Since   $0$ and $\pi$ are   the only  (strict local) minima of $V_{c,2},$   we find $\gamma>0,$ independent of $c$ close to $c_2,$ such that   
\[
V_{c,2}(x_2) \geq V_{c,2}(0)+\gamma= \frac{c}{4}+\gamma ,\quad \text{ whenever } \;   \lvert x_2  \rvert \geq \frac{\pi}{4} ,\; \; \lvert x_2 -\pi \rvert \geq \frac{\pi}{4}.
\]
Since $d_{c,\varrho} \leq \mathfrak{c}_1 \varrho^2,$ for $\varrho_0>0$ sufficiently small,  we have $0<d_{c,\varrho}<\gamma$ for every $0<\varrho<\varrho_0$ and hence  $V_{c,2}(x_2) \leq  {c}/{4}+d_{c,\varrho} <  {c}/{4}+\gamma.$
This proves the first assertion of the lemma.

 To prove the second assertion, we work near $x_2=0.$ The case near $x_2=\pi$ can be done in a similar way, using symmetries. 
Since $\lvert \xi_2 \rvert <\pi/4$ and    $c<c_3,$  we have a uniform estimate
\begin{align*}
V_{c,2}(\xi_2) -\frac{c}{4} &=\left( - \frac{c}{4} - \frac{1}{32}\cos^2 \xi_2 \right)  \sin^2 \xi_2 \\
&\geq   \left( - \frac{c_3}{4} - \frac{1}{32}  \right)  \sin^2 \xi_2\\
&\geq \mathcal{C}\xi_2 ^2  , \;\; \qquad \mathcal{C}:= \frac{4}{\pi^2} \left( - \frac{c_3}{4} - \frac{1}{32}  \right) >0  .
\end{align*}
It follows that 
\begin{equation}\label{eq:sdljxibound}
\xi_2^2 \leq \frac{d_{c,\varrho}}{\mathcal{C}}\leq \frac{\mathfrak{c}_1}{\mathcal{C}}\varrho^2.
\end{equation}
Moreover,  by \eqref{eq:sdlifhiesti} we find 
\[
P_1^2 +P_2^2=-2\tilde{V}_c(x_*(c)+\xi_1,\xi_2) \leq  2d_{c,\varrho} \leq 2\mathfrak{c}_1\varrho^2.
\]
Now \eqref{eq:P1P2} and the mean value theorem finish  the proof.
\end{proof}

After decreasing $\varrho_0>0$ further if necessary, we may assume that 
\begin{equation}\label{eq;crho-sd/2}
\mathfrak{C}\varrho_0 < \frac{r_*}{2} .
\end{equation}
It follows that every point in $\tilde{K}_c^{-1}(0)$ with $\varrho \leq \xi_1 = x_1 - x_*(c)< 4\varrho$ lies in one of the two  corresponding  coordinate $r_*/2$-balls centered at the critical points with $x_2=0$ and $x_2=\pi.$ In particular, every such point with $\lvert x_2 \vert<\pi/4$ lies in $\mathcal{B}_{r_*/2}(z_*(c)).$

Set 
\[
 \mathcal{U}_{c,\varrho}=\tilde{\Sigma}_c  \cap \{ \varrho \leq \xi_1 < 4 \varrho, \lvert \xi_2 \rvert<\pi/4 \} \cap \mathcal{B}_{r_*}(z_*(c)).
 \]
Define the radial vector field
\[
\mathcal{Z}_{c,\varrho}  = \frac{1}{2}(\zeta - \zeta_{c,\varrho} )\partial_{\zeta}= \frac{1}{2}[ (\xi_1-3\varrho)\partial_{\xi_1}+\xi_2\partial_{\xi_2}+\eta_1 \partial_{\eta_1}+(\eta_2-\Delta_{c,\varrho})\partial_{\eta_2}     ].
\]
 
\begin{proposition}  \label{prop:radislve} After decreasing $\varrho_0>0$ if necessary,  for    every $0<\varrho<\varrho_0$ and every $c_2 < c \leq c_2 +\varrho^3,  $ we have
\[
\mathcal{Z}_{c,\varrho} (\tilde{K}_c)>0\quad \text{ on } \; \mathcal{U}_{c,\varrho} .
\]
\end{proposition}
\begin{proof}

We compute
\begin{align*}
    2\mathcal{Z}_{c,\varrho}(\tilde{K}_c) &= P_1(\eta_1+A_2'(\xi_2)\xi_2) +P_2(\eta_2 - \Delta_{c,\varrho}+(\xi_1 - 3\varrho)A_1'(x_*(c)+\xi_1))\\
    &\;\; \; +(\xi_1 - 3\varrho) \tilde{V}_{c,1}'(x_*(c)+\xi_1)+\xi_2 V'_{c,2}(\xi_2).
\end{align*}
We estimate the terms  in the right-hand side above separately.

For the first term  $P_1(\eta_1+A_2'(\xi_2)\xi_2),$ we observe that the function  $F(\xi_2):= \xi_2 A_2'(\xi_2) - A_2(\xi_2)$ 
satisfies $F(0)=0$ and $F'(\xi_2) = \xi_2 A_2''(\xi_2) = - (\xi_2/2)\sin (2\xi_2).$ Since $\lvert \xi_2 \rvert<\pi/4,$ we have $\lvert \sin (2\xi_2) \rvert \leq 2\lvert \xi_2 \rvert$ and hence
\begin{align*}
   \lvert \xi_2 A_2'(\xi_2) - A_2(\xi_2)\rvert =  \lvert F(\xi_2) \lvert \leq \bigg| \int_0^{\xi_2} u^2 du \bigg| = \frac{1}{3}\lvert \xi_2 \rvert ^3 .
\end{align*}
Therefore, an application of Young's inequality provides 
\begin{equation}\label{eq:firstterm}
    \begin{aligned}
    P_1(\eta_1+A_2'(\xi_2)\xi_2)&=P_1 ^2 + P_1(  \xi_2 A_2'(\xi_2)-A_2(\xi_2))\\
    &\geq P_1^2 - \frac{1}{2}P_1^2 - \frac{1}{2}\lvert \xi_2A_2'(\xi_2) - A_2(\xi_2)\rvert^2\\
    &\geq \frac{1}{2}P_1^2 - \frac{1}{18}\xi_2^6. 
\end{aligned}
\end{equation}

We then consider the second term $P_2(P_2 +G(\xi_1)), $ where
\[
G(\xi_1) = (\xi_1 - 3\varrho)A_1'(x_*(c)+\xi_1) - A_1(x_*(c)+\xi_1)+A_1(x_*(c)+3\varrho).
\]
By integration by parts we have
\[
G(\xi_1) = (\xi_1 - 3\varrho)^2 \int_0^1 u A_1''(x_*(c)+3\varrho + u(\xi_1 - 3\varrho))du.
\]
Shrinking   $\varrho_0>0$ further if necessary, we may assume that $x_*(c) \leq x_*(c_2)+1/2$ and $\varrho_0\leq 1/8.$ Since $\varrho \leq \xi_1 <4 \varrho,$ we have for every $u \in [0,1]$
\[
0<x_*(c) + 3\varrho +u(\xi_1 - 3\varrho) \leq x_*(c_2)+1.
\]
Since $A_1''(\xi_1)=\frac{1}{2}\sinh(2\xi_1),$ we have
\[
\lvert G(\xi_1) \rvert \leq \mathcal{C}_0 \lvert \xi_1 - 3\varrho\rvert^2, \qquad \mathcal{C}_0:=\frac{1}{4}\sinh ( 2x_*(c_2)+2)>0.
\]
Then again by Young's inequality we obtain
\begin{equation}\label{eq:secondterm}
P_2(\eta_2 - \Delta_{c,\varrho}+(\xi_1 - 3\varrho)A_1'(x_*(c)+\xi_1)) \geq \frac{1}{2}P_2^2 - \frac{1}{2}\mathcal{C}_0^2 ( \xi_1 - 3\varrho)^4.
\end{equation}

For the third term $(\xi_1 - 3\varrho) \tilde{V}_{c,1}'(x_*(c)+\xi_1),$
Property ii) of Lemma \ref{lem:cappingfunction} provides
\begin{equation}\label{eq:thirdterm}
     (\xi_1-3\varrho)\tilde{V}_{c,1}'( x_*(c)+\xi_1 )  \geq \mathfrak{c}_0(\xi_1-3\varrho)^2.
\end{equation}

Finally,   for the last term we compute
\begin{equation}\label{eq:lastterm}
\begin{aligned}
\xi_2 V'_{c,2}(\xi_2) &= - \frac{1}{32}\xi_2 \sin(2\xi_2)( 8c+\cos (2\xi_2))\\
&\geq - \frac{1}{32} \frac{4}{\pi}\xi_2^2 (8c_3+1) \\
& \geq \frac{5}{4\pi}\xi_2^2.
\end{aligned}
\end{equation}

Combining \eqref{eq:firstterm}--\eqref{eq:lastterm} we get
\[
2\mathcal{Z}_{c,\varrho} (\tilde{K}_c)\geq \frac{1}{2}(P_1^2 +P_2^2) + \mathfrak{c}_0(\xi_1 - 3\varrho)^2 + \frac{5}{4\pi} \xi_2^2  - \frac{1}{2}\mathcal{C}_0^2 (\xi_1 - 3\varrho)^4 - \frac{1}{18}\xi_2^6.
\]
Shrinking $\varrho_0>0$ further if necessary, we may assume that  $4\mathcal{C}_0^2 \varrho_0^2 \leq \mathfrak{c}_0$ and $   4\pi{\mathfrak{c}_1^2\varrho_0^4} \leq  {45\mathcal{C}^2}    .$
Then, since $\varrho\leq \xi_1 < 4\varrho,$ we have $\lvert \xi_1 - 3\varrho\rvert \leq 2 \varrho,$ and hence 
\[
\frac{1}{2}\mathcal{C}_0^2 (\xi_1-3\varrho)^4 \leq  {2}\mathcal{C}_0^2   \varrho^2 (\xi_1 - 3\varrho)^2  \leq \frac{1}{2}\mathfrak{c}_0 (\xi_1 - 3\varrho)^2.
\]
Because of \eqref{eq:sdljxibound}, we also have
 $ \xi_2^6/18 \leq  \mathfrak{c}_1^2 \varrho^4  \xi_2^2/(18 \mathcal{C}^2) \leq  5  \xi_2^2/(8\pi).$
  Thus we   have 
  \[
  2\mathcal{Z}_{c,\varrho}(\tilde{K}_c) \geq \frac{1}{2}(P_1^2 +P_2^2) +\frac{1}{2}\mathfrak{c}_0 (\xi_1-3\varrho)^2 + \frac{5}{8\pi}\xi_2^2.
  \]
The right-hand side    vanishes only when $P_1=P_2=\xi_2=0 $ and  $\xi_1=3\varrho.$ We then obtain that
\begin{align*}
0&=P_1=\eta_1 +A_2(0) =\eta_1,\\
 0&=P_2=\eta_2 - A_1(x_*(c))+A_1(x_*(c)+3\varrho) = \eta_2 - \Delta_{c,\varrho},
\end{align*}
implying that    $\zeta = \zeta_{c,\varrho} .$ However,
\[
\tilde{K}_c(\zeta_{c,\varrho} )=\tilde{V}_{c,1}(x_*(c)+3\varrho)+V_{c,2}(0) = -\frac{c}{4}-d_{c,\varrho}+\frac{c}{4}=-d_{c,\varrho}<0,
\]
implying that $\zeta_{c,\varrho}  \notin \tilde{K}_c^{-1}(0).$ This completes the proof.
\end{proof}

The next lemma tells us that the set $\mathcal{U}_{c,\varrho}$ actually lies in the positive side of the function $\alpha_c.$

\begin{lemma}\label{lem:newlame1}
After decreasing $\varrho_0>0$ if necessary, we find a constant $\mathfrak{m}_0>0,$ independent of $\rho,c$ and $\varrho,$ satisfying
\[
\alpha_c( \zeta)  \geq \mathfrak{m}_0\varrho, \quad  \; \forall \zeta \in \mathcal{U}_{c,\varrho},
\]
provided $0<\varrho<\varrho_0$ and $c_2<c \leq c_2+\varrho^3.$
\end{lemma}
 \begin{proof} 
Recall from Section \ref{sec:firstin} that  $\alpha_{ c}(\zeta)  = \mathfrak{a}_c \xi_1 - \mathfrak{b}_c(\zeta),$
  where 
  \[
  \mathfrak{a}_c := \mu_c + \frac{1}{4}(\delta_c-1)>0, \qquad \mathfrak{b}_c(\zeta) := \eta_2 +\frac{\delta_c}{4}\xi_1 
  \]
  with $\delta_c = \cosh ( 2x_*)  >1$ and $\mu_{ c} =  \frac{ 4 \lambda_2 (c)    }{\delta_c - 1 }=\frac{4  V_{c,2}''(0)}{\delta_c-1}>0.$ By the Taylor theorem
\[
A_1( x_*(c)+ \xi_1) - A_1 ( x_*(c))  = \frac{\delta_c}{4} \xi_1 + O(\xi_1^2).
\]
uniformly for $c$ close to $c_2$ and $\lvert \xi_1 \rvert$ small enough. Then we can write  
\begin{equation}\label{eq:SPDF]}
 P_2 = \eta_2 - A_1(x_*(c))+A_1( x_*(c)+\xi_1)  = \mathfrak{b}_c (\zeta)+O(\xi_1^2).
\end{equation}
See \eqref{eq:P1P2}.
From the identity $\tilde{K}_c=0,$ we have
\begin{align*}
    P_2^2&= - P_1^2 -2\tilde{V}_{c,1}(x_*(c)+\xi_1) - 2V_{c,2}(\xi_2)\\
    &\leq -2V_{c,1}(x_*(c)+\xi_1) -2V_{c,2}(0)\\
    &=-2\left( V_{c,1}(x_*(c))+\frac{1}{2}\lambda_1(c)\xi_1^2 +O(\lvert \xi_1\rvert^3) + \frac{c}{4}\right) \\
    &=-2 B_*(c)  -\lambda_1(c)\xi_1^2 +O(  \xi_1 ^3)   ,
\end{align*}
where $\lambda_1(c)=V_{c,1}''(x_*(c))<0. $ Since by the mean value theorem $\lvert B_*(c)\rvert \leq C(c-c_2) \leq C\varrho^3$ for some $C>0,$  combining this with \eqref{eq:SPDF]} yields
\[
\lvert \mathfrak{b}_c(\zeta) \rvert \leq \sqrt{-\lambda_1(c)}\lvert \xi_1\rvert + \tilde{C}( \varrho^{3/2} +  \xi_1 ^{3/2} +\xi_1^2)
\]
for some $\tilde{C}>0.$ The argument in Remark \ref{rmk:separation} tells us that $\mathfrak{a}_{c_2} > \sqrt{-\lambda_1(c_2)},$ and hence by continuity we find a constant $\mathfrak{m}_0>0,$ independent of $\rho,c$ and $\varrho,$ such that
\[
\mathfrak{a}_c - \sqrt{-\lambda_1(c)}> 2 \mathfrak{m}_0.
\]
Since $\xi_1 \in [\varrho, 4\varrho),$ we find 
\[
\varrho^{3/2} = \sqrt{\varrho} \varrho \leq \sqrt{\varrho} \xi_1, \quad \;   \xi_1  ^{3/2} = \sqrt{\xi_1}\xi_1 <   2\sqrt{\varrho}\xi_1, \; \quad \xi_1^2  < 4 \varrho \xi_1,
\]
and hence 
\begin{align*}
    \alpha_c(\zeta)
     \geq  ( \mathfrak{a}_c   - \sqrt{-\lambda_1(c)})\xi_1- \tilde{C}( 3 + 4 \sqrt{\varrho}) \sqrt{\varrho}\xi_1 >2\mathfrak{m}_0 \xi_1- \tilde{C}( 3 + 4 \sqrt{\varrho}) \sqrt{\varrho}\xi_1 .
\end{align*}
Shrinking $\varrho_0>0$ if necessary, we may assume that $\tilde{C}( 3 + 4 \sqrt{\varrho_0}) \sqrt{\varrho_0}\leq \mathfrak{m}_0,$ so that 
\[
\alpha_c(\zeta) \geq 2\mathfrak{m}_0 \xi_1 - \mathfrak{m}_0 \xi_1 = \mathfrak{m}_0 \xi_1 \geq \mathfrak{m}_0\varrho, \quad \forall \zeta \in \mathcal{U}_{c,\varrho} .
\]
This completes the proof.
 \end{proof}

Pick    $ \mathfrak{n}_0 >1 $ so large that  $2/\mathfrak{n}_0 < \mathfrak{m}_0.$
Let $\mathfrak{n}\geq \mathfrak{n}_0.$ Note that 
\begin{equation}\label{eq:choiceofn}
\frac{2}{\mathfrak{n} }<\mathfrak{m}_0.
\end{equation}
Below, the second interpolation will be carried out in the  strip 
\[
\mathfrak{S}_{c,\varrho, \mathfrak{n}}:= \tilde{K}_c^{-1}(0) \cap  \mathcal{B}_{r_*}( z_*(c)) \cap \left\{\frac{\varrho}{\mathfrak{n}}\leq \alpha_c \leq \frac{2\varrho}{\mathfrak{n}}\right\} \cap \left\{  \lvert \xi_2 \rvert< \frac{\pi}{4}\right\}.
\]
Every $\zeta \in \mathfrak{S}_{c,\varrho, \mathfrak{n}}$ satisfies $\xi_1 < \varrho.$ Indeed, otherwise $\zeta \in \mathcal{U}_{c,\varrho},$ a contradiction to  Lemma \ref{lem:newlame1} and  the choice of $\mathfrak{n}. $  Furthermore, $\tilde{Y}_{\rm loc} =X_c$ on this strip and  $
  \mathfrak{S}_{c,\varrho, \mathfrak{n}} \cap \mathcal{U}_{c,\varrho}  =\emptyset.$
Moreover,  since  the capping $\psi$ vanishes for $\xi_1 \leq 2 \varrho,$ we see that 
\[
\tilde{K}_c = K_c \quad \text{  on } \mathfrak{S}_{c,\varrho, \mathfrak{n}}.
\]

 Set 
 \[
 \tilde{\mathfrak{S}}_{c,\varrho, \mathfrak{n}}:= \tilde{K}_c^{-1}(0) \cap  \mathcal{B}_{r_*}( z_*(c)) \cap \left\{\frac{\varrho}{\mathfrak{n}}\leq \alpha_c  \right\} \cap \left\{ \xi_1<\varrho,  \lvert \xi_2 \rvert< \frac{\pi}{4}\right\} .
 \]
 Note that $\tilde{K}_c= K_c$ on $\tilde{\mathfrak{S}}_{c,\varrho, \mathfrak{n}}$ and $ {\mathfrak{S}}_{c,\varrho, \mathfrak{n}} \subset \tilde{\mathfrak{S}}_{c,\varrho, \mathfrak{n}}.$
As in  Lemma \ref{lem:transtionregion}, the following statement provides   uniform bounds which are used later.    Its proof will be included in the appendix.

\begin{lemma}\label{lem:transtionregion2}

There exist positive constants $\mathfrak{E}_0, \mathfrak{R}_0$ and $\mathfrak{r}_0$, independent of $\mathfrak{n},\rho, c$ and $\varrho,$ such that the following holds. For each fixed $0<\rho<\rho_1$ and  every $\mathfrak{n} \geq \mathfrak{n}_0,$ 
after decreasing $\varrho_0=\varrho_0(\rho,\mathfrak{n})>0$ if necessary,   if   $0<\varrho<\varrho_0$ and $c_2<c \leq c_2+  \varrho^3,$  
\[
\lvert \zeta \rvert \leq \mathfrak{E}_0  {\varrho} \;\; \text{ and } \;\;  \xi_1 \geq \mathfrak{r}_0 \frac{\varrho}{\mathfrak{n}}, \qquad \forall \zeta \in \tilde{\mathfrak{S}}_{c,\varrho, \mathfrak{n}}.
\]
  Moreover, if $\zeta \in \mathfrak{S}_{c,\varrho, \mathfrak{n}},$ then we have a stronger bound
\[
\lvert \zeta \rvert \leq \mathfrak{R}_0 \frac{\varrho}{\mathfrak{n}}. 
\]
In particular,
\[
\mathfrak{r}_0 \frac{\varrho}{\mathfrak{n}} \leq \xi_1 \leq \lvert \zeta \rvert \leq \mathfrak{R}_0 \frac{\varrho}{\mathfrak{n}} , \qquad \forall \zeta \in {\mathfrak{S}}_{c,\varrho, \mathfrak{n}}.
\]
\end{lemma}

As mentioned earlier, the second interpolation will be performed in  $\mathfrak{S}_{c,\varrho, \mathfrak{n}},$ where $\tilde{Y}_{\rm loc}=X_c$ and $\tilde{K}_c = K_c.$ By Lemma \ref{lem:Weinstein}, the Liouville vector field $X_c$ is positively transverse to $\tilde{K}_c^{-1}(0)$ in this strip.  
So far, by Proposition \ref{prop:radislve}, we have positive transversality of $\mathcal{Z}_{c,\varrho}$ only on $\mathcal{U}_{c,\varrho}. $  However, the transition strip $\mathfrak{S}_{c,\varrho,\mathfrak{n}}$ is disjoint from $\mathcal{U}_{c,\varrho}.$ Therefore, we first establish the following additional lower bound.

 \begin{lemma}\label{lem:zoalsosmallstrip}
     With the same setup as in Lemma \ref{lem:transtionregion2}, 
       there exists $\mathfrak{e}_0>0,$ independent of $\rho, \mathfrak{n},c$ and $\varrho,$ such that for each fixed $0<\rho<\rho_1$ and every $\mathfrak{n} \geq \mathfrak{n}_0,$  
     after decreasing $\varrho_0=\varrho_0(\rho,\mathfrak{n})>0$  if necessary, we have
     \[
     \mathcal{Z}_{c,\varrho} (\tilde{K}_c)> \frac{\mathfrak{e}_0}{\mathfrak{n}} \varrho^2 \quad \text{ on } \;\tilde{\mathfrak{S}}_{c,\varrho, \mathfrak{n}}.
     \]
 \end{lemma}
\begin{proof} Recall that $\tilde{K}_c = K_c$ on $\tilde{\mathfrak{S}}_{c,\varrho,\mathfrak{n}}.$ Arguing as in the proof of  Proposition  \ref{prop:radislve}, we find a constant $C>0$ such that
\[
2\mathcal{Z}_{c,\varrho}(K_c) \geq (\xi_1 - 3\varrho)V_{c,1}' (x_*(c)+\xi_1) - C(\xi_1 - 3 \varrho)^4  - C \xi_2^6.  
\]
By Lemmas \ref{lem:cappingfunction} and  \ref{lem:transtionregion2}  we estimate
\begin{align*}
     (\xi_1 - 3\varrho) V_{c,1}'( x_*(c)+\xi_1)  &\geq -(\xi_1 - 3\varrho)a_0 \xi_1 > \frac{ 2 a_0 \mathfrak{r}_0}{\mathfrak{n}}\varrho^2\\
    C \left( (\xi_1 - 3\varrho)^4 + \xi_2^6 \right)  &\leq  \tilde{C}(  \varrho ^4 +   \varrho ^6 ), \quad \tilde{C}:= C( 81 +\mathfrak{E}_0^6).
\end{align*}
After decreasing $\varrho_0=\varrho_0(\rho, \mathfrak{n})>0$   if necessary, we may assume that
\[
 \tilde{C}(  \varrho ^4 +   \varrho ^6 ) < \frac{   a_0 \mathfrak{r}_0}{\mathfrak{n}}\varrho^2, \quad \; \varrho \in (0, \varrho_0).
\]
 Then the desired lower bound holds with $\mathfrak{e}_0:=a_0 \mathfrak{r}_0/2.$ This completes the proof.
\end{proof}


    We now perform the second interpolation between   the radial vector field $\mathcal{Z}_{c,\varrho} $ and the previously constructed Liouville vector field $\tilde{Y}_{\rm loc}=X_c$  in the   strip $\mathfrak{S}_{c,\varrho,\mathfrak{n}}.$  

\begin{figure}[t]
    \centering
    \begin{tikzpicture}
        \draw (5,0) to[out=90,in=10,looseness=1] (4,1.1);
        \draw (5,0) to[out=270,in=350,looseness=1] (4,-1.1);
        \draw (4,1.1) to[out=190,in=350,looseness=1] (2,1);
        \draw (4,-1.1) to[out=170,in=10,looseness=1] (2,-1);
        \draw (2,1) to[out=170,in=320] (1,1.5);
        \draw (2,-1) to[out=190,in=40] (1,-1.5);
        \draw (1,1.5) to[out=140,in=315] (0,2.5);
        \draw (1,-1.5) to[out=220,in=45] (0,-2.5);
        \fill[lightgray!20] (2,1) to[out=170,in=320] (1,1.5) -- (1,-1.5) to[out=40,in=190] (2,-1) -- (2,1);
        \fill[yellow!15] (5,0) to [out=90,in=10,looseness=1] (4,1.1) -- (4,-1.1) to[out=350,in=270,looseness=1] (5,0);
        \begin{scope}
        \clip (3,2) -- (4.1,2) -- (4.1,-2) -- (3,-2) -- (3,2);
        \fill[red!15] (4,1.1) to[out=190,in=350,looseness=1] (2,1) -- (2,-1) to[out=10,in=170,looseness=1] (4,-1.1) -- (4,-1.1);
        \end{scope}
        \begin{scope}
        \clip (2,2) -- (3,2) -- (3,-2) -- (2,-2) -- (2,2);
        \fill[green!15] (4,1.1) to[out=190,in=350,looseness=1] (2,1) -- (2,-1) to[out=10,in=170,looseness=1] (4,-1.1) -- (4,-1.1);
        \end{scope}
        \fill[lightgray!20] (1,1.5) -- (1,-1.5) to[out=220,in=45] (0,-2.5) -- (0,2.5) to[out=315,in=140] (1,1.5);
        \begin{scope}
        \clip (2.5,0) circle (1.3cm);
        \fill[orange!15] (1,-1.5) to[out=40,in=190] (2,-1) -- (2,1) to[out=170,in=320] (1,1.5) -- (1,-1.5) ;
        \end{scope}
        \draw[orange!70!black] (1.6,-0.3) node {$\mathcal{S}_{c,\rho}$};
        \draw[green!50!black] (2.5,-2) node {$-\rho<\alpha_c<\frac{\varrho}{\mathfrak{n}}$};
        \draw (2.5,-1.7) --(2.7,-0.8);
        \draw[red!60!black] (3.5,2) node {$\frac{\varrho}{\mathfrak{n}}\le \alpha_c\le\frac{2\varrho}{\mathfrak{n}}$};
        \draw (3.5,1.7) --(3.2,0.8);
        \draw[yellow!40!black] (6.2,-2.2) node {$\frac{2\varrho}{\mathfrak{n}}<\alpha_c$};
        \draw (6,-1.7) --(4.5,-0.8);
        \draw (2.5,0) node {$X_c$}; 
        \draw (3.5,-0.3) node {$\bar{Y}_{\mathrm{loc}}$};
        \draw (4.5,0.3) node {$\mathcal{Z}_{c,\varrho}$};
        \draw (1.7,0.3) node {$\tilde{Y}_\mathrm{loc}$};
        \draw (0.5,0.5) node {$\hat{Y}$};
        \draw[red!60!black] (3.5,0.5) node {$\mathfrak{S}_{c,\varrho,\mathfrak{n}}$};
    \end{tikzpicture}
    \caption{A sketch of the different Liouville vector fields in the neck region and their corresponding domains on the capped energy hypersurface}
    \label{fig:secondinterpolation}
\end{figure}
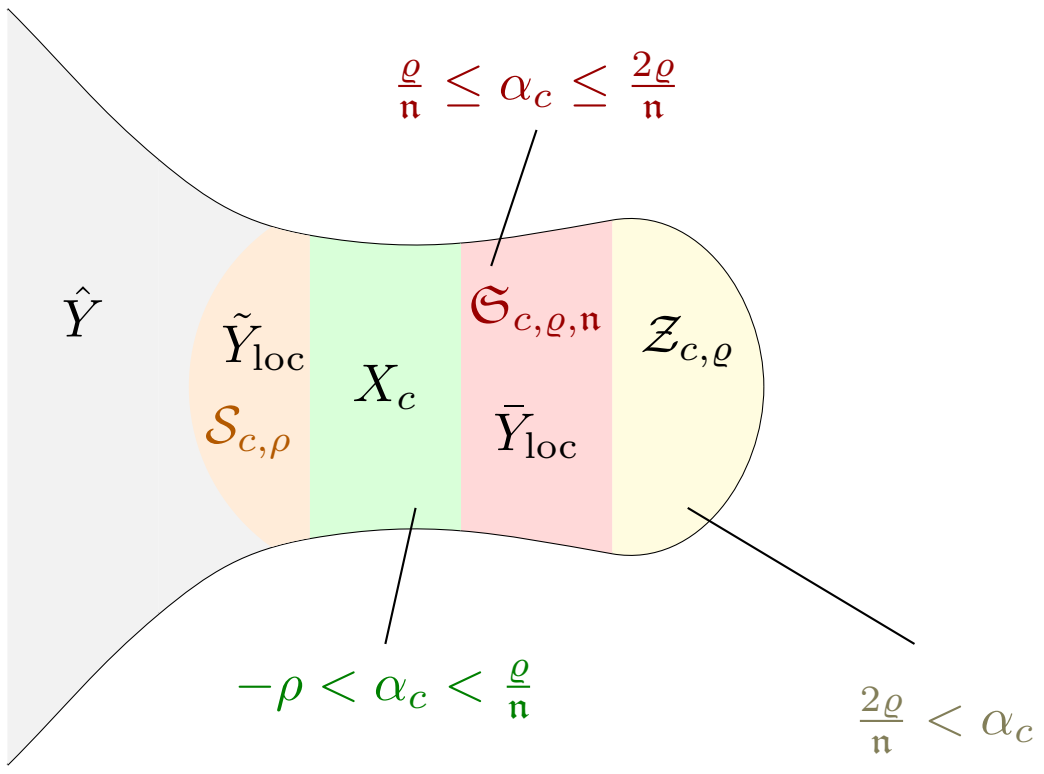

A direct computation gives $  \mathcal{Z}_{c,\varrho}  -X_c= X_{\Omega_{c,\varrho}},$ where
\begin{equation}\label{eq:secondomnega}
 \Omega_{c,\varrho}(\zeta)  = (  \frac{1}{2} -a) \xi_1 \eta_1 +( \frac{1}{2}-b)\xi_2 \eta_2 + \sigma_c \xi_1 \xi_2 -\frac{3}{2}\varrho \eta_1 + \frac{1}{2}\Delta_{c,\varrho}\xi_2.
\end{equation}
Pick a smooth non-decreasing function $g \colon \R \to [0,1]$ such that $g(t)=0$ for $t \leq 1,$ $g(t)=1$ for $t \geq 2$ and $0 \leq g'(t) \leq 2 $ and then  define
\[
g_{c,\varrho,\mathfrak{n}}(\zeta):=g \left(   \frac{\mathfrak{n}\alpha_c(\zeta)}{\varrho}      \right).
\]
We set
\[
\bar{Y}_{\rm loc}:=\tilde{Y}_{\rm loc}+X_{g_{c,\varrho,\mathfrak{n}}\Omega_{c,\varrho}}.
\]
See Figure \ref{fig:secondinterpolation}. Note  that $\bar{Y}_{\rm loc} = \tilde{Y}_{\rm loc} $ if $\alpha_c \leq \varrho/\mathfrak{n}$ and $\bar{Y}_{\rm loc} = \tilde{Y}_{\rm loc}+X_{\Omega_{c,\varrho}}=\mathcal{Z}_{c,\varrho}$ if $\alpha_c \geq 2\varrho /\mathfrak{n}.$

 \begin{lemma}\label{lem:latinterpolation}
There exists $\mathfrak{n}_1 \geq \mathfrak{n}_0,$ independent of $\rho,$ such that for each fixed $\rho \in (0,\rho_1)$ and every  $\mathfrak{n} \geq \mathfrak{n}_1,$ after decreasing  $\varrho_0=\varrho_0(\rho,\mathfrak{n})>0$ if necessary, the following holds. For every $ \varrho \in (0,\varrho_0)$ and every $c_2<c \leq c_2+  \varrho^3,$
the local Liouville vector field $\bar{Y}_{\rm loc}$ is positively transverse to $\tilde{K}_c^{-1}(0) \cap \mathcal{B}_{r_*}(z_*(c)) \cap \{ \lvert \xi_2 \rvert < \pi/4\}.$ 
 \end{lemma}
\begin{proof}
We first assume that $\alpha_c \leq \varrho / \mathfrak{n}.$ If   $\zeta$ satisfies $\xi_1 \geq\varrho,$ then $\zeta \in \mathcal{U}_{c,\varrho}. $ This contradicts Lemma \ref{lem:newlame1} and the choice of $\mathfrak{n},$ see \eqref{eq:choiceofn}. Hence, $\xi_1<\varrho,$ so   the capping   $\psi$ vanishes, implying that $\tilde{K}_c =K_c.$ 
Then the assertion follows by Lemma~\ref{lem:firstinterpolatoinsuccess}.

Suppose that  $\alpha_c \geq 2\varrho/ \mathfrak{n},$ so that $\bar{Y}_{\rm loc}=\mathcal{Z}_{c,\varrho}.$ If $\xi_1 \geq \varrho,$ then the assertion follows by Proposition \ref{prop:radislve}. If $\xi_1 <\varrho,$ then  the assertion follows by Lemma \ref{lem:zoalsosmallstrip}.

It remains to consider the case $ \varrho/ \mathfrak{n} \leq \alpha_c \leq 2\varrho/\mathfrak{n}.$ This region corresponds precisely to the strip $\mathfrak{S}_{c,\varrho,\mathfrak{n}}.$ As before, we compute
\[
\bar{Y}_{\rm loc}(\tilde{K}_c) = \bar{Y}_{\rm loc}(K_c)= (1-g_{c,\varrho, \mathfrak{n}})X_c(K_c)+ g_{c,\varrho, \mathfrak{n}} \mathcal{Z}_{c,\varrho}(K_c)+\Omega_{c,\varrho} X_{g_{c,\varrho, \mathfrak{n}}}(K_c).
\]
By Lemmas \ref{lem:Weinstein},  \ref{lem:newlame1}, \ref{lem:transtionregion2} and  \ref{lem:zoalsosmallstrip} we have
     \begin{equation}\label{eq:finsdl?}
\bar{Y}_{\rm loc}(\tilde{K}_c)  \geq \left( (1-g_{c,\varrho, \mathfrak{n}}) \frac{\mathfrak{w}\mathfrak{r}_0^2}{\mathfrak{n}^2} + g_{c,\varrho, \mathfrak{n}} \frac{\mathfrak{e}_0}{\mathfrak{n}} \right)\varrho^2+\Omega_{c,\varrho} X_{g_{c,\varrho, \mathfrak{n}}}(K_c).    
\end{equation}
To estimate the cutoff term, we first compute
\[
\Omega_{c,\varrho} X_{g_{c,\varrho,\mathfrak{n}}}(K_c)=\Omega_{c,\varrho}\frac{\mathfrak{n}}{\varrho}g' X_{\alpha_c}(K_c).
\]
Recall from the first interpolation step that $X_{\alpha_c}(K_c)=\nu_c E_{\rm lin}+R_{\alpha_c,c},$
where $E_{\rm lin}$ and $\nu_c$ are  as in \eqref{eq:Elin} and \eqref{eq:nuc}, respectively, and $R_{\alpha_c,c}=X_{\alpha_c}(\mathcal{R}_c)$ satisfies $\lvert R_{\alpha_c, c}(\zeta)\rvert \leq C_R \lvert \zeta \rvert^2$ for some $C_R>0,$ independent of $\mathfrak{n},c$ and $\varrho.$ To make use of the previous estimates, we set 
\begin{align*}
\Omega_{c,\varrho}^{\rm err}&=\Omega_{c,\varrho} - \frac{3\varrho}{\tanh 2x_*(c)}E_{\rm lin}\\
&= (  \frac{1}{2} -a) \xi_1 \eta_1 +( \frac{1}{2}-b)\xi_2 \eta_2 + \sigma_c \xi_1 \xi_2   + ( \frac{1}{2}\Delta_{c,\varrho}+\frac{3}{8}\delta_c \varrho)\xi_2.
\end{align*}
The first three terms are quadratic in $\zeta,$ and their coefficients are uniformly bounded. Thus, we find $D_0>0,$ independent of $\rho, \mathfrak{n},c$ and $\varrho,$ such that
\[
\lvert (  \frac{1}{2} -a) \xi_1 \eta_1 +( \frac{1}{2}-b)\xi_2 \eta_2 + \sigma_c \xi_1 \xi_2  \rvert \leq D_0 \lvert \zeta \rvert^2.
\]
For the last term, we use Taylor's theorem and obtain $\Delta_{c,\varrho} = -\frac{3}{4}\delta_c \varrho+ O(\varrho^2),$ so 
\[
\lvert ( \frac{1}{2}\Delta_{c,\varrho}+\frac{3}{8}\delta_c \varrho)\xi_2 \rvert \leq D_1 \varrho^2 \lvert \xi_2 \rvert \leq D_1 \varrho^2 \lvert \zeta\rvert
\]
for some $D_1>0,$ independent of $\rho, \mathfrak{n},c$ and $\varrho.$ Consequently,
\begin{equation}\label{eq:omegaerrdesti}
\lvert \Omega_{c,\varrho}^{\rm err}\rvert \leq D_{\rm err}( \lvert \zeta\rvert^2 +\varrho^2 \lvert \zeta \rvert ),  \;\; \quad D_{\rm err}:=\max\{ D_0,D_1\}.
\end{equation}

We find
\begin{align*}
\Omega_{c,\varrho}X_{\alpha_c}(K_c) &=  \frac{3\varrho}{\tanh 2x_*(c)} E_{\rm lin}R_{\alpha_c,c} + \nu_c E_{\rm lin}\Omega_{c,\varrho}^{\rm err} \\
& \quad + \Omega_{c,\varrho}^{\rm err}R_{\alpha_c,c}+\frac{3 \varrho \nu_c}{\tanh 2x_*(c)} E_{\rm lin}^2.
\end{align*}
We estimate each term separately.
For this, we set
\[
A_{\ell}:=\sup_{  I_{\ell}} \frac{1}{\lvert \tanh 2x_*(c) \rvert}, \quad \nu_{\ell}:=\sup_{  I_{\ell}}  \lvert \nu_c \rvert, \quad \beta_{\ell}:=\inf_{  I_{\ell}} \bigg| \frac{\nu_c}{\tanh 2x_*(c)}\bigg| .
\]
 
We consider the first term. 
We observe that
\[
\bigg| \frac{3\varrho}{\tanh 2x_*(c)} E_{\rm lin}R_{\alpha_c,c}  \bigg| \leq 3A_{\ell} \varrho \lvert E_{\rm lin}\rvert C_R \lvert \zeta\rvert^2  .
\]
An application of Young's inequality $ab \leq \varepsilon a^2 +  b^2/({4\varepsilon})$
with $a = \sqrt{\varrho}\lvert E_{\rm lin}\rvert,$ $b=3 A_{\ell}C_R\sqrt{\varrho}\lvert \zeta \rvert^2$ and $\varepsilon = \beta_{\ell}/4,$ we obtain
\begin{align*}
\bigg| \frac{3\varrho}{\tanh 2x_*(c)} E_{\rm lin}R_{\alpha_c,c}  \bigg|  &\leq   \frac{\beta_{\ell}}{4}  \varrho E_{\rm lin}^2 + \frac{ 9 A_{\ell}^2 C_R^2}{\beta_{\ell}}  \varrho\lvert \zeta \rvert^4  \\
&\leq   \frac{\beta_{\ell}}{4}  \varrho E_{\rm lin}^2 + \frac{ 9 A_{\ell}^2 C_R^2}{\beta_{\ell}}  \frac{\lvert \zeta \rvert^4}{\varrho}  ,
\end{align*}
 where the last inequality follows since $0<\varrho<\varrho_0<1.$

For the second term, we have using \eqref{eq:omegaerrdesti}
\[
\lvert \nu_c E_{\rm lin}\Omega_{c,\varrho}^{\rm err}\rvert  \leq \nu_{\ell} \lvert E_{\rm lin}\rvert D_{\rm err} (\lvert \zeta\rvert^2 + \varrho^2 \lvert \zeta \rvert ).
\]
Again by Young's inequality we have
\begin{align*}
\lvert \nu_c E_{\rm lin}\Omega_{c,\varrho}^{\rm err}\rvert   &\leq \frac{\beta_{\ell}}{4}\varrho E_{\rm lin}^2 + \frac{\nu_{\ell}^2 D_{\rm err}^2}{\beta_{\ell}\varrho} ( \lvert \zeta \rvert^2 + \varrho^2 \lvert \zeta \rvert )^2\\
&\leq \frac{\beta_{\ell}}{4}\varrho E_{\rm lin}^2 + \frac{\nu_{\ell}^2 D_{\rm err}^2}{\beta_{\ell}\varrho} ( 2\lvert \zeta \rvert^4 + 2\varrho^4 \lvert \zeta \rvert^2  ).
\end{align*}

We next find the following estimate for the third term:
\begin{align*}
    \lvert \Omega_{c,\varrho}^{\rm err} R_{\alpha_c,c}\rvert &\leq  D_{\rm err}(\lvert \zeta \rvert^2 + \varrho^2 \lvert \zeta \rvert) C_R \lvert \zeta \rvert^2 \\
    &= C_R D_{\rm err} (\lvert \zeta \rvert^4 + \varrho^2 \lvert \zeta \rvert^3)\\
    & \leq C_R D_{\rm err} (\frac{\lvert \zeta \rvert^4}{\varrho} + \varrho^2 \lvert \zeta \rvert^3),
\end{align*}
where the last inequality holds since $0<\varrho<\varrho_0<1.$

Finally, for the last term we find
\[
\frac{3 \varrho \nu_c}{\tanh 2x_*(c)} E_{\rm lin}^2 \geq  3 \varrho \beta_{\ell}  E_{\rm lin}^2.
\]

Consequently, we have
\begin{align*}
    \Omega_{c,\varrho}X_{\alpha_c}(K_c)&\geq  3 \varrho \beta_{\ell}  E_{\rm lin}^2 - \frac{\beta_{\ell}}{4}  \varrho E_{\rm lin}^2 - \frac{ 9 A_{\ell}^2 C_R^2}{\beta_{\ell}}  \frac{\lvert \zeta \rvert^4}{\varrho}- \frac{\beta_{\ell}}{4}\varrho E_{\rm lin}^2 \\
    &\quad - \frac{2\nu_{\ell}^2 D_{\rm err}^2}{\beta_{\ell}} (  \frac{\lvert \zeta \rvert^4}{\varrho} + \varrho^3 \lvert \zeta \rvert^2  )-C_R D_{\rm err} (\frac{\lvert \zeta \rvert^4}{\varrho} + \varrho^2 \lvert \zeta \rvert^3)\\
    &\geq   - \mathfrak{Q}_0 ( \frac{\lvert \zeta\rvert^4}{\varrho} + \varrho^3 \lvert \zeta\rvert^2 + \varrho^2 \lvert \zeta \rvert^3),
\end{align*}
where  
\[
\mathfrak{Q}_0:= \frac{ 9 A_{\ell}^2 C_R^2}{\beta_{\ell}} + \frac{2\nu_{\ell}^2 D_{\rm err}^2}{\beta_{\ell}} + C_R D_{\rm err}.
\]
Fix $0 <\rho<\rho_1$ and $\mathfrak{n} \geq \max\{1, \mathfrak{n}_0\},$ and decrease $\varrho_0 = \varrho_0(\rho,\mathfrak{n})$ so that $\varrho_0<\mathfrak{n}^{-2}.$  Lemma \ref{lem:transtionregion2} gives  $\lvert \zeta \rvert \leq \mathfrak{R}_0 \varrho/\mathfrak{n} $ on $\mathfrak{S}_{c,\varrho,\mathfrak{n}}.$ Since $\varrho^2 \mathfrak{n}^2\leq1$ and $\varrho^2 \mathfrak{n} \leq 1,$ we obtain 
\begin{align*}
    \frac{\lvert \zeta\rvert^4}{\varrho} + \varrho^3 \lvert \zeta\rvert^2 + \varrho^2 \lvert \zeta \rvert^3 & \leq \frac{\varrho^3}{\mathfrak{n}^4}(\mathfrak{R}_0^4 + \mathfrak{R}_0^2 \varrho^2 \mathfrak{n}^2 + \mathfrak{R}_0^3 \varrho^2 \mathfrak{n})\\
    &\leq ( \mathfrak{R}_0^4 + \mathfrak{R}_0^2 + \mathfrak{R}_0^3) \frac{\varrho^3}{\mathfrak{n}^4}.
\end{align*}
It follows that 
\[
\Omega_{c,\varrho}X_{\alpha_c}(K_c) \geq - \mathfrak{Q}_1 \frac{\varrho^3}{\mathfrak{n}^4},  
\]
where $ \mathfrak{Q}_1:= \mathfrak{Q}_0 ( \mathfrak{R}_0^4 + \mathfrak{R}_0^2 +\mathfrak{R}_0^3) $ is independent of $\rho, \mathfrak{n},c$ and $\varrho.$ It follows that
\[
\Omega_{c,\varrho}X_{g_{c,\varrho,\mathfrak{n}}}(K_c) =  \frac{\mathfrak{n}}{\varrho} g' \Omega_{c,\varrho}X_{\alpha_c}(K_c) \geq - 2 \mathfrak{Q}_1 \frac{\varrho^2}{\mathfrak{n}^3} .
\]
Plugging this into \eqref{eq:finsdl?} yields
\[
\bar{Y}_{\rm loc}(\tilde{K}_c)  \geq \left( (1-g_{c,\varrho,\mathfrak{n}}) \frac{\mathfrak{w}\mathfrak{r}_0^2}{\mathfrak{n}^2} + g_{c,\varrho,\mathfrak{n}} \frac{\mathfrak{e}_0}{\mathfrak{n}} - 2 \mathfrak{Q}_1 \frac{1}{\mathfrak{n}^3}\right)\varrho^2. 
\]
Since $\mathfrak{Q}_1, \mathfrak{w}, \mathfrak{r}_0$ and $\mathfrak{e}_0$ are independent of $\rho$ and   $\mathfrak{n},$ choose $\mathfrak{n}_1 \geq \mathfrak{n}_0,$ independent of $\rho,$ such that   $4\mathfrak{Q}_1 \leq   \mathfrak{w} \mathfrak{r}_0^2 \mathfrak{n} \leq \mathfrak{e}_0 \mathfrak{n}^2 $ for every $\mathfrak{n} \geq \mathfrak{n}_1.$ For each such $\mathfrak{n}$ and each fixed $0<\rho<\rho_1,$  we then decrease $\varrho_0=\varrho_0(\rho,\mathfrak{n})>0$ further if necessary, so that all the preceding conclusions hold.  It follows that 
\[
\bar{Y}_{\rm loc}(\tilde{K}_c)  \geq \left( (1-g_{c,\varrho,\mathfrak{n}}) \frac{\mathfrak{w}\mathfrak{r}_0^2}{\mathfrak{n}^2} + g_{c,\varrho,\mathfrak{n}}  \frac{\mathfrak{w}\mathfrak{r}_0^2}{\mathfrak{n}^2}  -  \frac{\mathfrak{w}\mathfrak{r}_0^2}{2\mathfrak{n}^2} \right)\varrho^2 =  \frac{\mathfrak{w}\mathfrak{r}_0^2}{2\mathfrak{n}^2} \varrho^2>0. 
\]
This completes the proof.
\end{proof}


\subsection{Extension of the local Liouville vector field}\label{sec:globalextn}
So far the construction has been carried out in the local coordinates. 
We now extend the construction to obtain a transverse Liouville vector field defined on the whole phase space.

We first collect   the parameters chosen so that the preceding conclusions apply. Let $\hat{Y}$ be the Liouville vector field as in Remark \ref{remark:slightlyaboveL2}.  Take the closed interval $I_*=[c_2, c_2+\varepsilon_*]$ as in Section \ref{sec:firstin} and choose $\rho_0>0$ and $r_*>0$ as in Lemma \ref{lem:firstinterpolatoinsuccess}. Then  for every $0<\rho<\rho_0$ and $c_2 <c<c_2+\rho^3,$ the local Liouville vector field $\tilde{Y}_{\rm loc}$ is positively transverse to $K_c^{-1}(0)$ in $\mathcal{B}_{r_*}(z_*(c)).$ After decreasing $r_*>0$ further if necessary, we may assume that the inequalities in \eqref{eq:SDFK} hold. Let $\rho_1 \in (0,\rho_0)$ be such that \eqref{eq:Qminunmd} and \eqref{eq:cdsd02} are valid. 
By   Lemma \ref{lem:locallemman}, there exists    $\delta_0 \in ( 0, \min_{c \in I_*} x_*(c))$ such that for every $0<\rho<\rho_1$ and every $c_2 <c<c_2+\rho^3,$ each point in $K_c^{-1}(0)$ with $\lvert x_1 - x_*(c)\rvert <\delta_0$ lies either in $\mathcal{B}_{r_*}(z_*(c))$ or in the corresponding coordinate $r_*$-ball centered at $(x_*(c), \pi, 0 ,-A_1(x_*(c))).$ Then we find $\varepsilon_{\ell} \in (0,\varepsilon_*) $ sufficiently small so that Lemma \ref{lem:additionalestimateforhaty} applies. We now take $\varrho_0$ sufficiently small so that for every $0<\varrho<\varrho_0$ and $c_2 < c \leq c_2+\varrho^3,$ the capping function in Lemma \ref{lem:cappingfunction} exists, and hence the perturbed Hamiltonian $\tilde{K}_c=\tilde{K}_{c,\varrho}$ is well-defined. The perturbed energy level $\tilde{\Sigma}_c=\tilde{K}_c^{-1}(0)$ is diffeomorphic to $S^1 \times S^2.$ Fix $ \rho \in (0,\rho_1) $ and $\mathfrak{n}_1>0$ sufficiently large as in  Lemma \ref{lem:latinterpolation}. Fix also $\mathfrak{n} \geq \mathfrak{n}_1$ and take $\varrho_0=\varrho_0(\rho,\mathfrak{n})$ sufficiently small, so that for every $0<\varrho<\varrho_0$ and $c_2 < c \leq c_2 +\varrho^3,$ the local Liouville vector field $\bar{Y}_{\rm loc}$ is positively transverse to $\tilde{K}_c^{-1}(0)  \cap \mathcal{B}_{r_*}(z_*(c)) \cap \{ \lvert \xi_2 \rvert < \pi/4\}.$

Set
\[
\mathcal{W}_{c,\rho,\varrho}  :=  \tilde{K}_c^{-1}(0) \cap \mathcal{B}_{r_*}( z_*(c)) \cap \{ \lvert \xi_2 \rvert  < \frac{\pi}{4}\} \cap \{ \alpha_c \geq - \frac{5}{2}\rho\}.
\]

 \begin{lemma}\label{lema:cloglbascne} Under the preceding choices, there exists $C_{\mathcal{W}}>0,$ independent of $c,\rho, \mathfrak{n}$ and $\varrho,$ such that 
 \[
 \lvert \zeta \rvert \leq C_{\mathcal{W}}\rho, \qquad \zeta \in \mathcal{W}_{c,\rho,\varrho}.
 \]
  
 \end{lemma}
\begin{proof}Pick $\zeta \in \mathcal{W}_{c,\rho,\varrho}.$  
 We first consider the case $\xi_1 \geq \varrho.$ By construction of $\tilde{K}_c$ we have $\xi_1 <4\varrho.$ Then  by Lemma \ref{lem:perturbedsall}  we have $\lvert \zeta \rvert \leq \mathfrak{C}\varrho < \mathfrak{C}\rho.$
 If $\xi_1<\varrho$ and $\alpha_c \geq \varrho/\mathfrak{n},$ then by Lemma \ref{lem:transtionregion2} we have $\lvert \zeta \rvert \leq \mathfrak{E}_0 \varrho<\mathfrak{E}_0\rho.$
Finally suppose that $\xi_1<\varrho$ and $-5\rho/2 \leq \alpha_c <\varrho / \mathfrak{n}.$   Since the capping function $\psi$ vanishes, $\tilde{K}_c=K_c$ at $\zeta.$ 
Recall from  the end of Section \ref{sec:firstin} that 
$\lvert Q_c(\zeta)\rvert+\alpha_c(\zeta)^2 \geq \frac{m}{2}\lvert \zeta\rvert^2 $
for some $m>0.$
See \eqref{eq:Qminunmd}. Using the Taylor expansion, with $C>0$ chosen uniformly for $c \in I_*,$ we have
\[
\lvert Q_c(\zeta)\rvert  \leq C (c-c_2) + C_{\rm Tay}\lvert \zeta \rvert^3 \leq C \varrho^3 +C_{\rm Tay}\lvert \zeta \rvert^3  <C \rho^3 +C_{\rm Tay}\lvert \zeta \rvert^3  .
\]
Since $\varrho<\rho$ and $\mathfrak{n}>1,$ we have $\lvert \alpha_c \rvert \leq  5\rho/2.$
Then by \eqref{eq:SDFK} we obtain 
\[
\frac{m}{4}\lvert \zeta \rvert^2 \leq C\rho^3 +\frac{25}{4}\rho^2 \leq C\rho^2 +\frac{25}{4}\rho^2 .
\]
Therefore, $\lvert \zeta \rvert \leq \mathfrak{C}_1 \rho$ with $\mathfrak{C}_1:= \sqrt{ (4C+25 )/m}.$
The choice $C_{\mathcal{W}}:=\max\{  \mathfrak{C}, \mathfrak{E}_0, \mathfrak{C}_1\}$ completes the proof.
 \end{proof}

Since $C_{\mathcal{W}}>0$ obtained in the preceding lemma is independent of $\rho$ and $\mathfrak{n},$ after decreasing $\rho_1>0$ if necessary, we may assume that
$C_{\mathcal{W}}\rho_1 < \min \{   {r_*}/{2},  {\pi}/{8}\}.$ 
This procedure does not affect any of the preceding conclusions. We now fix $0<\rho<\rho_1$ and $\mathfrak{n} \geq \mathfrak{n}_1 $ and choose $\varrho_0=\varrho_0(\rho, \mathfrak{n})>0$ as above. Then  the preceding lemma provides
\[
\mathcal{W}_{c,\rho,\varrho} \subset \mathcal{B}_{r_*/2}(z_*(c)) \cap \{ \lvert \xi_2 \rvert < \frac{\pi}{8}\} 
\]
for every $0<\varrho<\varrho_0$ and $c_2<c \leq c_2+\varrho^3.$

Set
\[
\mathcal{N}_{c,\rho,\varrho}({z_*(c)}):=\tilde{K}_c^{-1}(0) \cap\mathcal{B}_{r_*}(z_*(c))\cap \{ \lvert \xi_2 \rvert <\frac{\pi}{4} \} \cap \{ \alpha_c >-\frac{9}{4}\rho\}. 
\]
Note that $ \mathcal{N}_{c,\rho,\varrho}({z_*(c)})   \subset \mathcal{W}_{c,\rho,\varrho}$ and that on the complement $   \mathcal{W}_{c,\rho,\varrho} \setminus \mathcal{N}_{c,\rho,\varrho}({z_*(c)})$ we have 
\[
-\frac{5\rho}{2} \leq \alpha_c \leq -\frac{9\rho}{4} <-2\rho, \qquad \bar{Y}_{\rm loc}=\tilde{Y}_{\rm loc}=\hat{Y}=Y_1. 
\]

Let $\Lambda=\{ {\rm id}, \sigma , \tau, \tau \sigma \},$ where $\sigma:=\sigma_{\rm el}(-1,\cdot)$ and   $\tau(x_1, x_2, y_1, y_2)=(x_1, x_2+\pi, y_1, y_2).$
The maps $\sigma$ and $\tau$ are commuting symplectic involutions. The construction above implies that  $\tilde{K}_c$ and $\hat{Y}$ are $\Lambda$-invariant.

For each $\lambda \in \Lambda,$ we define
\[
\bar{Y}_{\rm loc}^{\lambda} := \lambda_* \bar{Y}_{\rm loc}.
\]
Note that $\bar{Y}_{\rm loc}^{\rm id} :=  \bar{Y}_{\rm loc}.$ We also set
\begin{align*}
\mathcal{W}_{c,\rho,\varrho}^{\rm global}&:=\bigcup_{\lambda \in \Lambda} \lambda( \mathcal{W}_{c,\rho,\varrho}) \\  \mathcal{M}_{c,\rho,\varrho}&:= \tilde{K}_c^{-1}(0) \setminus \bigcup_{\lambda \in \Lambda} \lambda( \mathcal{N}_{c,\rho,\varrho}(z_*(c))). 
\end{align*}
 Note that $  \mathcal{W}_{c,\rho,\varrho}^{\rm global} \cup \mathcal{M}_{c,\rho,\varrho}=\tilde{K}_c^{-1}(0).$ Moreover, the four subsets $\lambda(\mathcal{W}_{c,\rho,\varrho}), \lambda \in \Lambda,$ have pairwise disjoint closures. 
By the definition of $\mathcal{M}_{c,\rho,\varrho},$ we have
\[
 z\in \lambda( \mathcal{W}_{c,\rho,\varrho}) \cap \mathcal{M}_{c,\rho,\varrho} \quad \Rightarrow \quad \lambda^{-1}z \in \mathcal{W}_{c,\rho,\varrho} \setminus \mathcal{N}_{c,\rho,\varrho}(z_*(c)).
 \]
Therefore, using the $\Lambda$-invariance of $\hat{Y},$ for every $\lambda \in \Lambda$ we have
\[
\bar{Y}_{\rm loc}^{\lambda}=\hat{Y} \quad \text{ on } \lambda( \mathcal{W}_{c,\rho,\varrho}) \cap \mathcal{M}_{c,\rho,\varrho}.
\]
The intersection  $\lambda(\mathcal{W}_{c,\rho,\varrho}) \cap\mathcal{M}_{c,\rho,\varrho}$ is compact for  each $\lambda\in\Lambda$ and satisfies
\[
\lambda(\mathcal{W}_{c,\rho,\varrho})\cap\mathcal{M}_{c,\rho,\varrho}\subset
\lambda\left(\mathcal{B}_{r_*}(z_*(c)) \cap  \{ \lvert \xi_2 \rvert < \frac{\pi}{4} \} \cap \{\alpha_c<-2\rho\} \right).
\]
Therefore, there exist   a $\Lambda$-invariant open neighborhood $\mathcal{U}_{\mathcal M}$ of $\mathcal{M}_{c,\rho,\varrho}$ and  an open neighborhood $\mathcal{U}_{\mathcal W}$ of $\mathcal{W}_{c,\rho,\varrho}$ such that the four subsets
$\lambda(\mathcal U_{\mathcal W}) ,  \lambda\in\Lambda,$ are pairwise disjoint.
Moreover, they are chosen so that the following holds: the closure $\overline{\mathcal{U}_{\mathcal{W}}}$ is compact and contained in $\mathcal{B}_{r_*}(z_*(c)) \cap \{ \lvert \xi_2 \rvert <\pi/4\},$ and for every $\lambda \in \Lambda$ we have
\[
\mathcal{U}_{\mathcal{M}} \cap \lambda(\mathcal{U}_{\mathcal{W}}) \subset \lambda ( \mathcal{U}_{\mathcal{W}} \cap \{ \alpha_c< -2\rho\}).
\]
 Then for every $p \in \mathcal{U}_{\mathcal{M}} \cap \lambda(\mathcal{U}_{\mathcal{W}}),$ we have $\lambda^{-1}(p)\in \mathcal{U}_{\mathcal{W}} \cap \{ \alpha_c <-2\rho\}.$
In this region, the cutoff functions $\ell_c$ from the first interpolation step and $g_{c, \varrho, \mathfrak{n}}$ from the second interpolation step vanish identically. Using the $\Lambda$-invariance of $\hat{Y},$ we obtain 
\[
\bar{Y}_{\rm loc}^\lambda=\hat{Y} \qquad \text{ on } \; \mathcal{U}_{\mathcal{M}} \cap \lambda(\mathcal{U}_{\mathcal{W}}).
\]
We then define
\[
\bar{Y}:=  \begin{cases} \hat{Y} & \text{ on } \mathcal{U}_{\mathcal{M}}, \\
\bar{Y}_{\rm loc}^{\lambda} & \text{ on } \lambda(\mathcal{U}_{\mathcal{W}}), \;\; \lambda \in \Lambda .
\end{cases}
\]
This is a well-defined smooth $\Lambda$-invariant Liouville vector field defined on the open neighborhood $\mathcal{U}_{\mathcal{M}} \cup \bigcup_{\lambda \in \Lambda}\lambda(\mathcal{U}_{\mathcal{W}})$ of $\tilde{K}_{c}^{-1}(0).$

We now extend $\bar{Y}$ to the whole phase space. Recall that the local Liouville vector field $\bar{Y}_{\rm loc}$ is defined by $\bar{Y}_{\rm loc} = \hat{Y} + X_{\mathcal{F}_c},$ where 
\[
\mathcal{F}_c =\mathcal{F}_{c,\varrho,\mathfrak{n}}:= \ell_c E_c + g_{c,\varrho, \mathfrak{n}}\Omega_{c,\varrho},
\]
  with $E_c$ and $\Omega_{c, \varrho}$  as in \eqref{eq:firstinterlE} and \eqref{eq:secondomnega}, respectively,  is smooth on $\mathcal{U}_{\mathcal{W}}.$ 
 Let $\mathcal{V}_{\mathcal{W}} \subset \mathcal{U}_{\mathcal{W}}$ be a smaller open neighborhood of $\mathcal{W}_{c,\rho,\varrho}$ such that $\overline{\mathcal{V}_{\mathcal{W}}} \subset \mathcal{U}_{\mathcal{W}}.$ We then choose  a compactly supported smooth function $\mathfrak{f} \colon \mathcal{U}_{\mathcal{W}} \to [0,1]$ satisfying $\mathfrak{f} \equiv 1 $ on $\mathcal{V}_{\mathcal{W}} $ and extend $\mathfrak{f}\mathcal{F}_c$ to a smooth function   on $T^*(\R \times S^1)$  by setting it equal to zero outside $\mathcal{U}_{\mathcal{W}}.$  We now define 
 \[
 Y:= \hat{Y}+X_{\mathcal{G}},\qquad 
\mathcal{G}:= \sum_{\lambda \in \Lambda} (\mathfrak{f}\mathcal{F}_c) \circ \lambda^{-1}.
\]
Note that $Y$ is a  $\Lambda$-invariant Liouville vector field defined on $T^*(\R \times S^1).$ By construction, $Y$ agrees with $\bar{Y}$ on a neighborhood of $\tilde{K}_{c}^{-1}(0).$

\begin{theorem}\label{thm:perturbation} 
Fix $0<\rho<\rho_1$ and $\mathfrak{n} \geq \mathfrak{n}_1,$ and let $\varrho_0=\varrho_0(\rho,\mathfrak{n})>0$ be as above. For every $0<\varrho<\varrho_0$ and every $c_2<c<c_2+\varrho^3,$   the perturbed regularized energy level $\tilde{\Sigma}_c = \tilde{K}_c^{-1}(0)$ is compact and diffeomorphic to $S^1 \times S^2.$ Moreover, there exists a $\Lambda$-invariant   Liouville vector field $Y$ defined on $T^*(\R \times S^1)$   which is positively  transverse to $\tilde{\Sigma}_c.$ 
\end{theorem} 
\begin{proof}
Let $z \in \tilde{\Sigma}_c.$ We first assume that 
\begin{equation}\label{eq:firstasumtppat}
z \in \mathcal{M}_{c,\rho,\varrho} \setminus \bigcup_{\lambda\in \Lambda} \lambda\left( \mathcal{B}_{r_*}(z_*(c)) \cap \{ \lvert \xi_2 \rvert < \frac{\pi}{4}\}\right).
\end{equation}
We claim that 
\[
\tilde{K}_c = K_c \qquad \text{on } \mathcal{M}_{c,\rho,\varrho}.
\]
For every $ z \in \tilde{\Sigma}_c \setminus   K_c^{-1}(0),$   by construction, there is $\lambda \in \Lambda$ such that 
\[
2\varrho < \xi_1=\xi_1(\lambda^{-1}z) < 4\varrho.
\]
Lemma \ref{lem:perturbedsall} together with \eqref{eq;crho-sd/2} implies that
\[
\lambda^{-1}z \in \tilde{\Sigma}_c \cap \mathcal{B}_{r_*/2}(z_*(c)) \cap \{ \lvert \xi_2 \rvert<\frac{\pi}{4}\} \cap \{ \varrho \leq \xi_1 < 4 \varrho\} \subset \mathcal{U}_{c,\varrho}.
\]
See also the paragraph immediately following Lemma \ref{lem:perturbedsall}.
Then by Lemma \ref{lem:newlame1}
\[
\alpha_c(\lambda^{-1}z) \geq \mathfrak{m}_0\varrho >0>-\frac{9}{4}\rho,
\]
and hence 
\[
\lambda^{-1}z \in \mathcal{N}_{c,\rho,\varrho}(z_*(c)).
\]
It follows that
\[
 {\tilde{\Sigma}_c} \setminus K_c^{-1}(0) \subset \bigcup_{\lambda \in \Lambda}\lambda( \mathcal{N}_{c,\rho,\varrho}(z_*(c))).
\]
The claim follows by  definition of $\mathcal{M}_{c,\rho,\varrho}.$ Thus,   $z \in K_c^{-1}(0).$ Moreover, 
since the capping functions are non-negative and vanish on $\mathcal{M}_{c,\rho,\varrho},$ the differentials of $\tilde{K}_c$ and $K_c$ also agree on $\mathcal{M}_{c,\rho,\varrho}.$

By the choice of $r_*,$ every $\zeta \in \mathcal{B}_{r_*}(z_*(c))$ satisfies $\lvert \xi_2 \rvert < \pi/4.$  Then assumption \eqref{eq:firstasumtppat} implies that $z$ lies outside the four coordinate $r_*$-balls. Since Lemma \ref{lem:perturbedsall} (and its symmetric versions) says that every point of $\tilde{\Sigma}_c$ with $\lvert x_1 \rvert \geq x_*(c)+\varrho$ lies in one of these four balls, it follows that $z$ satisfies $\lvert x_1 \rvert< x_*(c)+\varrho.$ Suppose that $\lvert x_1 \rvert > x_*(c)-\delta_0.$ Since $\varrho <\varrho_0< \delta_0/4$ (see \eqref{eq:varrho0intro}), we have
\[
x_*(c)- \delta_0 < \lvert x_1 \rvert < x_*(c)+\varrho < x_*(c) +\delta_0,
\]
implying that $\lvert \lvert x_1 \rvert -x_*(c) \rvert <\delta_0.$ Then by the choice of $\delta_0>0,$ the point $z$ lies in one of the four coordinate balls, see Lemma \ref{lem:locallemman} (and its symmetric versions), a contradiction. 
Therefore, $\lvert x_1 \rvert\leq x_*(c) - \delta_0,$ and hence the point $z$ is contained in the central part $\mathcal{K}_c,$ see \eqref{eq:centralpart}. Then, since $Y = \hat{Y}$ at $z,$ Lemma \ref{lem:additionalestimateforhaty} provides
\[
Y(\tilde{K}_c)(z)=\hat{Y}(K_c)(z) \geq \sigma_0>0.
\]

 We now assume that there is $\lambda \in \Lambda$ such that 
 \[
 z \in  \mathcal{M}_{c,\rho,\varrho} \cap    \lambda\left( \left[\mathcal{B}_{r_*}(z_*(c)) \cap \{ \lvert \xi_2 \rvert < \frac{\pi}{4}\}  \right]   \setminus  \mathcal{W}_{c,\rho,\varrho}\right).
\] 
By definition of $Y,$ we have $Y = \hat{Y}$ at $z.$ Moreover, since
\[
\alpha_c(\lambda^{-1}z) <-\frac{5}{2}\rho<-2\rho, \quad \lambda^{-1}z \in    \left[\mathcal{B}_{r_*}(z_*(c)) \cap \{ \lvert \xi_2 \rvert <  {\pi}/{4}\}  \right]   \setminus  \mathcal{W}_{c,\rho,\varrho} ,
\]
we have $\hat{Y} = \tilde{Y}_{\rm loc}$ at $\lambda^{-1}z.$ Since $\tilde{K}_c = K_c$ and $d\tilde{K}_c = dK_c,$ now the $\Lambda$-invariance of $\hat{Y}$ and $K_c$ and Lemma \ref{lem:firstinterpolatoinsuccess} give positive transversality.

Finally, suppose that $z \in \lambda(\mathcal{W}_{c,\rho,\varrho})$ for some $\lambda \in \Lambda,$ so that $Y = \bar{Y}_{\rm loc}^{\lambda}$ at $z.$  Now positive transversality follows from Lemma \ref{lem:latinterpolation} and the $\Lambda$-invariance of $\tilde{K}_c.$ 
This completes the proof.     
\end{proof}



\section{Birkhoff regularization and its relation to elliptic coordinates}
The purpose of this section is to introduce Birkhoff's approach to regularizing the PCR3BP and to show that it is closely related to the regularization via elliptic coordinates by constructing an explicit symplectomorphism between the  corresponding regularized Hamiltonian systems.

As with elliptic coordinates, Birkhoff's approach \cite{Birk15regularization} allows a simultaneous regularization around two primaries at once. To perform the regularization we define 
	\[
	\mathfrak{B}:=T^*\mathfrak{b} \colon T^*\left(\mathbb{C}\setminus \left\{-\tfrac{1}{2}, 0, \tfrac{1}{2}\right\}\right) \to T^*\mathbb{C}; \ \ (z,w)\mapsto \left(\frac{z^2+\frac{1}{4}}{2z}, \frac{2\bar{z}^2w}{\bar{z}^2-\frac{1}{4}}\right),
	\label{Bmap}
	\]
	which is the cotangent lift of the map
	\[
	\mathfrak{b}\colon  \mathbb{C}\setminus \{0\} \to \mathbb{C}; \ \ z \mapsto \frac{1}{4}\left(2z+\frac{1}{2z}\right).
	\label{bmap}
	\]
    Note that in these coordinates the position $z=0$ corresponds to a point $q=\infty$. Before pulling back the Hamiltonian system with this map, the coordinate of the PCR3BP has to be shifted such that both primaries are at the position $\left(\pm\frac{1}{2},0\right).$ For $\mu=\frac{1}{2}$ this is automatically fulfilled. Then the corresponding Birkhoff regularized Hamiltonian for the energy level $c$ is defined to be
    \begin{equation}
        \begin{aligned}
          &  H_B^c\colon   T^*\left(\mathbb{C}\setminus\{0\}\right) \to \mathbb{R}  ;\\
	  &  H_B^c(z,w)=   \mathfrak{B}^*\left(\Big|q-\frac{1} {2}\Big|\Big|q+\frac{1}{2}\Big|\big(H(q,p)-c\big)\right)(z,w).
    \end{aligned}
    \end{equation}
    Here, $H$ denotes the Hamiltonian of the PCR3BP.
	The resulting Hamiltonian is   not singular at $z=\pm\frac{1}{2}$ anymore.

    Suppose that $c_1<c<c_2.$  Denote by $\overline{\Sigma}_c^{\mathfrak B}$ the     bounded component of $({H^c_B})^{-1}(0),$ which is a closed manifold double covering the bounded component of $H^{-1}(c)$ with the two fiber circles 
	\[
    \Lambda_c^{E}:=\overline{\Sigma}_c^\mathfrak{B} \cap T^*_{\frac{1}{2}}\mathbb{C},\qquad	\Lambda_c^{M}:=\overline{\Sigma}_c^\mathfrak{B} \cap T^*_{-\frac{1}{2}}\mathbb{C}
	\]
	 added over the primaries   at $\left(\frac{1}{2},0\right)$ and    $\left(-\frac{1}{2},0\right),$ respectively. 
    
	As in the elliptic coordinates, the Birkhoff regularized Hamiltonian system gains an additional $\mathbb{Z}_2$-symmetry. More precisely, the $\Z_2$-action  
    \begin{align*}
        &\sigma_B\colon   \Z_2\times T^*(\C\setminus\{0\}) \to T^*(\C\setminus\{0\})\\
    &\sigma_B(1,(z,w)):=  (z,w)\ ,\quad \sigma_B( -1, (z,w)):= \left(\frac{1}{4z}, -4\bar{z}^2w\right)
    \end{align*}
  preserves the standard symplectic form coming from $T^*\mathbb{C}$, and the map $\mathfrak{B}$ is invariant under this action, i.e.~$\mathfrak{B}(z,w)=\mathfrak{B}\left(\sigma_B(-1,(z, w))\right)$. Therefore, the Birkhoff regularized Hamiltonian $H_B^c$ is invariant under this action.

\begin{figure}[h]
		\begin{center}
			\begin{tikzpicture}
				\draw (-6,0) --(3,0);
				\draw[->] (-1.5,-1.5) -- (-1.5,1.5);
				\draw[purple!80!black][ultra thick] (-6,0)--(-3,0);
				\draw[green!40!black][ultra thick] (0,0)--(3,0);
				\draw[gray!80!black] (-3, 0) node{$\bullet$};
				\draw[gray!80!black] (-3, 0) circle[radius=0.1];
				\fill[gray!80!black] (-3, 0) circle[radius=0.1];
				\draw[blue!80!black] (0, 0) node{$\bullet$};
				\draw[blue!80!black] (0, 0) circle[radius=0.1];
				\fill[blue!80!black] (0, 0) circle[radius=0.1];
				\draw[purple!80!black] (-4.5,-1.3) node{$\Gamma_c^{M}$};
				\draw[purple!80!black] (-6,-0.3) to[out=270, in=90, looseness=1] (-4.5,-1);
				\draw[purple!80!black] (-3,-0.3) to[out=270, in=90, looseness=1] (-4.5,-1);
				\draw[green!40!black] (3,-0.3) to[out=270, in=90, looseness=1] (1.5,-1);
				\draw[green!40!black] (0,-0.2) to[out=270, in=90, looseness=1] (1.5,-1);
				\draw[green!40!black] (1.5,-1.3) node{$\Gamma_c^{E}$};
				\draw[gray!50!black] (-2.7,0.5) node{$\Lambda_c^{M}$};
				\draw[blue!80!black] (0.4,0.5) node{$\Lambda_c^{E}$};
				\draw (-1.55,1.7) node{$q_2$};
				\draw  (3.25,0 ) node{$q_1$};

			\end{tikzpicture}
		\end{center}
		\caption{Interpretation of the Legendrians $\Gamma_c^E$, $\Lambda_c^{E}$, $\Gamma_c^M$ and $\Lambda_c^{M}$ in the position space  of the unregularized system}
        \label{fig:Legendrians}
	\end{figure}
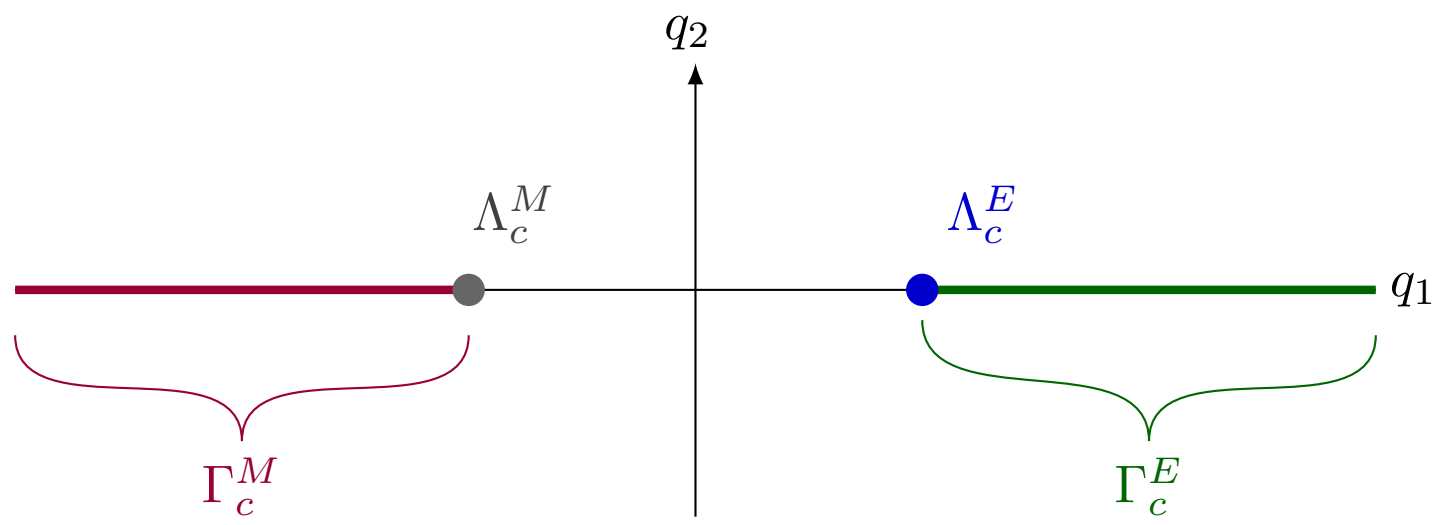

    The unregularized  system admits the   anti-symplectic involution 
    \begin{equation*}\label{eq:unregularirho}
        \varrho \colon  \R^4 \to \R^4,  \qquad (q_1,q_2,p_1,p_2) \mapsto (q_1,-q_2,-p_1,p_2)
    \end{equation*}
    which  leaves $H$ invariant. Since the map $\mathfrak{b}$  is a (branched) double covering, the involution $\varrho$ lifts to two anti-symplectic involutions 
    \begin{equation}\label{eq:twoBirkhoffrho}
        \varrho_1(z,w):= (\bar{z}, -\bar{w}),\qquad \varrho_2(z,w):= \left(\frac{1}{4\bar{z}}, 4z^2\bar{w}\right)
    \end{equation}
    in the Birkhoff regularized system. Note that $\varrho_2(z,w) = \sigma_B(-1, \varrho_1(z,w)). $
    For simplicity we consider only $\varrho_1.$ Its fixed point set ${\rm Fix}(\varrho_1)=(\R\setminus \{ 0 \}) \times i \R$ is a Lagrangian of $\C^2,$  and its intersection with  $\overline{\Sigma}_c^{\mathfrak B}$ decomposes as 
    \[
        \mathrm{Fix}(\varrho_1) \cap \overline{\Sigma}_c^{\mathfrak B} = \Gamma_c^E \sqcup \Gamma_c^M,
    \]
    where $\Gamma_c^E$ and $\Gamma_c^M$ are two connected Legendrian submanifolds. The submanifold  $\Gamma_c^E$ intersects $\Lambda_c^{E},$ but not $\Lambda_c^M,$ while  $\Gamma_c^M$ intersects $\Lambda_c^{M}$, but not $\Lambda_c^E.$ In the position space  of the unregularized   system, $\Gamma_c^E$ corresponds to the ray $(\frac{1}{2},\infty)$ on the $q_1$-axis and $\Gamma_c^M$ to $(-\infty, -\frac{1}{2} )$, see Figure~\ref{fig:Legendrians}.

Now we construct a symplectomorphism between the Hamiltonian systems in elliptic   and Birkhoff coordinates. Since the latter is formulated using complex numbers, we shall also rewrite the elliptic coordinates transformation (see Section~\ref{sec:ellipticcoord}) in the same way. Under   the usual identification $ \R^2 \to \C  ,\; (v_1,v_2) \mapsto v_1+iv_2,$  
the elliptic coordinates become
\[
    q=\frac{1}{2}\cosh (x)  , \qquad p=\frac{2y}{\overline{\sinh(x)}} = \frac{2y \sinh(x)}{\lvert \sinh(x)\rvert^2}
\]

\begin{lemma}    \label{lem:Phi}

    The map 
    \[
        \Phi \colon  (\mathbb{R}+i(\mathbb{R}/2\pi\mathbb{Z}))\times\mathbb{C} \to (\mathbb{C}\setminus \{0\})\times \mathbb{C}, \quad (x,y)\mapsto (z,w)=\left(\frac{1}{2}e^x, 2e^{-\bar{x}}y\right)
    \]
    is a  symplectomorphism with respect to the standard symplectic forms.  Furthermore, it satisfies the equation
    \[
        \Phi^*H^c_B = K_c.
    \]
\end{lemma}
\begin{proof}
    It is easy to see that the complex exponential function maps $\mathbb{R}+i[0,2\pi)$ bijectively to $\mathbb{C}\setminus \{0\}$. Since this in particular implies that   $2e^{-\bar{x}}$ is never zero, we  conclude that the whole map $\Phi$ is a bijection.  
    It has a smooth inverse, 
    and hence is a diffeomorphism. 

    In the complex coordinates $(z,w)$  one can write the standard symplectic form as $\omega = \sum_{i=1}^2 dw_i\wedge dz_i =\textrm{Re}(d\bar{w}\wedge dz).$
    With this we calculate
    \[
        \Phi^*\omega = \textrm{Re}\left(d\left(2e^{-x}\bar{y}\right)\wedge d\left(\frac{1}{2}e^x\right)\right) = \textrm{Re}(d\bar{y}\wedge dx) = \omega,
    \]
    where $x=x_1+ix_2$ and $y=y_1+iy_2$. This proves the first part.

    For the second claim we first compute the composition $ \mathfrak{B} \circ \Phi.$ If $\sinh(x) \neq 0,$  equivalently $ z \neq \pm \frac{1}{2},$ then we have
    \begin{align*}
        \frac{z^2 +\frac{1}{4}}{2z }&= \frac{\frac{1}{4}e^{2x}+\frac{1}{4}}{e^x}=\frac{1}{4}\left(e^x+e^{-x}\right)=\frac{1}{2}\cosh{(x)}\\
        \frac{2\bar{z}^2 w }{\bar{z}^2 -\frac{1}{4}}&= \frac{2(\frac{1}{4}e^{2\bar{x}})(2e^{-\bar{x}}y)}{\frac{1}{4}e^{2\bar{x}}-\frac{1}{4}} =\frac{4y}{  e^{\bar{x}}-e^{-\bar{x}} } =\frac{2y}{\overline{\sinh(x)}}
    \end{align*}
Therefore, $\mathfrak{B} \circ \Phi$ is precisely the transformation into elliptic coordinates. 
  Combining this with the definition of $H^c_B,$ a direct computation shows that, away from the collisions, $\Phi^* H^c_B = K_c.$ 
Since both sides extend smoothly over the collisions, this finishes the proof.
\end{proof}

We now show that   $\Phi$ is compatible with the $\Z_2$-symmetries.

\begin{lemma}\label{lem:sigmaB}
    The symplectomorphism  $\Phi$  is $\Z_2$-equivariant, i.e.
    \[
        \sigma_B(\pm1,\ \cdot\ ) \circ \Phi = \Phi\circ\sigma_{\rm el}(\pm1,\ \cdot\ ).
    \]
    
\end{lemma}
\begin{proof}
    The case of $+1$ is trivial. For $-1,$ 
we find that
\begin{align*}
        \sigma_B(-1,\Phi(x,y))&=\left(\frac{1}{4}\frac{1}{\frac{1}{2}e^x}, -4\frac{1}{4}e^{2\bar{x}}2 e^{-\bar{x}}y\right) \\ &= \left(\frac{1}{2}e^{-x}, -2e^{\bar{x}}y\right) \\ & = \Phi\left(\sigma_{\rm el}(-1,(x,y))\right)
    \end{align*}
    which completes the proof.
\end{proof}

In \cite{kang2025consecutive} it was shown that there exists a positively  transverse Liouville vector field in a neighborhood of $\bar{\Sigma}^\mathfrak{B}_c$ for energies $c$ slightly above $c_1$. The $\Z_2$-equivariant symplectomorphism constructed above allows us to pull  this Liouville vector field back to elliptic coordinates. We shall denote this Liouville vector field by $Y_B.$ Then the connected component $\Sigma_c \subset K_c^{-1}(0)$ carries two contact forms 
\begin{equation}\label{eq:twpoconact}
 \lambda_B:= \iota_{Y_B} \omega |_{\Sigma_c} \quad \text{ and } \quad   \lambda_\mathrm{el}:=\iota_{\hat{Y}}\omega |_{\Sigma_c},
 \end{equation}
that are related in the following way.

\begin{proposition} \label{prop:elToBirk}
    For $c=c_1+\varepsilon$ with $\varepsilon>0$ small enough, there exists a $\Z_2$-equivariant contactomorphism
    \[
    \Psi \colon  \left(\Sigma_c,\xi_\mathrm{el}=\mathrm{ker}(\lambda_\mathrm{el})\right) \to \left(\Sigma_c,\xi_B=\mathrm{ker}(\lambda_B)\right),
    \]
    such that for any closed submanifold $\Lambda$ of $\Sigma_c$ that is Legendrian with respect to  both $\lambda_\mathrm{el}$ and $\lambda_B,$ we have  $\Psi(\Lambda)=\Lambda$.
\end{proposition}
\begin{proof}
    The proof is mainly based on the Moser trick: Since both Liouville vector fields $\hat{Y}$ and $Y_B$ are positively transverse to  $\Sigma_c$, we can define a smooth family of positively transverse Liouville vector fields
    \[
        Y_t:=(1-t)\hat{Y} + t Y_B, \qquad t\in [0,1]
    \]
    in a  neighborhood of $\Sigma_c$. This induces a smooth family of contact forms
    \[
        \lambda_t:=\iota_{Y_t}\omega |_{\Sigma_c}=(1-t) \lambda_\mathrm{el}+t\lambda_B.
    \]
Let $ \xi_t = \ker \lambda_t.$ Since $\mathrm{d}\lambda_t$ is non-degenerate when restricted to $\xi_t$, there is a unique vector field $V_t\in \xi_t$ such that 
\[
d\lambda_t( V_t, \cdot)\big\vert_{\xi_t} = -\partial_t \lambda_t\big\vert_{\xi_t} + \mu_t\lambda_t\big\vert_{\xi_t},
\]
where $\mu_t:=( 
\lambda_B-\lambda_\mathrm{el})(X_{K_c})/\lambda_t(X_{K_c})$. Note that $X_{K_c}$ is a positive reparameterization of any Reeb vector field induced by  $\lambda_{t}$. Furthermore, the tangent bundle of $\Sigma_c$ splits as $T\Sigma_c= \xi_t \oplus \R X_{K_c}$ 
for all $t\in[0,1]$. Since  $d\lambda_t( V_t, X_{K_c}) = -\partial_t \lambda_t(X_{K_c}) + \mu_t\lambda_t(X_{K_c})=0,$
$V_t$ now interpreted as a vector field on $\Sigma_c$ still satisfies its defining equation. Denote by $\{ \Psi_t \}_{t \in [0,1]}$ the flow of $V_t$ such that $\Psi_0={\rm Id}.$ We then find
\[
\frac{d}{dt}\Psi_t^*\lambda_t = \Psi_t^* ( \partial_t \lambda_t + \mathcal{L}_{V_t}\lambda_t) = \Psi_t^* ( \partial_t \lambda_t +  \iota_{V_t}d\lambda_t + d(\lambda_t(V_t))) =\Psi_t^*\mu_t \cdot \Psi_t^*\lambda_t. 
\]
At any point $v_x\in T\Sigma_c$ this defines an ODE in $t$ with the unique solution
\[
    (\Psi_t^*\lambda_t)(v_x)= \exp\left(\int_0^t\Psi_{\tau}^*\mu_\tau(x)d\tau\right)\cdot\lambda_\mathrm{el}(v_x).
\]
It follows that $T\Psi_t(\xi_\mathrm{el})=\xi_t.$ 
Note that since both contact forms $\lambda_{\rm el}$ and $\lambda_B$ are invariant under $\sigma_{\rm el},$ the family $\lambda_t$ is also $\sigma_{\rm el}$-invariant. Then by uniqueness, the vector field $V_t$ is $\sigma_{\rm el}$-invariant, and hence its flow $\Psi_t$ is $\Z_2$-equivariant. Then $\Psi:=\Psi_1$ is the desired contactomorphism.

Pick a closed submanifold $\Lambda \subset \Sigma_c$ that is Legendrian with respect to both $\lambda_B$ and $\lambda_{\rm el}.$ Then it is a Legendrian submanifold with respect to $\lambda_t$ for all $t \in [0,1].$  For $v \in T_x\Lambda,$ 
\[
d\lambda_t( V_t, v) = -(\partial_t \lambda_t)(v) + \mu_t(x)\lambda_t(v)=0,
\]
Since   $T_x \Lambda$ is a Lagrangian subspace of $(\xi_t(x), d\lambda_t|_{\xi_t(x)})$ for every $x \in \Lambda$ and every $t \in[0,1],$ this implies that $V_t(x) \in T_x \Lambda$ for every $x \in \Lambda.$ In particular, $\Lambda$ is preserved along the flow of $V_t,$ implying that $\Psi(\Lambda)=\Lambda.$ This completes the proof. 
\end{proof}

This relation between the Birkhoff regularization and the elliptic coordinates will be one of the key ingredients in the proof of Theorem~\ref{thm:ccorbits} in Section \ref{sec:proofsection}.


\section{Consecutive collision orbits and Rabinowitz Floer homology}

\subsection{The basics of Lagrangian Rabinowitz Floer homology} \label{sec:basiscsRFH}

 Rabinowitz Floer homology (RFH) is an invariant associated to a contact-type hypersurface in a   symplectic manifold.  It captures information about the existence of periodic orbits \cite{cieliebak2009a} or Reeb chords with Lagrangian boundary conditions   \cite{merry2014a}. The latter  is usually referred to as Lagrangian RFH. To avoid unnecessary technical discussion we restrict  ourselves   to the specific setting used in the proof of Theorem~\ref{thm:ccorbits}.


 Let $M=S^1\times\R^3,$ with circle coordinate $x_2,$  be endowed with its standard symplectic structure $\omega = dy_1 \wedge dx_1 + dy_2 \wedge dx_2,$  let $L_1, L_2 \subset M $ be  Lagrangian submanifolds of the form $\{t_j\} \times V_j,$ where $V_j \subset \R^3$ is a two-dimensional linear subspace, and let  $\Sigma\subset M$ be a closed hypersurface with a     Liouville vector field $X_L$ that is positively transverse to $\Sigma$. This, in particular, implies that $\lambda:= \iota_{X_L}\omega \in \Omega^1(M)$ is a primitive of $\omega,$ and  the restriction $\alpha:=\lambda|_{\Sigma}$ is a contact form on $\Sigma.$  

Choose  a Hamiltonian  $H \colon M \to \R,$ constant outside a compact set, such that 
\begin{equation}\label{eq:RFHHamade}
    H^{-1}(0)=\Sigma \qquad \text{and} \qquad  \alpha(X_H)= dH(X_L)|_{\Sigma}>0, 
\end{equation}
so that   $X_H|_\Sigma=\alpha(X_H) R  $ is a positive multiple of the Reeb vector field  $R=R_\alpha$  of $\alpha.$  Such a Hamiltonian is called \textit{adapted to} $(\Sigma,\alpha).$ 

\begin{remark}\label{rmk:XHnotR}
    As in \cite[Section 2]{merry2014a}, we do not impose the normalization $X_H|_\Sigma=R.$ Keeping the original Hamiltonian near $\Sigma$ will be useful  for the index estimates in the next section. 
\end{remark}

We first consider  $L= L_1= L_2.$ Assume that $L$ intersects $\Sigma$ transversely.
We set
\[
\mathscr{P}(M,L ):= \{  x\in C^\infty([0,1],M)  \mid  x(0) ,x(1)\in L \} 
\]
and  define the \textit{Rabinowitz action functional} $\mathcal{A}_H\colon  \mathscr{P}(M,L ) \times\R \to \R$ as
\[
\mathcal{A}_H(x,\tau)=   \int _0^1 x^*\lambda -\tau \int _0^1 H(x(t))dt  - l (x(1))+l (x(0)),
\]
where $l $ is a smooth function on $L $ such that $\lambda\big\vert_{L }=\mathrm{d}l .$ 
 The critical points of this functional are precisely the pairs $(x,\tau)$ satisfying
\begin{align*}
    \begin{cases}
        \partial_tx(t)=\tau X_H(x(t)),\quad t\in[0,1]\\
    H(x(t))=0.
    \end{cases}
\end{align*}
Since $H$ is adapted to $(\Sigma,\alpha)$ and time-independent, the second  equation implies that    $x(t) \in \Sigma$ for every $t \in [0,1].$ We thus have three types of critical points:  
If $\tau>0,$ the reparametrized path $\gamma_x \colon [0,\tau] \to \Sigma,$ $\gamma_x(t)=x(t/\tau),$ is a Hamiltonian chord with endpoints in $L,$  and hence it becomes a Reeb chord after a suitable positive time-reparametrization.   If $\tau<0,$ then the path $x$ follows the negative direction. Finally, critical points with  $\tau=0 $ correspond precisely to points of   $\Sigma \cap L ,$ via constant paths.  Thus, they   form a one-dimensional Morse-Bott critical manifold, identified with $\Sigma \cap L  .$
See \cite[Lemma 2.13]{merry2014a}.

The construction of RFH is analogous to Morse homology, where  the action functional  $\mathcal{A}_H$ can heuristically be thought of as a Morse function on an infinite-dimensional space. The generators are its critical points or, in the Morse-Bott case, the critical points of auxiliary Morse functions on its critical manifolds. 

We recall the  {Robbin--Salamon index}. The coordinate vector fields provide a global symplectic trivialization of $TM ,$ denoted by $\mathfrak{T}.$ Since the Lagrangian $L$ is of the form $\{t_0\} \times V,$ where $V \subset \R^3$ is a two-dimensional linear subspace,  under this trivialization, its tangent spaces are identified  with the fixed Lagrangian subspace $\{ 0 \} \times V \subset \R \times \R^3.$ By abuse of notation, we identify $V$ with $\{0\} \times V.$ 
Let $z=(x,\tau)$ be a critical point with $\tau \neq 0.$ Consider  
\[
\psi_z(t):= \mathfrak{T}_{x(t)} \circ D\Phi_H^{\tau t}(x(0)) \circ \mathfrak{T}_{x(0)}^{-1}  , \qquad t \in [0,1].
\]
This also includes chords with $\tau<0.$ 
We define the {\it full Robbin--Salamon index} by
\[
\mu_{\rm RS}(z):=\mu_{\rm RS}(\psi_z V, V),
\]
using the convention of Robbin and Salamon \cite{RS93}, in which each regular endpoint crossing contributes half the signature of its crossing form. The Hamiltonian chord is nondegenerate if $\psi_z(1)V \cap V = \{0\}.$ In this case, $\mu_{\rm RS}(z) \in \Z$ since $\psi_z(0)V=V$ and $\dim V = 2.$

 For the grading, we use the convention due to Merry \cite[Definition 2.22]{merry2014a}, with a shift by $+\frac{1}{2}.$ For a nondegenerate critical point $z =(x,\tau)$ with $\tau\neq0,$ we set
 \[
 \lvert z \rvert = \mu_{\rm RS}(z) - \frac{1}{2}\chi(z)+\frac{1}{2} \in \Z,
 \]
where $\mu_{\rm RS}$ denotes the Robbin--Salamon index,  and   $\chi(z) \in \{\pm 1\}$ is the correction term defined in \cite[equation (2.11)]{merry2014a}. 
We recall this correction term in  Section \ref{sec:compat}. Under the hypotheses and with the compatible trivializations of Proposition  \ref{lem:compa},  
\[
\lvert z\rvert =\mu_{\mathrm{RS}} (z)-\frac{1}{2}\chi(z)+\frac{1}{2}
=\mu_{\mathrm{RS}}^{\rm tr} (z)+\frac{1}{2},
\]
where $\mu_{\rm RS}^{\rm tr}(z)$ denotes the {\it transverse Robbin--Salamon index}. In   our application, we need not   compute $\chi(z)$ explicitly. 
The grading of a constant generator $(p,0) $   is defined using an auxiliary Morse function on $\Sigma \cap L. $ Since we will not need its explicit description,   we refer to \cite[Definition 2.22]{merry2014a} for details, keeping the shift by $+\frac{1}{2}.$ Thus, the chain group ${\rm CRF}_k$ is generated over $\Z_2$ by the degree-$k$ generators, with the standard completion, see \cite[Section 2.3]{merry2014a}.

 To define the boundary operator, we   count  finite-energy solutions $ (u,\eta) \colon  \R \to \mathscr{P}(M,L ) \times\R$  to the Floer equation
\begin{equation}
    \begin{cases}
        \partial_su- J(u)\left(\partial_tu-\eta X_H(u)\right)=0,\\
        \partial_s\eta- \int_0^1 H(u(s,t))dt=0,
    \end{cases}
    \label{eq:floer}
\end{equation}
where $J$ is an almost complex structure compatible with the symplectic form $d\lambda$. From the perspective of Morse homology one can interpret these solutions $(u,\eta)$ as  negative gradient flow lines.  The relevant metric on $\mathscr{P}(M,L) \times\R$ is 
\[
    g_J|_{ (x,\tau)}((\hat{x}_1,\hat{\tau}_1),(\hat{x}_2,\hat{\tau}_2))= \int _0^1 d\lambda(\hat{x}_1(t),J(x(t))\hat{x}_2(t)) dt + \hat{\tau}_1 \cdot\hat{\tau}_2
\]
As in  Morse homology, we identify solutions that differ by a translation in the $s$-component. 

    
    For $z_\pm \in  \mathrm{crit}(\mathcal{A}_H),$  let $\mathcal{M}(z_-, z_+)$ denote the moduli space of finite-energy solutions of \eqref{eq:floer} from $z_-$ to $z_+,$  modulo translation in $s.$  For generic data, its dimension is $\lvert z_-\rvert - \lvert z_+\rvert -1.$    The boundary   operator $\partial  \colon {\rm CRF}_k(M,\Sigma,L) \to {\rm CRF}_{k-1}(M,\Sigma,L)$ is given by  
    \[ 
                 \partial  z_- :=   \sum\limits_{\substack{ z_+ \in {\rm crit}(\mathcal{A}_H)  \\ \lvert z_+\rvert = \lvert z_- \rvert -1}} \#_2 \mathcal{M}(z_-, z_+)  \cdot  z_+,
    \]  
where the symbol $\#_2$ indicates the count modulo two. For Morse--Bott critical manifolds, the generators are the auxiliary Morse critical points and the same formula counts rigid cascades of Floer and Morse trajectories. 
The standard compactness and gluing arguments give $\partial^2=0$ for generic Floer data, showing that $({\rm CRF}_*(M, \Sigma, L), \partial)$ is a chain complex. See \cite[Section 2.2--2.3]{merry2014a} for more details. 


\begin{definition}
    The Lagrangian Rabinowitz Floer homology is the homology associated to the Rabinowitz Floer complex, i.e.
    \[
        {\rm RFH}_k(M,\Sigma,L ):= \frac{\ker  (\partial |_{{\rm CRF}_k(M,\Sigma,L )} )}{\mathrm{im} (\partial|_{{\rm CRF}_{k+1}(M,\Sigma,L )}) }
    \]
\end{definition}


Each Floer trajectory connecting two critical points gives a homotopy between their underlying paths, with both endpoints free to move in $L.$ The boundary operator therefore   preserves the    path components of $\mathscr{P}(M,L ).$ 
We denote by $[\mathrm{pt.}]$ the component containing the constant paths.   The homology of the corresponding subcomplex will be denoted by   ${\rm RFH}_*(M,\Sigma,L  ,[\mathrm{pt.}]).$ 

\medskip

 Let a   finite group $G$ act  on $(M,\lambda)$  such that it leaves $\lambda$, $L  $ and $H$ invariant. We work in a $G$-invariant path component and choose all the Floer data equivariantly.  If the induced action on $\mathscr{P}(M,L )$ acts freely on ${\rm crit}(\mathcal{A}_H)$, one can  define an equivariant version of RFH, denoted by ${\rm RFH}_*^{G}(M,\Sigma,L),$ by taking generators modulo the induced $G$-action. Each generator is a  $G$-orbit of generators of the ordinary complex.  Since the induced symmetry preserves the RFH grading, $\lvert [z]\rvert :=\lvert z \rvert $ is well-defined.  To get a well-defined boundary operator, we choose a $G$-invariant  almost complex structure $J,$ satisfying 
\[
    \mathrm{d}\phi_g  \circ J(x) = J(\phi_g(x))  \circ \mathrm{d}\phi_g \qquad \forall g\in G,
\]
where $\phi_g$ denotes the action of $G$ on $M$. According to \cite[Section 4.5]{ruck24b} such a $J$ can be chosen to satisfy the required transversality conditions (note that a finite group does not need to act freely on all of $M$ for this proof to hold).   
Then  equivariance of the Floer data induces a free $G$-action on the sets of Floer trajectories asymptotic to critical points. This  allows us to define the boundary operator by
\[ 
                 \partial [z_-] =    \sum\limits_{\substack{[z_+] \in {\rm crit}(\mathcal{A}_H)/G  \\  \lvert z_+ \rvert = \lvert z_-\rvert -1 }} \left(  \sum\limits_{ \hat{z} \in [z_+]  }    \#_2  \mathcal{M}(z_-, \hat{z})   \right)\cdot [z_+]
    \]
The corresponding homology ${\rm RFH}_*^{G}(M,\Sigma,L)$ is   defined analogously to   ordinary RFH. For more details, see \cite[Section 5]{ruck24b}.

The key feature of RFH is that it is invariant under deformations of the contact  hypersurface $\Sigma $ induced by good homotopies in the sense of \cite[Definition 2.6]{merry2014a}, with $(M, \lambda)$ and $L $ fixed. 
In the equivariant setting, the homotopies and auxiliary choices are taken to be $G$-equivariant.

\medskip

We now turn to the case $L_1 \neq L_2.$ To describe chords between $L_1$ and $L_2,$ we use the perturbed functional, first introduced in \cite{albers2010a}, in its Lagrangian version developed in \cite[Section 2.4]{merry2014a}. See also \cite[Section 4.2]{kang2025consecutive}.

\begin{figure}[ht!]
    \centering
        \includegraphics[width=0.8\textwidth]{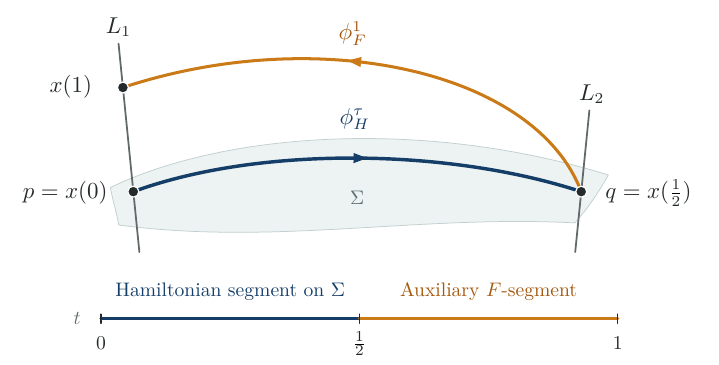}
        \caption{An illustration for a critical point $(x, \tau)$ of $\mathcal{A}_H^F$ with $\tau>0.$ The first half follows $H$ on $\Sigma$ from $L_1 $ to $L_2.$ The second half follows $F$ back to $L_1$ and need not remain on $\Sigma.$}
        \label{fig:traj}
\end{figure}

Let $F\colon [0,1] \times M \to \R$ be a compactly supported Hamiltonian  such that 
\[
(\phi_F^1)^{-1}(L_1)\cap \Sigma = L_2 \cap \Sigma.
\]
Thus,  for every $p \in \Sigma,$ we have $\phi_F^1(p) \in L_1$ if and only if $ p \in L_2.$ Choose a smooth   function  $\beta \colon [0,1] \to [0, \infty)$ such that  
\[
{\rm supp}(\beta) \subset (0, \tfrac{1}{2}), \qquad \int_0^1 \beta(t)dt=1 
\]
 and a smooth  nondecreasing function $\kappa \colon [0,1]\to[0,1]$ such that $\kappa(\frac{1}{2})=0$ and $\kappa =1 $ near $1.$ The Hamiltonian $\dot{\kappa}(t)F_{\kappa(t)}$ generates the flow $t \mapsto \phi_F^{\kappa(t)}.$
 We then define the perturbed action functional $\mathcal{A}_{H}^{F} \colon \mathscr{P}(M, L_1) \times \R \to \R$ by
 \begin{align*}
 \mathcal{A}_H^F(x, \tau):=   &\;   \int _0^1 x^*\lambda -\tau \int _0^1 \beta(t)H(x(t))dt\\
 & \;- \int_0^1 \dot{\kappa}(t) F_{\kappa(t)} (x(t))dt  - l_1(x(1))+l_1(x(0)),
 \end{align*}
where $l_1 \colon L_1 \to \R$ satisfies $\lambda|_{L_1}=dl_1.$  A critical point $(x,\tau)$ satisfies 
\begin{align*}
    \begin{cases}
        \dot{x}=\tau \beta  X_H(x )+ \dot{\kappa} X_{F_{\kappa }}(x  ), \\
    \int_0^1\beta(t)H(x(t))dt=0.
    \end{cases}
\end{align*}
Together with the choices of the functions $\beta$ and $\kappa$ this implies that 
\[
x(0) \in \Sigma \cap L_1, \quad x(\tfrac{1}{2}) \in \Sigma, \quad x(1) =\phi_F^1 ( x(\tfrac{1}{2})) \in L_1.
\]
Note that the last inclusion is equivalent to $x(\tfrac{1}{2}) \in L_2.$ See Figure \ref{fig:traj}. Thus, if $\tau>0,$ then the restriction of $x$ to $[0,\tfrac{1}{2}]$ gives a Hamiltonian chord on $\Sigma$ from $L_1$ to $L_2.$ So, it becomes a Reeb chord from $L_1$ to $L_2$ after a suitable positive time-reparametrization. If $\tau<0,$ then this restriction follows the Hamiltonian flow in the negative direction. Finally, critical points with $\tau=0$ are pairs $(x_p,0),$ where $p \in \Sigma \cap L_1\cap L_2$ and $x_p(t)=\phi_F^{\kappa(t)}(p), t \in [0,1].$ Note that $x_p(t)=p$ for $t \in [0,\tfrac{1}{2}].$

As in the single Lagrangian case, we work in the path component $[{\rm pt.}]$ of $\mathscr{P}(M, L_1) $ containing the constant paths.  The perturbed Rabinowitz Floer complex and its differential are defined as in \cite[Section 2.4]{merry2014a},   with the same  grading shift by $\tfrac{1}{2}$ as above. In particular, for every nondegenerate critical point $z =(x,\tau),$ its grading  $\lvert z \rvert$ satisfies  
\begin{equation}\label{eq:fadslfsih}
\lvert \lvert z \rvert - \mu_{\rm RS}(z) \rvert \leq 1,
\end{equation}
where $\mu_{\rm RS}(z)$ denotes the full Robbin--Salamon index of the   perturbed Hamiltonian chord  $x$, computed by linearizing it, with $\tau$ fixed.

We denote the corresponding homology by $   \mathrm{RFH}_*(M,\Sigma,L_1,L_2,[{\rm pt.}];F).$  By continuation invariance, the  homology is independent of the choice of $F,$  up to a degree-preserving continuation isomorphism. So, we suppress $F$ from the notation and write $   \mathrm{RFH}_*(M,\Sigma,L_1,L_2,[{\rm pt.}] ).$ In particular,  continuation to the unperturbed functional gives an isomorphism
\[
\mathrm{RFH}_k(M,\Sigma,L_1,L_2,[{\rm pt.}] )  \cong \mathrm{RFH}_k(M,\Sigma,L_1,[{\rm pt.}]).
\]
In the $G$-equivariant setting, we take all the data described above equivariantly. The same construction then provides an isomorphism
\[
\mathrm{RFH}_k^{G}(M,\Sigma,L_1,L_2,[{\rm pt.}] )  \cong \mathrm{RFH}_k^{G}(M,\Sigma,L_1,[{\rm pt.}]).
\]
 See \cite[Section 4.2]{kang2025consecutive} and \cite[Section 5]{Ruck24a}.

\subsection{Comparison of the indices}\label{sec:compat}

For the index estimate in the next section, it suffices to bound the difference between the full and transverse Robbin--Salamon indices of a Hamiltonian chord. We prove the exact comparison formula in  Proposition \ref{lem:compa} since it determines the grading shift relating our RFH grading to the Reeb chord grading used by Kang and Ruck in their model computation \cite{kang2025consecutive}. The proposition gives the chord analogue of the computation for periodic orbits established in    \cite[Lemma 2.9 and Corollary 2.10]{MP11}. We include a self-contained proof, as we have not found an explicit reference for the chord case.

  Let   $(M, \omega),$ $\Sigma,$ $X_L,$ $L,$ $\lambda ,$   $\alpha$ and $H$ be as given  at the beginning of Section \ref{sec:basiscsRFH}. Assume that $\lambda|_L=0$ near $\Sigma  .$  Then    $X_L$ is tangent to $L $ near $\Sigma$ and transverse to $\Sigma,$ so the intersection $\mathcal{L} = \Sigma \cap L$ is a Legendrian submanifold of $(\Sigma, \alpha)$ and  
  \begin{equation}\label{eq:splittinglan}
  TL|_{\mathcal{L}} = T\mathcal{L} \oplus\R X_L.
  \end{equation}

  Let $z = (x,\tau)$ with $\tau\neq0$ be a critical point of the unperturbed action functional $\mathcal{A}_H.$  
Let $\mathfrak{T} \colon TM \to M \times \R^4$ denote the  global symplectic trivialization given by the coordinate vector fields as above, under which the tangent spaces of $L$ are identified with  a fixed Lagrangian subspace $V \subset \R^4.$  Define 
\[
\psi_x(t):=\mathfrak{T}_{x(t)}\circ D\phi_H^{\tau t}(x(0)) \circ \mathfrak{T}_{x(0)}^{-1}, \quad   0 \leq t \leq 1.
\]
It is non-degenerate if $\psi_x(1)V \cap V = \{ 0 \}.$ We then define  $\mu_{\rm RS}(x)=\mu_{\rm RS}(\psi_x V, V).$

We now recall Merry's correction term  (see \cite[equation (2.11)]{merry2014a}).   Let $(x_s, \tau_s ), s \in (-\varepsilon, \varepsilon),$ be a smooth family of  Hamiltonian chords on the energy levels $H^{-1}(e_s),$ with endpoints in $L,$ where $\tau_s \neq 0$ and $(x_0, \tau_0, e_0)=(x, \tau,0).$ Set $T_s:=\lvert \tau_s \rvert>0.$ For such a family with $e_0' \neq 0$ and $T_0'\neq0,$ the correction term  is defined by
\begin{equation}\label{eq:crorr}
\chi(z):= {\rm sign}(\tau) \cdot {\rm sign}\left( - \frac{e_0'}{T_0'}\right)= {\rm sign}\left( - \frac{e_0'}{\tau_0'}\right)
\end{equation}

To define the transverse Robbin--Salamon index consistently with $\mathfrak{T},$  we first choose a symplectic trivialization $\tilde{\Phi} \colon x^*(\ker \alpha) \to [0,1] \times \R^2$ satisfying
\begin{equation}\label{eq:tildephiinitial}
\tilde{\Phi}_{x(t)} ( T_{x(t)}\mathcal{L}) = \R, \qquad t=0,1,
\end{equation}
where $\R:=\R \times \{ 0 \} \subset \R^2.$   We then
extend it by the frame $\{ X_L, R\}.$ Explicitly, we define the symplectic trivialization $\hat\Phi$ of $x^*TM $   by  
\begin{equation}\label{eq:sdilthat}
\hat\Phi_{x(t)} \colon T_{x(t)}M \to \R^2 \oplus \R^2, \quad \hat\Phi_{x(t)}( v + aX_L +bR )=( \tilde{\Phi}_{x(t)}(v), a, b),
\end{equation}
where $v \in \ker \alpha_{x(t)}$ and the target is equipped with the symplectic form $dv_1 \wedge dv_2 + da \wedge db.$
By \eqref{eq:splittinglan}, it sends both endpoint tangent planes of $L$   to $\R \oplus \R.$

Choose an ordered basis $\{e_1, e_2\}$ of $V$ and complete it to a symplectic basis $\{e_1, g_1, e_2, g_2\}$ of $(\R^4, \omega)$ with $\omega(e_i, g_j)=\delta_{ij}.$ Define the  linear map 
\begin{equation}\label{eqLAlinser}
A \colon \R^2 \oplus \R^2 \to \R^4, \quad (v_1, v_2, a, b) \mapsto v_1 e_1 + v_2 g_1 + a e_2 + b g_2.
\end{equation}
 Then   $A^*\omega = dv_1 \wedge dv_2 + da \wedge db$ and $A(\R \oplus \R)=V.$
We now define the symplectic trivialization $\Phi \colon x^*TM \to [0,1]\times \R^4$ by $\Phi_{x(t)}:=A \circ \hat{\Phi}_{x(t)}$ for every $t \in [0,1].$ Explicitly, 
\[
\Phi_{x(t)}( v + a X_L  + bR ) = A ( \tilde{\Phi}_{x(t)}(v), a,b).
\]
Note that $\Phi_{x(t)}(T_{x(t)}L)=V$ at $t=0,1.$  By a standard argument, we may choose $\tilde{\Phi}$ so that its extension $\Phi$ is homotopic to the coordinate trivialization $\mathfrak{T}$ through symplectic trivializations that send $T_{x(t)}L$ to $V$ for $t=0,1$ throughout the homotopy. See, for instance, the proof of \cite[Theorem 4.5]{KKK}.

Define
\[
\psi_x^{\rm tr}(t):= \tilde{\Phi}_{x(t)}  \circ \pi \circ D\phi_H^{\tau t}(x(0))|_{\ker \alpha_{x(0)}} \circ \tilde{\Phi}_{x(0)}^{-1},
 \]
 where $\pi \colon T\Sigma \to \ker \alpha$ denotes the projection along the Reeb direction, 
and
\[
\mu_{\rm RS}^{\rm tr}(x):=\mu_{\rm RS}(\psi_x^{\rm tr} \R, \R).
\]
The corresponding Reeb chord $x$ is non-degenerate if $\psi_x^{\rm tr}(1)\R \cap \R = \{ 0 \}.$

\begin{remark}
    As noted in Remark \ref{rmk:XHnotR}, $X_H|_{\Sigma}=fR,$ where the positive function $f = \alpha(X_H)$ may vary on $\Sigma.$  
    Therefore, in the above definition of $\psi_x^{\rm tr},$ we need the projection $\pi.$ Indeed, even for $v \in \ker \alpha_{x(0)},$  the vector $D\phi_H^{\tau t}(x(0))v$ may have a component along $R(x(t)).$ The projection $\pi$ then removes this component.  
\end{remark}

 \begin{proposition}\label{lem:compa}
Let  $z = (x,\tau),$ with $\tau \neq 0,$ be  a critical point  of an unperturbed action functional. Assume that both the corresponding Reeb  and Hamiltonian chords are non-degenerate. 
  Then with the compatible trivializations chosen above, we have
\[
\mu_{\rm RS}(x)  - \frac{1}{2}\chi(z)=\mu_{\rm RS}^{\rm tr}(x).
\]
In particular,  
        \[
        \lvert \mu_{\rm RS}(x) - \mu_{\rm RS}^{\rm tr}(x)\rvert \leq \frac{1}{2}.
        \] 
\end{proposition}

\begin{proof}
By the choice of $\tilde{\Phi},$ its extension $\Phi$ is homotopic to $\mathfrak{T}$ through symplectic trivializations sending both endpoint tangent planes of $L$ to $V.$ By homotopy invariance,   we can compute the full Robbin--Salamon index using $\Phi$ instead of $\mathfrak{T}.$

We identify, using the Liouville flow, a small neighborhood of $\Sigma$ with $(1-\varepsilon, 1+\varepsilon)\times \Sigma,$ so that $\lambda = r \alpha$ and $\Sigma = \{ r =1\}. $ Since $X_L$ is tangent to $L,$ we have  $L=(1-\varepsilon, 1+\varepsilon) \times \mathcal{L}$ in this neighborhood. We now set $H_0=r-1$ and $H_s = (1-s)H_0 + sH, \; 0 \leq s \leq 1.$ Then $X_{H_s}=( (1-s)+s\alpha(X_H))R$ along $\Sigma.$ Therefore, we may reparametrize $x$ as a Hamiltonian chord of $H_s,$ with the same orientation. For each $s \in [0,1],$ let $(\hat{x}_s, \hat{\tau}_s)$ be such  a chord, depending smoothly on $s.$ Note that $(\hat{x}_1 , \hat{\tau}_1)=(x,\tau).$ 
For each $s,$ define
\[
\Psi_s(t):=\Phi_{\hat{x}_s(t)}\circ D\phi_{H_s}^{\hat{\tau}_s t }(x(0)) \circ \Phi_{x(0)}^{-1}, \quad  0\leq t \leq 1.
\]
Then $\mu_{\rm RS}(x)=\mu_{\rm RS}(\Psi_1 V, V).$ On the other hand, since the  flow of $H_0$ is $(r, z) \mapsto (r, \phi_R^t(z)),$ we have $\Psi_0(t)=A ( \Psi_0^{\rm tr}(t) \oplus I_2)A^{-1},$
where  $A$ is as in \eqref{eqLAlinser} and 
\[
\Psi_0^{\rm tr}(t):=\tilde{\Phi}_{\hat{x}_0(t)}\circ D\phi_{H_0}^{\hat{\tau}_0 t}(x(0))|_{\ker \alpha_{x(0)}} \circ \tilde{\Phi}_{x(0)}^{-1}.
\]
By naturality and the product property   we have  $\mu_{\rm RS}^{\rm tr}(x)=\mu_{\rm RS}(\Psi_0 V, V).$ 

For each $s,$ define the Lagrangian plane 
\[
E_s:= D\phi_{H_s}^{\hat{\tau}_s }(x(0)) T_{x(0)}L \subset T_{x(1)}M.
\]
For every $s \in [0,1],$ we have $\Psi_s(0)V=V$ and $\Psi_s(1)V=\Phi_{x(1)}E_s.$ 
For each $s,$ we   follow the path $t \mapsto \Psi_s(t)V, \; 0 \leq t \leq 1,$ and then $\Phi_{x(1)}E_u,\; s \leq u \leq 1.$  These concatenated paths form a homotopy with the two Lagrangian planes $V$ and $\Phi_{x(1)}E_1$ fixed. Then homotopy invariance and additivity   give
\[
\mu_{\rm RS}(\Psi_1 V, V)=\mu_{\rm RS}(\Psi_0 V, V) +\mu_{\rm RS}( s \mapsto \Phi_{x(1)}E_s, V).
\]

We compute  $\mu_{\rm RS}( s \mapsto \Phi_{x(1)}E_s, V)=\mu_{\rm RS}(s \mapsto E_s, T_{x(1)}L).$
To apply the localization property of the Robbin--Salamon index  \cite[Theorem 2.3]{RS93}, we shall choose linear symplectic   coordinates on $T_{x(1)}M$ identifying $T_{x(1)}L$ with $\R^2 \times \{ 0 \}$ in which     the planes $E_s$ are represented by the graphs of symmetric matrices $Q_s.$  
Since each $H_s$ is time-independent, 
\[
E_s \cap T_{x(1)}\Sigma = D\phi_{H_s}^{\hat{\tau}_s}(x(0))T_{x(0)}\mathcal{L}.
\]
For each $s,$ by the identity $X_{H_s} =\alpha(X_{H_s})R$ on $\Sigma$ we find a smooth function $\theta_s,$ defined on a neighborhood of  $x(0) $ in $\Sigma,$ such that 
\[
\phi_{H_s}^{\hat{\tau}_s}(z) = \phi_R^{\theta_s(z)}(z), \quad \theta_s( x(0))=\hat{\tau}_0,
\]
from which we get
\begin{align*}
D\phi_{H_s}^{\hat{\tau}_s}(x(0))v&=D\phi_R^{\theta_s(x(0))}(x(0))v + d\theta_s(x(0))v \cdot R\left( \phi_R^{\theta_s(x(0))} (x(0))\right) \\
&=D\phi_R^{\hat{\tau}_0}(x(0))v + d\theta_s(x(0))v \cdot R(  x(1))
\end{align*}
for every $ v \in T_{x(0)}\mathcal{L}.$ Since the first term on the right-hand side belongs to $\ker\alpha_{x(1)},$ we have
\[
\pi ( E_s \cap T_{x(1)}\Sigma) = D\phi_R^{\hat{\tau}_0}(x(0))T_{x(0)}\mathcal{L} =:\ell ,
\]
where $\pi \colon T\Sigma \to \ker\alpha$ denotes the projection along the Reeb direction. In particular, the line $\ell$ is independent of $s.$

Let $\mathfrak{l}  \subset \ker \alpha_{x(1)}$ be a line transverse to both $\ell$ and $T_{x(1)}\mathcal{L},$ and set $W:=\mathfrak{l} \oplus \R R(x(1)),$ which is a Lagrangian complement to $T_{x(1)}L.$ Observe that it is also a Lagrangian complement to every $E_s.$ Indeed, if $v \in E_s \cap W,$ then $v \in T_{x(1)}\Sigma,$ and hence $\pi(v) \in \ell \cap \mathfrak{l} = \{ 0 \}.$ Thus $v \in \R R(x(1)).$ Since $\lambda|_L=0$ and $\lambda(X_{H_s})>0$ on  $\Sigma,$ we have $X_{H_s}( x(0)) \notin T_{x(0)}L,$ from which we get
\[
X_{H_s}(x(1)) \notin E_s = D\phi_{H_s}^{\hat{\tau}_s}(x(0))T_{x(0)}L.
\]
It follows that $R(x(1)) \notin E_s,$ and hence $v =0. $ This proves the claim.

We now take $v \in T_{x(1)}\mathcal{L}$ and  $w \in \mathfrak{l}$ such that $\omega(v,w)=1. $ Since $W$ is a Lagrangian complement to both $T_{x(1)}L$ and $E_s,$ there exists a unique linear map $\mathcal{R}_s \colon T_{x(1)}L \to W$ such that $E_s=\{ u + \mathcal{R}_s(u) \mid u \in T_{x(1)}L\}.$ Denote by $Q_s$ the matrix representation of $\mathcal{R}_s$  with respect to the basis $\{v, X_L(x(1))\}$ of $T_{x(1)}L$ and $\{w, R(x(1))\}$ of $W.$ It is straightforward to check that  $Q_s$ is symmetric, so we may write
\[
Q_s = \begin{pmatrix} a_s & b_s \\ b_s & c_s \end{pmatrix}
\]
 Then $E_s$ can be  written as
 \begin{equation}\label{eq:spanes}
 E_s = {\rm span}\{ v+ a_s w + b_s R( x(1)), X_L(x(1))+ b_s w + c_s R(x(1))\}.
 \end{equation}
It follows that $E_s \cap T_{x(1)}\Sigma = \R ( v+ a_s w + b_s R( x(1))).$ Then 
\[
\ell = \pi ( E_s \cap T_{x(1)}\Sigma) =  \R ( v + a_s w),
\]
from which we see that $a_s$ is independent of $s.$ Write $a=a_s.$ Since $T_{x(1)}\mathcal{L} = \R v,$ non-degeneracy as a Reeb chord implies that $ a \neq 0. $ 
The quadratic form corresponding to the matrix $Q_1$ is 
\[
a \xi^2 + 2b_1 \xi \eta + c_1 \eta^2=a ( \xi + \tfrac{ b_1}{a}\eta)^2 +  ( c_1 - \tfrac{b_1^2}{a})\eta^2,
\]
implying that 
\[
{\rm sign}Q_1 = {\rm sign} (a) + {\rm sign}( c_1 - \tfrac{b_1^2}{a}).
\]

From the expression of the Hamiltonian flow of $H_0,$ we have
\[
D\phi_{H_0}^{\hat{\tau}_0}(x(0))X_L(x(0))=X_L(x(1)).
\]
Thus $X_L(x(1)) \in E_0. $ On the other hand, since $X_L(x(1)) \in T_{x(1)}L,$ the definition of the linear map $\mathcal{R}_0$ provides $X_L(x(1))+ \mathcal{R}_0(X_L(x(1))) \in E_0,$ from which we obtain $\mathcal{R}_0(X_L(x(1))) \in E_0 \cap W = \{0\}.$ Consequently, $b_0=c_0=0,$ and
\[
Q_0 = \begin{pmatrix} a & 0 \\ 0 & 0 \end{pmatrix}
\]
 It follows that ${\rm sign}Q_0 = {\rm sign} (a).$ Therefore,
\[
\mu_{\rm RS}( s \mapsto \Phi_{x(1)}E_s, V)= \frac{1}{2}( {\rm sign}Q_1 - {\rm sign} Q_0) = \frac{1}{2} {\rm sign}( c_1 - \tfrac{b_1^2}{a}).
\]

  We now find a suitable family of Hamiltonian chords used to evaluate the correction term. 
  Since $x$ is non-degenerate as a Hamiltonian chord, the map 
\[
T_{x(0)}L \to T_{x(1)}M/T_{x(1)}L, \quad v \mapsto [ D\phi_H^{\tau}( x(0))v]
\]
is an isomorphism. Set $\tau_s = \tau+s$ for sufficiently small $\lvert s \rvert.$ By the implicit function theorem, we find a smooth family $p_s \in L$   such that  $p_0 = x(0)$ and  $q_s:=\phi_H^{\tau_s}(p_s) \in L.$ Define
\begin{equation}\label{eq:xsests}
x_s(t):=\phi_{H}^{\tau_s t}(p_s), \quad e_s:=H(p_s), \quad T_s:=\lvert \tau_s \rvert.
\end{equation}
Then $H(x_s(t))=H(p_s)=e_s$ and $T_0'={\rm sign}(\tau) \tau_0'={\rm sign}(\tau) \neq 0.$ Assume by contradiction that  $e_0'=0.$ Then $ dH(x(0))p_0' =0.$ Thus,  
\[
p_0' \in T_{x(0)}L \cap \ker dH( x(0)) = T_{x(0)}L \cap T_{x(0)}\Sigma = T_{x(0)}\mathcal{L}.
\]
In a similar way, $q_0' \in T_{x(1)}\mathcal{L}.$ Since   
\begin{equation}\label{q0prime}
q_0'=D\phi_H^{\tau}(x(0))p_0'+   \tau_0' X_H(x(1))  \in T_{x(1)}L ,
\end{equation}
this implies that  
\[
D\phi_R^{\hat{\tau}_0}(x(0))p_0'=\pi  ( D\phi_H^{\tau}(x(0))p_0') \in \ell \cap T_{x(1)}\mathcal{L}=\{0\},
\]
where the intersection is trivial since $x$ is     non-degenerate as a Reeb chord. Since $D \phi_R^{\hat{\tau}_0}$ is invertible,  we conclude that  $p_0'=0.$ Thus,   $\tau_0' X_H (x(1)) \in T_{x(1)}\mathcal{L},$ contradicting $\tau_0'\neq 0$ and $X_H \notin T_{x(1)}L.$ Consequently, $e_0'\neq 0.$ So, using this family, we can evaluate   the correction term $\chi(z)$ as in \eqref{eq:crorr}.

It remains to show that ${\rm sign}( c_1 -  {b_1^2}/{a}) = \chi(z).$  To this end, we first observe
\[
dH(x(1))( D\phi_H^\tau (x(0))p_0')=dH(x(0))p_0'=e_0'
\]
See \eqref{eq:xsests}. Since $X_L \pitchfork \Sigma,$ there exist unique $u \in T_{x(1)}\Sigma$ and $k \in \R$ such that 
$D\phi_H^\tau(x(0))p_0' = u + k X_L(x(1)),$ 
so that  $dH(x(1))\left(D\phi_H^\tau(x(0))p_0'\right) =  k \mathfrak{c}, $ 
where $\mathfrak{c}:=dH(x(1))(X_L(x(1)))=\alpha(X_H)(x(1))>0.$ This implies that $k = e_0'/\mathfrak{c}.$ Since $D\phi_{H}^\tau( x(0))p_0' \in D\phi_H^\tau(x(0)) T_{x(0)}L=E_1,$ there exists   $d \in \R$ such that
\[
D\phi_{H}^\tau( x(0))p_0' = d\left(v+ a  w + b_1 R( x(1))\right)+  \frac{e_0'}{\mathfrak{c}}\left( X_L(x(1))+ b_1 w + c_1 R(x(1))\right).
\]
Since $T_{x(1)}L={\rm span}\{v, X_L(x(1))\}$ and $X_H(x(1))=\mathfrak{c} R(x(1)),$ together with  \eqref{q0prime}, we find
\[
 a d + \frac{e_0'b_1}{\mathfrak{c}}=0, \qquad  b_1 d + \frac{ e_0' c_1}{\mathfrak{c}}+ \mathfrak{c} \tau_0'=0,
\]
from which we get
\[
e_0'\left( c_1 - \frac{b_1^2}{a}\right) + \mathfrak{c}^2 \tau_0'=0, \qquad c_1 - \frac{b_1^2}{a}=-\mathfrak{c}^2\frac{\tau_0'}{e_0'}.
\]
Since $\mathfrak{c}>0,$ $e_0'\neq0$ and $\tau_0'\neq0,$ it follows that
\[
{\rm sign}\left( c_1 - \frac{b_1^2}{a}\right)  = {\rm sign}\left( - \frac{\tau_0'}{e_0'}\right) = {\rm sign}\left( - \frac{e_0'}{\tau_0'}\right)={\rm sign}(\tau){\rm sign}\left(-\frac{e_0'}{T_0'}\right)=\chi(z).
\]
This completes the proof.
\end{proof}

\begin{remark}
    The preceding result does not depend on the product structure of $M$ or on $L$ being a plane: The same argument applies to Hamiltonian chords on contact type hypersurfaces in more general exact symplectic manifolds, provided that $\lambda|_L=0$ near the hypersurfaces and the trivializations are compatible as above.

\end{remark}



\section{Proof of Theorem~\ref{thm:ccorbits}}\label{sec:proofsection}
To prove the theorem, we first reformulate it in the language of Lagrangian submanifolds. As before, for the elliptic regularized system, we use the coordinates $(x_1, x_2, y_1, y_2)$ with $x_2 \in S^1=\R/2\pi\Z.$ The collision loci  over the first primary (Earth) and the second primary (Moon) are the Lagrangian submanifolds 
\[
    L^E_\mathrm{col}:= \{x_1=0,\;x_2=0\}   , \quad L^M_\mathrm{col}:= \{x_1=0, \;x_2=\pi\}  
\]
of $M = S^1 \times \R^3.$ Therefore, we can describe a consecutive collision orbit  as a Hamiltonian trajectory with endpoints in   $L^E_\mathrm{col}$ or $L^M_\mathrm{col}$.  

To describe symmetric consecutive collision orbits,  consider the   anti-symplectic involution 
\[
    \rho \colon  \R^4 \to \R^4, \qquad (q_1,q_2,p_1,p_2) \mapsto (q_1,-q_2,-p_1,p_2)
\]
of the unregularized PCR3BP. As in the Birkhoff regularization (see \eqref{eq:twoBirkhoffrho}), it admits two anti-symplectic lifts to the regularized system via elliptic coordinates:
\[
\rho_a(x_1, x_2, y_1, y_2)=(x_1, -x_2, -y_1, y_2), \quad \rho_b(x_1, x_2, y_1, y_2)=(-x_1, x_2, y_1, -y_2).
\]
 We   focus   on $\rho_a.$  Below, by a symmetric consecutive collision orbit, we mean a consecutive collision orbit $x \colon [0,T] \to M $ satisfying $x(t)=\rho_a(x(T-t)),$   $t \in [0,T].$ 

 As in the Birkhoff regularization, the  fixed point set  ${\rm Fix}(\rho_a)$  is  a Lagrangian submanifold with two connected components,  denoted by $L_{\rho_a}^E $ and  $L_{\rho_a}^M,$ near $x_2=0$ and $x_2=\pi,$ respectively.   Explicitly,
 \[
    L_{\rho_a}^E= \{x_2=0, \; y_1=0\}, \qquad L_{\rho_a}^M= \{x_2=\pi, \; y_1=0\}.
\] 
We shall apply the argument  from Section \ref{sec:basiscsRFH} with $(L_1, L_2)=(   L_{\rm col}^E,  L_{\rho_a}^E).$
The hypersurface is taken to be $\Sigma=\Sigma_c$ for $c\in (c_1,c_2)$ and $\Sigma=\tilde{\Sigma}_c$ for $c $ slightly above $c_2.$
The same argument applies to the case $(L_1,L_2)=( L_{\rm col}^M, L_{\rho_a}^M).$


Define $\tilde{F}(x_1,y_1):=(\pi/4)(x_1^2 +y_1^2).$ Since $\phi_{\tilde{F}}^1(x_1, x_2, y_1,y_2)=(y_1, x_2, -x_1, y_2),$ 
\[
(\phi_{\tilde{F}}^1)^{-1}(L_{\rm col}^E)=L_{\rho_a}^E.
\]
Consider the compact set $\mathcal{K}=\{ \phi_{\tilde{F}}^t(p) \mid p \in \Sigma,\; 0\leq t \leq 1\}.$ Let $\ell \colon M \to [0,1]$ be a compactly supported smooth function such that $\ell \equiv 1$ in a small neighborhood of $\mathcal{K}.$ Define a compactly supported Hamiltonian 
\[
F \colon [0,1] \times M \to \R, \quad F(t,p):=\ell(p)\tilde{F}(p).
\]
By construction, we have
\begin{equation}\label{eq:compareKandiltdeK}
( \phi_F^1)^{-1}( L_{\rm col}^E) \cap \Sigma = L_{\rho_a}^E \cap \Sigma , \qquad \phi_F^t(p)=\phi_{\tilde F}^t(p), \quad p \in \Sigma, \; 0 \leq t \leq 1.
\end{equation}
We choose the cut-off function $\ell$ and all the auxiliary Floer data equivariantly with respect to  $\sigma=\sigma_{\rm el}(-1, \cdot).$  Since  $L_{\rm col}^E$ is  contained in $\{x_2=0\}=\{0\} \times \R^3 \subset S^1 \times \R^3,$ an element of $\mathscr{P}(M, L_{\rm col}^E)$ belongs to the component   $[{\rm pt.}]$ if and only if its projection to the $S^1$ is contractible.  Recall that,   for a critical point $(x,\tau)$ of the  perturbed action functional $\mathcal{A}_H^F$, the second half of $x$ follows the  flow of $F.$ But   by \eqref{eq:compareKandiltdeK}, its $x_2$-coordinate remains constant.   It follows that $x$ belongs to   $[{\rm pt.}]$ if and only if its first half,  from  $L_{\rm col}^E$ to  $L_{\rho_a}^E,$ projects to $S^1$ as a contractible loop.  Thus, 
we can  describe a half of a symmetric consecutive collision orbit around Earth as a Hamiltonian trajectory   that starts in $L^E_\mathrm{col}$ and ends in $L_{\rho_a}^E.$ The other half is obtained by reflection with respect to $\rho_a.$
The same holds for the Moon.

\medskip

\subsection{The subcritical case}  We first consider  energies $c \in (c_1, c_2).$ We search for nonconstant Reeb chords on $\Sigma_c$ from $\Sigma_c \cap L_{\rm col}^E$ to $\Sigma_c \cap L_{\rho_a}^E,$ following the argument in \cite[Section 4]{kang2025consecutive}. 
Let $\lambda_{\rm el}$ be the contact form  on $\Sigma_c$ as in \eqref{eq:twpoconact}. Its Reeb flow is a positive time-reparametrization of the Hamiltonian flow of $K_c$ on $\Sigma_c.$

Consider the one-form 
\[
\lambda_0 = y_2 dx_2 + x_1 dy_1 - y_1 dx_1
\]
on $S^1 \times \R^3$ with symplectic form   $d\lambda_0 = 2 dx_1 \wedge dy_1 + dy_2 \wedge dx_2.$   We denote by $S^2$ the unit sphere in $\R^3$ with coordinates $(x_1,y_1, y_2),$ so that $S^1 \times S^2$ is a closed hypersurface in $M=S^1 \times \R^3.$ It is of contact type with $\lambda_{\rm st}:=\lambda_0|_{S^1 \times S^2},$  and the involution $\sigma=\sigma_{\rm el}(-1,\cdot)$ restricts to an involution preserving $\lambda_{\rm st}.$

Suppose that  $c>c_1$ is sufficiently close to $c_1.$ 

\begin{lemma}\label{lem:hatPsi}
    There is  a $\sigma$-equivariant co-orientation preserving  contactomorphism
\[
    \hat{\Psi}  \colon  (\Sigma_c, \ker \lambda_\mathrm{el}) \to (S^1\times S^2, \ker \lambda_\mathrm{st}), 
\]
which    maps the Legendrians $\Sigma_c\cap L^E_\mathrm{col}$ and $\Sigma_c\cap L_{\rho_a}^E$ to meridians (with respect to the $y_2$-axis) in $\{0\}\times S^2,$  denoted by  $\Lambda_0$ and $\Lambda_1,$ respectively.   
\end{lemma}
\begin{proof}The   contactomorphism $\Psi\colon  (\Sigma_c, \ker \lambda_\mathrm{el}) \to  (\Sigma_c, \ker \lambda_B)$ from  Proposition \ref{prop:elToBirk} leaves the two Legendrians $\Sigma_c \cap L_{\rm col}^E$ and $\Sigma_c \cap L_{\rho_a}^E$ invariant. 
 Let $\Phi$ be the  equivariant symplectomorphism, defined in Lemma \ref{lem:Phi}. By definition we have
\[
    \Phi\left(\Sigma_c\cap L^E_\mathrm{col}\right)= \Lambda_c^E, \qquad  \Phi\left(\Sigma_c\cap L_{\rho_a}^E\right)= \Gamma_c^E.
\]
Consequently, the composition $\Phi \circ \Psi \colon ( \Sigma_c, \ker \lambda_{\rm el}) \to \left( \Phi(\Sigma_c), \ker ( (\Phi^{-1})^*\lambda_B)\right)$ is an  equivariant contactomorphism satisfying 
\[
    ( \Phi \circ \Psi)\left(\Sigma_c\cap L^E_\mathrm{col}\right)= \Lambda_c^E, \qquad  ( \Phi \circ \Psi)\left(\Sigma_c\cap L_{\rho_a}^E\right)= \Gamma_c^E.
\]
Composing this with the equivariant contactomorphism  constructed in  \cite[Proposition~3.3]{kang2025consecutive} gives  the desired  contactomorphism. This completes the proof. 
\end{proof}

Let $\hat{\Psi}$ be the contactomorphism obtained in the preceding lemma.
The involution $\rho_a$ induces the anti-contact involution $\hat{\Psi} \circ \rho_a \circ \hat{\Psi}^{-1}$ on $(S^1 \times S^2, \ker \lambda_{\rm st}),$  where $\Lambda_1$ is a connected component of its fixed point set. 

Consider  the contact form $(\hat{\Psi}^{-1})^*\lambda_{\rm el}$ on $S^1 \times S^2.$ By construction, it is a positive multiple of $\lambda_{\rm st}.$ Therefore, we may write  $(\hat{\Psi}^{-1})^*\lambda_{\rm el} = e^g \lambda_{\rm st} $ for some $\sigma$-invariant $g \in C^{\infty}(S^1 \times S^2) .$  Let $Z$ be the Liouville vector field associated with $\lambda_0,$ defined by  $\iota_{Z} d\lambda_0 =\lambda_0.$ Explicitly, $Z= \frac{1}{2}( x_1 \partial_{x_1} + y_1 \partial_{y_1}) + y_2 \partial_{y_2}.$ Every nonconstant flow line of $Z$ intersects the hypersurface $S^1 \times S^2$ precisely once and transversely. We then obtain an equivariant embedding
\[
\iota_g \colon S^1 \times S^2 \hookrightarrow S^1 \times \R^3, \qquad p \mapsto \phi_Z^{g(p)} (p).
\]
We set $\Sigma_g:=\iota_g(S^1 \times S^2).$  In the coordinates defined by the flow of $Z,$ the hypersurface $\Sigma_g$ is the graph of $g.$ Note that $\iota_g^*\lambda_0 = e^g\lambda_{\rm st},$ so  the composition $\iota_g \circ \hat{\Psi}$ pulls back the contact form $ \lambda_0|_{\Sigma_g}$ to $\lambda_{\rm el}.$

Let $\mathcal{L}_0, \mathcal{L}_1 \subset \{ 0 \} \times \R^3$ be the Lagrangian planes in $(S^1 \times \R^3, d\lambda_0)$ containing the meridians $\Lambda_0, \Lambda_1,$ respectively. Then $\mathcal{L}_j \cap (S^1 \times S^2) = \Lambda_j, \;j=0,1.$ Applying a rotation in the $(x_1, y_1)$-plane, which preserves $\lambda_0$ and commutes with $\sigma$, if necessary, we may assume that $\mathcal{L}_0=\{x_2=y_1=0\}.$
Since the vector field $Z$ is tangent to these planes, we have 
\[
\iota_g(\Lambda_j) = \Sigma_g \cap \mathcal{L}_j, \quad j=0,1.
\]
Since $Z$ is transverse to $\Sigma_g,$ the above intersections are transverse.

We now apply the perturbed construction from Section \ref{sec:basiscsRFH} to $(S^1 \times \R^3,d\lambda_0),$ with hypersurface $\Sigma_g$ and Lagrangians $\mathcal{L}_0$ and $\mathcal{L}_1.$ As in \cite[Section 4.2]{kang2025consecutive}, we choose a compactly supported $\sigma$-invariant Hamiltonian $G$ by cutting off the Hamiltonian generating a rotation in the $(x_1, y_1)$-plane which maps $\mathcal{L}_1$ to $\mathcal{L}_0.$ Thus, we have
\[
( \phi_G^1)^{-1}(\mathcal{L}_0)\cap \Sigma_g = \mathcal{L}_1 \cap \Sigma_g.
\]
Then by equivariant continuation to the zero perturbation and \cite[Proposition 4.1]{kang2025consecutive}  
\[
{\rm RFH}_k^{\Z_2}(S^1 \times \R^3, \Sigma_g, \mathcal{L}_0, \mathcal{L}_1, [{\rm pt.}])  \cong {\rm RFH}_k^{\Z_2}(S^1 \times \R^3, \Sigma_g, \mathcal{L}_0,   [{\rm pt.}])  \cong \Z_2
\]
for every $k \in \Z.$
We now  extend this argument to all $c \in (c_1, c_2).$

\begin{lemma}\label{leam:scondfi} For every $c \in (c_1, c_2),$ there is   an equivariant contactomorphism from  $\left( \Sigma_{c }, \ker(\lambda_{c }|_{\Sigma_{c }})\right)$ to $(S^1 \times S^2, \ker \lambda_{\rm st})$ 
which maps $\Sigma_{c } \cap L_{\rm col}^E$ and $\Sigma_{c } \cap L_{\rho_a}^E$ to $\Lambda_0$ and $\Lambda_1,$ respectively. 
\end{lemma}
\begin{proof}
Fix $c_* \in (c_1, c_2)$ and choose $c_0 \in (c_1, c_*)$ sufficiently close to $c_1,$ so that the preceding argument applies. Since $K_c|_{L_{\rm col}^E}=\frac{1}{2}( y_1^2 +y_2^2-1), $ we have
  \[
  \Sigma_c \cap L_{\rm col}^E = \Sigma_{c_0} \cap L_{\rm col}^E, \qquad c \in [c_0,c_*].
  \]
   Since the gradient $\nabla K_c$ (taken with respect to the standard product metric) is nowhere zero and normal to $\Sigma_c$ for every $c \in [c_0,c_*],$ the flow of a  suitable $c$-dependent multiple of $\nabla K_c$ gives a smooth family of diffeomorphisms $\Xi_c \colon \Sigma_{c_0} \to \Sigma_c,   c \in [c_0, c_*],$ such that $\Xi_{c_0}={\rm id}$ and $\Xi_c$ is equivariant with respect to both $\sigma$ and $\rho_a.$ Moreover, since $\partial_c K_c=0$ along $L_{\rm col}^E,$ each $\Xi_c$ fixes the common Legendrian circle $\Sigma_{c_0} \cap L_{\rm col}^E.$  
Choose a smooth family of $\sigma$-invariant and $\rho_a$-anti-invariant primitives $\lambda_c$ of $\omega,$ for $c \in [c_0,c_*],$ such that 
\[
 \lambda_c(X_{K_c})|_{\Sigma_c}>0, \quad  \lambda_c|_{L_{\rm col}^E}=0, \quad \lambda_{c_0}|_{\Sigma_{c_0}}=\lambda_{\rm el}.
 \] 
We then obtain a smooth family of contact forms $\alpha_c:=\Xi_c^*(\lambda_c|_{\Sigma_c})$ on $\Sigma_{c_0}.$   Note that $\Sigma_{c_0} \cap L_{\rm col}^E$ is Legendrian for every $\alpha_c.$ 
By the Moser trick as in the proof of Proposition \ref{prop:elToBirk}, we find a family of diffeomorphisms 
\[
\psi_c \colon \Sigma_{c_0} \to \Sigma_{c_0} , \quad \psi_{c_0}={\rm id}
\]
such that each $\psi_c^* \alpha_c$ is a positive multiple of $\alpha_{c_0}.$ Each $\psi_c$ preserves  $\Sigma_{c_0}\cap L_{\rm col}^E$ and commutes with both $\sigma$ and $\rho_a.$ Therefore, for every $ c \in [c_0,c_*],$ the composition
\[
\Xi_c \circ \psi_c \colon ( \Sigma_{c_0}, \ker (\lambda_{c_0}|_{\Sigma_{c_0}})) \to (\Sigma_c, \ker ( \lambda_c|_{\Sigma_c}))
\]
  is an equivariant    contactomorphism which preserves $\Sigma_{c_0} \cap L_{\rm col}^E$  and maps $\Sigma_{c_0} \cap L_{\rho_a}^E$ to $\Sigma_c \cap L_{\rho_a}^E.$   We therefore have an equivariant contactomorphism 
  \[
  \hat{\Psi} \circ ( \Xi_{c_*} \circ \psi_{c_*})^{-1} \colon \left( \Sigma_{c_*}, \ker(\lambda_{c_*}|_{\Sigma_{c_*}})\right) \to (S^1 \times S^2, \ker \lambda_{\rm st})
  \]
which maps $\Sigma_{c_*} \cap L_{\rm col}^E$ and $\Sigma_{c_*} \cap L_{\rho_a}^E$ to $\Lambda_0$ and $\Lambda_1,$ respectively. Since $c_*$ was arbitrary,  this completes the proof.
\end{proof}

For each $ c \in (c_1, c_2),$ let $\hat{\Phi}$ be the contactomorphism from Lemma \ref{leam:scondfi} and let $\iota_g$ and $\Sigma_g$ be as above. 
Let $H$ be an adapted Hamiltonian associated with $\Sigma_g$ as in Section \ref{sec:basiscsRFH}. 
If the  perturbed  action functional $\mathcal{A}_H^G$ is Morse, then the preceding discussion implies that there are infinitely many critical points. In the degenerate case, choose  compactly supported $\sigma$-invariant time-dependent Hamiltonians $G _n \to G$ for which   $\mathcal{A}_{H }^{G_n}$ is Morse. By the same reasoning as above, for each $n,$ the associated RFH is isomorphic to $\Z_2$ in every degree. Applying the approximation and  spectral invariant argument in \cite[Corollary 4.2]{kang2025consecutive}, based on \cite[Section~4.3]{frauenfelder2019a}, we conclude that $\mathcal{A}_H^G$ also has infinitely many critical points. 
The critical points with $\tau=0$ are the two constant paths at  $\Sigma_g \cap \mathcal{L}_0 \cap \mathcal{L}_1  = (\iota_g \circ \hat{\Phi})( 0,0,0,\pm 1).$
 It follows that  there are infinitely many critical points with $\tau \neq 0.$ Applying $(\iota_g \circ \hat{\Phi})^{-1}$ to their Hamiltonian segment yields, up to reparametrization, infinitely many nonconstant chords of $K_c$ on $\Sigma_c$ from $L_{\rm col}^E$ to $L_{\rho_a}^E.$  Since $\rho_a$ reverses the  flow of $K_c,$ preserves $L_{\rm col}^E$ and fixes $L_{\rho_a}^E,$ it gives a bijection between negative- and positive-time Hamiltonian chords. Hence there are infinitely many positive-time chords with these boundary conditions. 
 As in \cite[Corollary 4.2]{kang2025consecutive} these chords imply that either there exists a symmetric periodic collision orbit, whose iterations may account for the infinitely many critical points,  or there are infinitely many distinct symmetric consecutive collision orbits. The same argument applies to $L_{\rm col}^M$ and $L_{\rho_a}^M.$ 
This proves the theorem for the subcritical  case.

 \begin{remark}\label{rem:z2compute} For the supercritical argument, we show that  for every $c \in (c_1, c_2)$
\[
{\rm RFH}_k^{\Z_2}(M, \Sigma_c, L_{\rm col}^E, L_{\rho_a}^E,[{\rm pt.}]) \cong {\rm RFH}_k^{\Z_2}(M, \Sigma_c, L_{\rm col}^E,  [{\rm pt.}]) \cong \Z_2, \quad k \in \Z .
\]

Fix $d \in (c_1, c_2)$  and choose $\theta<1$ sufficiently close to $1$ so that the Liouville vector field $\hat{Y}$ from Remark \ref{remark:slightlyaboveL2} is positively transverse to $\Sigma_d.$ Define 
\[
E \colon M \to \R, \quad (x_1, x_2, y_1, y_2) \mapsto \frac{1}{2}( x_1 ^2 + y_1^2) + y_2^2
\]
and set $\Sigma_*:=E^{-1}(1).$ Observe that, after taking $\theta<1$ closer to $1,$ we may assume that $\hat{Y}$ is positively transverse to $\Sigma_*$ as well. Indeed, we compute
\[
Y_1(E)= \frac{1}{2} x_1 \tanh(2x_1) + \tanh^2 ( 2x_1) y_1^2 + 2y_2^2 \geq 0,
\]
which vanishes precisely when $x_1=y_2=0.$ Therefore, we have $Y_1(E)=0$ on $\Sigma_* $ precisely along the circles $\{(0,x_2,\pm \sqrt{2} ,0) \mid x_2 \in S^1 \},$ contained in $\mathcal{V} = \{ x_1^2 +p_2^2 <s_-\}.$  See \eqref{eq:VandUI}. Therefore, on the complement $\Sigma_* \setminus \mathcal{V},$ we conclude $\hat{Y}(E)>0,$ provided $\theta<1$ is sufficiently close to $1.$ On $\Sigma_* \cap \mathcal{V},$ we have
\[
\hat{Y}(E)=Y_\theta(E) = \frac{\theta}{2} x_1 \tanh (2x_1) + ( 1- \theta \sech^2(2x_1) ) y_1^2 +2y_2^2>0.
\]
It follows that $\hat{Y}$ is positively transverse to $\Sigma_*.$

We claim that $\Sigma_d$ can be represented by a graph over $\Sigma_*$ in the  coordinates defined by the flow of $\hat{Y}.$ 
Consider the circle $\mathcal{C} =\{ x_1 = y_1 = y_2=0\},$ which consists of rest points of $\hat{Y}.$ One easily checks that $\mathcal{C}$ lies in the interiors of the compact domains bounded by $\Sigma_d$ and $\Sigma_*.$   To describe the nonconstant flow lines of $\hat{Y},$ we set $u = \sinh(2x_1),$ so that along the flow we have
\[
\dot{u} = \varphi_\theta u, \qquad \dot{y}_2 = y_2, \qquad \theta \leq \varphi_\theta \leq 1.
\]
The flow of $\hat{Y}$ is complete. Indeed, these equations show that $u$, and hence $x_1,$ as well as $y_2,$ remain bounded on every bounded time interval. The equation for $y_1$ is linear with a smooth coefficient depending only on $x_1 $ and $y_2,$ so $y_1$ also remains bounded on each such interval. Moreover, $x_2$ takes values in $S^1.$

In backward time,  $u,y_2\to 0,$ and hence $x_1 \to 0.$  Since $p_2=  y_2 + A_1(x_1) \to 0,$ every nonconstant trajectory eventually stays in $\mathcal{V},$ where $\hat{Y} = Y_\theta.$ Thus, on $\mathcal{V},$ the vector field $\hat{Y}$  has no $x_2$-component, so $x_2$ is constant along each trajectory while it remains in $\mathcal{V}.$  Moreover,   $\dot{y}_1 = (1-\theta \sech^2(2x_1))y_1,$ with $1-\theta \sech^2 (2x_1) \geq 1-\theta>0.$ Therefore, every nonconstant trajectory converges in backward time to   $\mathcal{C}.$

In forward time,  we have $\lvert u(t)\rvert \geq e^{\theta t} \lvert u(0)\rvert$ and $y_2(t) = e^t y_2(0)$ for $t \geq 0.$
If $u(0) \neq 0$ or $y_2 (0) \neq 0,$ then $\lvert x_1(t)\rvert$ or $\lvert y_2\rvert,$ respectively, tends to infinity.  
If $u(0)=y_2(0)=0,$ then $x_1(0)=y_2(0)=0$ for all $t,$ so the trajectory remains in $\mathcal{V}$ and $y_1(t)=e^{(1-\theta)t}y_1(0).$ In this case, a trajectory starting outside $\mathcal{C}$ has $y_1(0)\neq0,$ so $\lvert y_1 (t)\rvert \to \infty.$ Thus every trajectory starting outside $\mathcal{C}$ leaves every compact set in forward time.

It follows that each nonconstant flow line meets both $\Sigma_*$ and $\Sigma_d.$ Since $\hat{Y}$ points outward at every intersection, each flow line meets each hypersurface precisely once.   We can therefore define $g(p)$ for each $p \in \Sigma_*$ as the unique time such that $\phi_{\hat{Y}}^{g(p)} (p) \in \Sigma_d.$ Transversality to $\Sigma_d$ implies that $g$ is smooth. Thus, in the coordinates determined by the flow, the hypersurface $\Sigma_d$ is the graph of $g,$ as claimed.

 For $0 \leq s \leq 1,$ the graphs $p \mapsto \phi_{\hat{Y}}^{sg(p)}(p)$ give an equivariant homotopy from $\Sigma_*$ to $\Sigma_d.$ The pullback of $\lambda = \iota_{\hat{Y}}\omega$ is $e^{sg}\lambda|_{\Sigma_*},$ and $\hat{Y}$ is tangent to $L_{\rm col}^E$ and transverse to every graph. Equivariant continuation then gives 
 \[
 {\rm RFH}_k^{\Z_2}( M, \Sigma_d, L_{\rm col}^E, [{\rm pt.}]) \cong {\rm RFH}^{\Z_2}_k(M, \Sigma_*, L_{\rm col}^E, [{\rm pt.}]).
 \]
 To compute the latter, we  consider the symplectomorphism
\[
\Theta \colon (M, \omega) \to (M, d\lambda_0), \quad (x_1, x_2, y_1, y_2) \mapsto (\frac{y_1}{\sqrt{2}}, x_2, \frac{x_1}{\sqrt{2}}, y_2).
\]
A direct computation shows that
\[
\Theta^*\lambda_0 = \iota_{\hat{Y}}\omega + d\ell, \qquad \ell:=\frac{1}{2}\varphi_\theta \tanh(2x_1)y_1 - \frac{1}{2}x_1y_1.
\]
Note that $\Theta(\Sigma_*)=S^1 \times S^2$ and  $\Theta(L_{\rm col}^E)=\mathcal{L}_0.$ Choose $H$ adapted to $\Sigma_*$ with respect to $\iota_{\hat{Y}}\omega$ as in Section \ref{sec:basiscsRFH}.
We then compare the two action functionals 
\begin{align*}
   & \mathcal{A}_H  \colon \mathscr{P}(M, L_{\rm col}^E) \times \R \to \R \; (\text{with respect to $\iota _{\hat{Y}} \omega$}) \\
   & \mathcal{A}_{H \circ \Theta^{-1}}  \colon \mathscr{P}(M, \mathcal{L}_0)  \times \R \to \R \; (\text{with respect to $\lambda_0$})
    \end{align*}
    Since $(\iota_{\hat{Y}}\omega )|_{L_{\rm col}^E}=0$ and $\lambda_0|_{\mathcal{L}_0}=0,$ we can take the auxiliary functions to be the zero function. A direct computation shows that $   \mathcal{A}_{H \circ \Theta^{-1}} (\Theta \circ x, \tau) =   \mathcal{A}_H(x,\tau).$
Transporting the Floer data by   $\Theta $   identifies  the two unperturbed complexes. The map commutes with $\sigma$ and preserves the component $[{\rm pt.}].$  It follows that
\[
{\rm RFH}_k^{\Z_2}(M, \Sigma_* , L_{\rm col}^E, [{\rm pt.}]) \cong {\rm RFH}_k^{\Z_2}(M, S^1 \times S^2, \mathcal{L}_0, [{\rm pt.}]) \cong \Z_2, \;\;k \in \Z.
\]

 \end{remark}

\subsection{The supercritical case}\label{secL:72} We now turn to the case  $c\in (c_2,c_2+\varepsilon),$ with $\varepsilon>0$ sufficiently small, where  we have   the capped energy hypersurface $\tilde{\Sigma}_c$. By Theorem~\ref{thm:perturbation} there exists a globally defined Liouville vector field $Y$ that coincides with  $\hat{Y}$  away from the open sets $\lambda(\mathcal{U}_\mathcal{W}),$  $\lambda \in \Lambda$. Fix $d \in (c_1, c_2) $ and choose the parameters in the definition of $\hat{Y},$   so that $\hat{Y}$ is positively transverse to $\Sigma_d.$ See Remark \ref{remark:slightlyaboveL2}. Since $K_d(z_*(c_2))>0,$ the coordinate neighborhoods $\lambda ( \mathcal{B}_{r_*}(z_*(c))),$ $\lambda \in \Lambda,$   can be taken to be   contained in $\{K_d >0\},$ and hence  outside the compact region bounded by $\Sigma_d,$ after shrinking these neighborhoods and  taking smaller $\varepsilon>0$ if necessary. Since $\overline{\mathcal{U}_{\mathcal{W}}} \subset \mathcal{B}_{r_*}(z_*(c)),$ the modification of $\hat{Y}$  (to $Y$) is supported outside this compact region. Hence, $Y = \hat{Y}$ near $\Sigma_d,$  so $Y$ is positively   transverse to $\Sigma_d$. Moreover,  since  $x_1 = \pm x_*(c_2) \neq 0$ at the equilibria   and $L_{\rm col}^E \subset \{x_1 =0\},$  after shrinking the coordinate neighborhoods if necessary, we may assume they are disjoint from  $L_{\rm col}^E.$   Consequently,   $Y=\hat{Y}$ also near $L_{\rm col}^E.$

As in the subcritical case, we would like to express $\tilde{\Sigma}_c$ as a graph over $\Sigma_d$ using the Liouville flow of $Y.$ However, $Y$ vanishes at the equilibria, and hence this procedure is not straightforward. We resolve  this issue by modifying $Y$ near the equilibria, but keeping it fixed near   $\Sigma_d \cup  \tilde{\Sigma}_c \cup  L_{\rm col}^E.$


  \begin{proposition}\label{prop:weinstein}
With suitable choices of  $\rho, \varrho,  \mathfrak{n} $ and  the cutoff function $\mathfrak{f}$ of Section \ref{sec:4}, 
there is $0<\varepsilon \leq \varrho^3 $  such that for each $c \in (c_2, c_2+\varepsilon)$ the Liouville vector field $Y$ admits a $\Lambda$-equivariant  modification supported in the four   coordinate neighborhoods $\lambda(\mathcal{B}_{r_*}( z_*(c))),$ $\lambda \in \Lambda,$   away from $\Sigma_d \cup \tilde{\Sigma}_c \cup L_{\rm col}^E.$ The modified vector field  $Y'$ is    Liouville   with respect to $\omega.$ Moreover, if a flow line of $Y'$ intersects either $\Sigma_d$ or $\tilde{\Sigma}_c,$ then it intersects each of these hypersurfaces precisely once.
  \end{proposition}
\begin{proof}[Proof Sketch] (For the detailed proof, see Appendix~\ref{appendproof}). As mentioned above the    unmodified Liouville vector field $Y$ has two equilibrium points in each of the four distinct neck regions, one caused by the Weinstein-like vector field $X_c$ (an index-one saddle) and one caused by the radial vector field $\mathcal{Z}_{c,\varrho}$ (a minimum). We construct a   Weinstein cobordism  containing these four pairs. Let $\mathfrak{D}_c$ be the compact domain bounded by $\tilde{\Sigma}_c.$ Then the cobordism is given by  
\[
    W_c:=\left\{z\in \mathfrak{D}_c  :\ h\le \tilde{K}_c(z)\le 0 \right\},
\]
for a suitable regular value $h<0$, with the standard symplectic structure,   Morse function  $\tilde{K}_c$ and the gradient-like vector field is $Y$. The $\Lambda$-action on $W_c$ is free. In the quotient $W_c/\Lambda,$ there is a single minimum and a single index-one saddle, connected by a unique trajectory, up to reprarametrization. The cancellation theorem of Cieliebak and Eliashberg (see Theorem \ref{thm:cancle}) applies to $W_c/\Lambda_c.$   Under this homotopy the saddle point and the minimum degenerate into a single birth-death point and then disappear. Since the closure of the stable manifold of each saddle point lies in its coordinate neighborhoods, the homotopy can be chosen, so that its $\Lambda$-equivariant lift is supported in the four coordinate neighborhoods, away from $\Sigma_d \cup \tilde{\Sigma}_c \cup L_{\rm col}^E.$ After cancellation, the modified Morse function increases at a uniformly positive rate along $Y'$ on $W_c,$ so the flow crosses the cobordism in finite time. Since $Y'=\hat{Y}$ below $W_c $ and  positive transversality imply that every flow intersecting either $\Sigma_d$ and $\tilde{\Sigma}_c$ meets each exactly once.
\end{proof}

By abuse of notation, denote  by $Y$ the   modified Liouville vector field $Y'$ in  Proposition \ref{prop:weinstein}.  Denote by $\phi_Y^t$ its flow. For each $s \in [0,1]$ we define the embedding
\[
\iota_s \colon \Sigma_d \hookrightarrow M=S^1 \times \R^3,\quad \iota_s(p):=\phi_Y^{s f(p)}(p),
\]
where   $f(p)$ denotes the time required for the flow of $Y,$ starting at $p \in \Sigma_d, $  to reach $\tilde{\Sigma}_c.$ By Proposition \ref{prop:weinstein}, the function  $f$ is   smooth and  $\sigma$-invariant. 
Write $\Sigma^s = \iota_s(\Sigma_d),$ so that $\Sigma^0 = \Sigma_d$ and $\Sigma^1=\tilde{\Sigma}_c.$

Set $\lambda = \iota_Y \omega$ and $\alpha_d= \lambda |_{\Sigma_d}.$ Since $\lambda(Y)=0$ and $(\phi_Y^r)^*\lambda = e^r \lambda,$ we have  $\iota_s^*\lambda = e^{sf}\alpha_d.$ Since $e^{sf} \alpha_d$ is a contact form on $\Sigma_d,$ this implies that $\lambda|_{\Sigma^s}$ is a contact form on $\Sigma^s$ for every $s \in [0,1].$ Since $Y$ is tangent to $L_{\rm col}^E,$ its flow preserves $L_{\rm col}^E.$ Moreover,  since 
$\Sigma_d \cap L_{\rm col}^E = \tilde{\Sigma}_c \cap L_{\rm col}^E = \{ (0,0,y_1,y_2) \mid y_1^2+y_2^2=1\},$ 
we have $f =0$ along this Legendrian circle, implying that 
$\iota_s(p)=p$ for all  $p \in \Sigma_d \cap L_{\rm col}^E .$ 
Hence, $\Sigma^s \cap L_{\rm col}^E = \Sigma_d \cap L_{\rm col}^E $ for all $s \in [0,1].$ Since $Y$ is tangent to $L_{\rm col}^E$ and transverse to $\Sigma^s,$ this intersection is transverse. 
By Remark~\ref{rmk:functioanlnunchange}, replacing the original primitive on $\tilde{\Sigma}_c$ and $\Sigma_d$ by the common primitive $\lambda = \iota_Y\omega$ preserves the corresponding action functionals. Thus, we have   
\begin{equation}\label{eq:RFHsuper}
\begin{aligned}
{\rm RFH}^{\Z_2}_k(M, \tilde{\Sigma}_c, L_{\rm col}^E, L_{\rho_a}^E,[{\rm pt.}]) & \cong {\rm RFH}_k^{\Z_2}( M, \tilde{\Sigma}_c, L_{\rm col}^E, [{\rm pt.}])\\
&\cong {\rm RFH}_k^{\Z_2}(M, \Sigma_d, L_{\rm col}^E,[{\rm pt.}])\\
&\cong \Z_2, \qquad k \in \Z,
\end{aligned}
\end{equation}
where the first isomorphism follows from perturbation invariance, the second  from equivariant continuation along the family $(\Sigma^s)_{s \in [0,1]}$ and the last   from Remark \ref{rem:z2compute}. 

By the same argument as in the subcritical case, there are infinitely many chords on $\tilde{\Sigma}_c.$ However, some of these orbits  may pass through the neck and enter the region where we modified the Hamiltonian   in Section~\ref{sec:capping}. To obtain chords of the original Hamiltonian, we use the following index estimate to show that, for any given $N$ and all $c>c_2$ sufficiently close to $c_2,$ at least $N$ of these chords remain in the unmodified region. 
Let $\mathfrak{T} \colon TM \to M \times \R^4$ denote the  global symplectic trivialization   as in  Section \ref{sec:basiscsRFH}, under which    the tangent spaces of the collision Lagrangian $L_{\rm col}^E$ are identified with the fixed Lagrangian subspace   $V = {\rm span} \{ \partial_{y_1}, \partial_{y_2} \}.$

\begin{proposition}\label{prop:iundexesti}
    For every $N \in \N,$ there exist  $\varepsilon_N>0$ and an open set $\mathcal{B}_N,$ which is independent of $c \in (c_2, c_2+\varepsilon_N)$ and consists of four disjoint small balls centered at the equilibria of $K_{c_2},$ such that  for every  $c \in (c_2, c_2+\varepsilon_N),$ every   trajectory $x \colon [0,\tau]\to \tilde{\Sigma}_c$ of $\tilde{K}_c,$  satisfying $x(0) \in L_{\rm col}^E$ and $x([0,\tau]) \cap \mathcal{B}_N \neq \emptyset,$ satisfies
    \[
    \mu_{\rm RS}( \psi V  , V)>N, \qquad \psi(t):= \mathfrak{T}_{x(t)} \circ D\phi_{\tilde{K}_c}^t(x(0)) \circ \mathfrak{T}_{x(0)}^{-1}, \; 0 \leq t \leq \tau.
    \]

\end{proposition}
\begin{proof} The proof   is based on the argument due to Liu and Salom\~ao. We include some details since we will use some arguments used in this proof later. 

Fix $N \in \N.$ Let $c>c_2$ be sufficiently close to $c_2.$  
Let $z_*(c)$ be an equilibrium point corresponding to the second or third Lagrange point. As in \cite[Section 7]{LS25finite} we choose ($c$-dependent) local symplectic  coordinates $\hat{z}=(\hat{x}_1, \hat{x}_2, \hat{y}_1, \hat{y}_2)$  by translating $z_*(c)$ to the origin and then applying a linear symplectic transformation, in which $\omega=d\hat{x}_1\wedge d\hat{y}_1 + d\hat{x}_2 \wedge d\hat{y}_2$ and   $K_c = K_c(z_*(c))+Q+R,$ where
\[
    Q(\hat{x},\hat{y})=\frac{\hat{\lambda}_1}{2}\left(\hat{x}_1^2-\hat{y}_1^2\right) +\frac{\hat{\lambda}_2}{2}\left(\hat{x}_2^2+\hat{y}_2^2\right)
\]
where $\hat{\lambda}_j=\hat{\lambda}_j(c)>0, j=1,2,$ depend smoothly on $c,$   and $R$ vanishes up to second order at the origin.   After shrinking the coordinate neighborhood, denoted by $\mathcal{U}_c,$ if necessary, we may assume that   $\tilde{K}_c=K_c$ on $\mathcal{U}_c$  and that $\mathcal{U}_c \cap L_{\rm col}^E=\emptyset.$ Fix $T>0$ to be chosen below.
Since   $\phi_{K_{c_2}}^t(z_*(c_2))=z_*(c_2)$ for all $t,$ by  continuous dependence of the Hamiltonian flow on the initial condition and energy,  we find a sufficiently small $\varepsilon_N>0$ and  a small $r>0,$ independent of $c \in (c_2, c_2+\varepsilon_N),$  such that 
\[
B_r(z_*(c_2)) \subset \mathcal{U}_c, \qquad \phi_{K_c}^t(p) \in \mathcal{U}_c
\]
for every $c \in (c_2,c_2+\varepsilon_N),$  every $p \in B_r(z_*(c_2)) $ and every $t \in [-T,0].$ We do the same job near the other three equilibria using symmetry.  Denote by $\mathcal{B}_N$ the union of the four  $r$-balls, which is independent of $c \in (c_2, c_2+\varepsilon_N). $ Below, the neighborhoods $\mathcal{U}_c, $ the radius $r$ and $\varepsilon_N>0$  will be chosen sufficiently small, depending on $T.$

Choose a Hamiltonian chord $x \colon [0,\tau] \to \tilde{\Sigma}_c$ starting on $L_{\rm col}^E$ and entering $\mathcal{B}_N$ at time $T_0 \in (0,\tau].$ The preceding discussion tells us that $x ( [T_0-T,T_0])\subset \mathcal{U}_c.$ Moreover,  $T_0 >T $ since otherwise we have $x(0) \in \mathcal{U}_c,$ contradicting $\mathcal{U}_c \cap L_{\rm col}^E  = \emptyset.$ 
Let $\psi(t)$ be as in the statement of the proposition.
By additivity of the Robbin--Salamon index,  we have
\begin{align*}
\mu_{\rm RS}(\psi V  , V) &= \mu_{\rm RS}(\psi V |_{[0,T_0-T]}, V) \\
&\;\;\; +  \mu_{\rm RS}(\psi V |_{[ T_0-T, T_0]  }, V)+\mu_{\rm RS}(\psi V |_{  [T_0,\tau]}, V). 
\end{align*}
Since the Hessian of $\tilde{K}_c$ restricted to the $(y_1,y_2)$-plane is the identity matrix, the computation in the proof of \cite[Proposition 7.1]{LS25finite} shows that every crossing of $\psi(t)V$ with $V$ is positive.  In particular,   the   contributions  from $[0,T_0-T]$ and $[T_0,\tau]$ to the Robbin--Salamon index  are nonnegative.

We claim that $\mu_{\rm RS}(\psi V |_{[T_0-T,T_0]}, V) >N.$
 Define  
\[
\psi_0(t):=\mathfrak{T}_{x(T_0-T+t)} \circ D\phi_{\tilde{K}_c}^t (x(T_0-T)) \circ \mathfrak{T}_{x(T_0-T)}^{-1}, \quad 0 \leq t \leq T.
\]
Since $x([T_0-T,T_0])\subset \mathcal{U}_c$ and $\tilde{K}_c=K_c$ on $\mathcal{U}_c,$ it is the linearized flow of $K_c$ along   $x|_{[T_0-T,T_0]},$ with $\psi_0(0)={\rm Id}.$ Since  $\psi(T_0-T+t)V  = \psi_0(t) \circ \psi(T_0-T)V  ,$
we have
\begin{align*}
\mu_{\rm RS}( \psi V|_{[T_0-T,T_0]} , V) &=\mu_{\rm RS}(  \psi_0  ( \psi(T_0-T)V ), V)\\
&= \mu_{\rm RS}( {\rm Gr}(\psi_0), ( \psi(T_0-T)V ) \times V),
\end{align*}
where  the second identity follows from  \cite[Theorem 3.2]{RS93}. Here,    $\Delta$   denotes the diagonal in  $(\R^4 \oplus \R^4, -\omega \oplus \omega).$
It follows that 
\[
\mu_{\rm RS}( \psi V|_{[T_0-T,T_0]} , V)  \geq \mu_{\rm RS}({\rm Gr}(\psi_0), \Delta)-6.
\]
 Indeed,  let $L$ be a common Lagrangian complement to $\Delta, \psi(T_0-T)V \times V$ and ${\rm Gr}(\psi_0(T)).$ In the symplectic coordinates corresponding to the splitting $\R^4 \oplus \R^4 = \Delta \oplus L,$ we have
  $\Delta = \R^4 \times \{0\},$ $ \psi(T_0-T)V \times V={\rm Gr}(P)$ and  ${\rm Gr}(\psi_0(T)) = {\rm Gr}(Q)$ for some symmetric matrices $P$ and $Q.$ Then by \cite[Theorem 3.5]{RS93}, we have
\begin{align*}
&\mu_{\rm RS}( {\rm Gr}(\psi_0), ( \psi(T_0-T)V ) \times V) - \mu_{\rm RS}({\rm Gr}(\psi_0), \Delta)\\
=& \;   \frac{1}{2} ( {\rm sign}(Q-P) +{\rm sign}(P) - {\rm sign}(Q) )\\
\geq&\; -6,
\end{align*}
where ${\rm sign}(S)$ denotes the signature of the symmetric matrix $S.$

  It remains to show that $\mu_{\rm RS}({\rm Gr}(\psi_0), \Delta)> N+6.$
To this end,  we specify the choice of $T$ and the coordinate neighborhood $\mathcal{U}_c.$ For each $c,$ let $ A(c)$ be the change of coordinates matrix, namely, $\hat{z} =A(c)( z -z_*(c)).$ Let
\[
\Xi(t)=A(c_2) \mathfrak{T}_{z_*(c_2)} \circ D\phi_{K_{c_2}}^t(z_*(c_2)) \circ \mathfrak{T}_{z_*(c_2)}^{-1} A(c_2)^{-1}, \quad 0 \leq t \leq T
\]
denote the   linearized Hamiltonian flow at $z_*(c_2)$, expressed in the local coordinates at $c=c_2.$  Applying the computation in \cite[equation (7.7)]{LS25finite} to the quadratic part $Q,$ we have
\[
\mu_{\rm RS}( {\rm Gr} ( \Xi ), \Delta) \geq 2 \bigg\lfloor \frac{\hat{\lambda}_2(c_2) T}{2\pi}  \bigg\rfloor -1.
\]
Choose  $T>0$ such that
\[
2 \bigg\lfloor \frac{\hat{\lambda}_2(c_2) T}{2\pi}  \bigg\rfloor -1> N+6, \qquad \hat{\lambda}_2(c_2)T \notin 2 \pi \Z.
\]
The eigenvalues of $\Xi(T)$  are $e^{\pm \hat{\lambda}_1(c_2) T}$ and $e^{\pm i\hat{\lambda}_2 (c_2)T}.$   The first two eigenvalues are not equal to $1$ since $\hat{\lambda}_1(c_2)>0$ and $T>0.$ The choice of $T$ ensures that the other two eigenvalues are not equal to $1$ as well.  Therefore, $\Xi(T)-{\rm Id}$ is invertible, or equivalently, ${\rm Gr}(\Xi(T)) \cap \Delta = \{ 0 \}.$ 
In the local coordinates at   $c\neq c_2,$ the linearized flow along $x|_{[T_0-T,T_0]}$ is given by $\eta(t)=A(c) \psi_0(t) A(c)^{-1}, \; 0 \leq t \leq T.$ By smoothness of $K_c$ and the coordinate transformations and by continuous dependence of solutions on the coefficient matrices in the linearized equations,  after shrinking $\mathcal{U}_c$ and taking smaller $\varepsilon_N>0$ if necessary,  we may assume that $\eta(t)$ is uniformly close to $\Xi(t)$ on $[0,T].$ Therefore, the two linearized flows can be joined by   a homotopy through  symplectic paths starting at $I.$   Moreover, this homotopy can be chosen so that the time-$T$ matrices have no eigenvalue $1$ during the homotopy, implying that  their graphs at $t=T$ remain transverse to $\Delta$ during the homotopy. 
Homotopy invariance and naturality of the Robbin--Salamon index then give
\[
\mu_{\rm RS}({\rm Gr}(\psi_0), \Delta) = \mu_{\rm RS}({\rm Gr}(\eta), \Delta) = \mu_{\rm RS}({\rm Gr}(\Xi), \Delta) > N+6.
\]
This  completes the proof.
\end{proof}

We now construct a perturbed   action functional $\mathcal{A}_H^F \colon \mathscr{P}(M, L_{\rm col}^E) \times \R \to \R$ and relate the above-obtained index estimate to the grading. 
 Choose $H  = \tilde{K}_c $ near $\tilde{\Sigma}_c$ and extend it by constants. Let $F$ be the auxiliary compactly   supported Hamiltonian, constructed at the beginning of Section \ref{sec:proofsection}, so that  
\[
(\phi_F^1)^{-1}( L_{\rm col}^E) \cap \tilde{\Sigma}_c = L_{\rho_a}^E \cap \tilde{\Sigma}_c.
\]
We take the cutoff function $\ell$ to be independent of $c$ for $c>c_2$ sufficiently close to $c_2.$ We also choose $\beta  $ and  $\kappa$ as before.

We first assume that $\mathcal{A}_H^F$ is Morse. For   a critical point $z=(y,\tau),$ set $x_z(t):=\phi_{\tilde{K}_c}^{\tau t}(y(0)), \; t \in [0,1].$ Note that for $\tau>0,$ $x_z$ is  a Reeb chord from $L_{\rm col}^E$ to $L_{\rho_a}^E$ up to a reparametrization of time. The first half of $y$ is a time-reparametrization of $x_z,$ so their full Robbin--Salamon indices agree. 
The linearized flow of $x_z$ in the trivialization $\mathfrak{T}$ is given by
\[
\psi_z(t):=\mathfrak{T}_{x_z(t)} \circ D\phi_{\tilde{K}_c}^{\tau t}( y(0)) \circ \mathfrak{T}_{x_z(0)}^{-1}, \quad 0 \leq t \leq 1.
\]
Denote by $\hat{\psi}_z(t)$ the full   linearized flow along $y$, with $\tau$ fixed, expressed in $\mathfrak{T}.$ Let $R(t),$ $0\leq t \leq 1,$ denote the linearized flow of $\tilde{F} = (\pi/4)(x_1^2 + y_1^2)$ also in $\mathfrak{T}.$ By additivity, we have $\mu_{\rm RS}(\hat{\psi}_z V, V) = \mu_{\rm RS}(\psi_z V, V) + \mu_{\rm RS} ( R(t) \psi_z(1)V, V).$ 
An application of \cite[Theorems 3.2 and 3.5]{RS93} as above tells us that
\begin{align*}
\lvert \mu_{\rm RS}(R(t)\psi_z(1)V, V) \rvert = \lvert \mu_{\rm RS}({\rm Gr}(R), \psi_z(1)V \times V) \rvert  \leq \lvert \mu_{\rm RS}({\rm Gr}(R), \Delta) \rvert +6.
\end{align*}
It is not hard to compute 
\begin{equation}\label{eq:gradelta}
\mu_{\rm RS}({\rm Gr}(R), \Delta)=1.
\end{equation}
Consequently, together with the preceding estimate and   \eqref{eq:fadslfsih}, this provides 
\begin{equation}\label{eq:boundeforgraindin}
\lvert \lvert z \rvert - \mu_{\rm RS}( {\psi}_zV, V)  \rvert \leq 8.
\end{equation}

We now complete the proof in the Morse case. To obtain Reeb chords from the generators detected by RFH, we first observe that  their Lagrange multipliers are positive. Recall from the proof of Proposition \ref{prop:iundexesti} that  $\mu_{\rm RS}(\psi_z V, V) \leq 0$ if  $\tau \leq 0.$
This, together with \eqref{eq:boundeforgraindin}, implies that $\lvert z \rvert \leq 8$ if $\tau \leq 0.$ It follows that every generator of degree bigger than $8$ has positive multiplier.

Fix $N \in \N$ and choose integers $k_1 < \cdots <k_N$ with $k_1>8.$ We also choose  $\Theta>k_N+9. $  We then apply Proposition \ref{prop:iundexesti} with $\Theta $ in place of $N,$ to obtain $\varepsilon_{\Theta}>0$ and an open set $\mathcal{B}_{\Theta}$ for which the assertion of that proposition holds. 
At any point where the trajectory crosses $x_1 = \pm x_*(c), $ we have $\tilde{K}_c = K_c.$ Since  $B_*(c)=V_{c,1}(\pm x_*(c))+ {c}/{4},$
see Remark \ref{remark:valuev1xstar}, we have 
\[
p_1^2 + p_2^2 + 2(V_{c,2}(x_2) - \frac{c}{4}) = - 2B_*(c) \to 0 \quad \text{ as } \quad c \to c_2^+,
\]
where $p_1 = y_1 + A_2(x_2) $ and   $p_2 = y_2 + A_1(x_1).$ Since all terms on the left-hand side are nonnegative, we have
\[
p_1 \to 0, \qquad p_2 \to 0, \qquad V_{c,2}(x_2) - \frac{c}{4} \to 0 .
\]
It follows that as $c \to c_2^+,$ the points in $\tilde{\Sigma}_c \cap \{ x_1 = \pm x_*(c)\}$ approach the four equilibria.  Consequently, for $\varepsilon_N \in (0, \varepsilon_{\Theta})$ sufficiently small, we have
\begin{equation}\label{eq:crossinghappen}
\tilde{\Sigma}_c \cap \{ x_1 = \pm x_*(c) \} \subset \mathcal{B}_{\Theta} \quad \text{ for } \; c \in (c_2, c_2+\varepsilon_N).
\end{equation}
By \eqref{eq:RFHsuper}, the RFH complex has a generator in each chosen degree $k_j,$ $j=1,2, \ldots, N ,$ and each such generator has positive multiplier. Assume by contradiction that its $\tilde{K}_c$-trajectory segment   enters $\mathcal{B}_{\Theta},$ so that, by Proposition \ref{prop:iundexesti} and \eqref{eq:boundeforgraindin} we obtain
\[
\lvert  z\rvert \geq \mu_{\rm RS}(\psi_z V, V) - 8 > \Theta - 8 > k_N,
\]
a contradiction. Therefore, these Hamiltonian chords remain in the unmodified region. Reflecting these chords using $\rho_a$ gives at least $N$ distinct symmetric consecutive collision orbits unless a symmetric periodic collision orbit exists. See also \cite[Corollary 4.2]{kang2025consecutive}. This proves the assertion in the Morse case.

\medskip

We next consider the degenerate case. As in the subcritical case, the spectral invariant argument  yields infinitely many chords on $\tilde{\Sigma}_c.$ 
To obtain $N$ chords of the original Hamiltonian $K_c,$ we use Proposition  \ref{prop:iundexesti} and \eqref{eq:crossinghappen} to show that the Hamiltonian segments of generators in the chosen degrees avoid the region where the Hamiltonian was modified. We must also show that these generators have positive-time limits in pairwise distinct $\sigma$-orbits.

Fix $N \in \N$ and choose integers $k_1 < \cdots < k_N$ with $k_1 > 10$ and $k_{j+1} -k_j >4.$  Choose an integer $\Theta> k_N+9,$ and let $B_\Theta$ and $\varepsilon_\Theta$ be as in  Proposition  \ref{prop:iundexesti}. Choose $0<\varepsilon_N <\varepsilon_\Theta$ so that \eqref{eq:crossinghappen}  holds, and fix $c \in (c_2, c_2+\varepsilon_N).$

We now choose smooth time-dependent $\sigma$-invariant Hamiltonians $F_n \to F$ in $C^{\infty},$ with common compact support, such that $\mathcal{A}_H^{F_n}$ is Morse. Indeed, since  $\sigma$ acts freely along the trajectories of $F$ starting on $\tilde{\Sigma}_c,$ the local perturbations used in the genericity argument of \cite[Theorem 2.28]{merry2014a} can be extended to $\sigma$-invariant perturbations. Therefore, we can choose $F_n$ so that  each  $\mathcal{A}_H^{F_n}$ is Morse.

We first study how much the Robbin--Salamon index can change under a small perturbation of a Lagrangian path. Let $W$ be a fixed Lagrangian subspace of a symplectic vector space, and let $\Lambda$ be a Lagrangian path on $[0,1].$ For every sufficiently $C^0$-close path $\Lambda'$ with $\Lambda'(0)=\Lambda(0),$ we then have
\begin{equation}\label{eq:pathlagpertubindex}
    \lvert \mu_{\rm RS}(\Lambda',W) - \mu_{\rm RS}(\Lambda, W) \rvert \leq \tfrac{1}{2} \dim ( \Lambda(1)\cap W).
\end{equation}
To see this, we choose a small homotopy $\Lambda^s$ from $\Lambda$ to $\Lambda'$ fixing the common initial plane. Then by homotopy invariance and additivity we obtain 
\[
\mu_{\rm RS}(\Lambda',W)-\mu_{\rm RS}(\Lambda, W)=\mu_{\rm RS}( s \mapsto \Lambda^s(1), W).
\]
Since $\Lambda'$ is $C^0$-close to $\Lambda$, we can choose the homotopy $\Lambda^s,$ so that $\Lambda^s(1)$ stays sufficiently close to $\Lambda(1).$ Then by applying  \cite[Proposition 2.3.1(g)]{Zhu} gives 
 \begin{align*}
 0 &\leq {\rm Mas}_+\{ \Lambda^s(1), W\} \\
 &= \mu_{\rm RS}(s \mapsto \Lambda^s(1), W) - \tfrac{1}{2}\left( \dim ( \Lambda(1) \cap W) + \dim ( \Lambda'(1) \cap W) \right) \\
 &\leq \dim ( \Lambda(1) \cap W) - \dim ( \Lambda'(1) \cap W) ,
 \end{align*}
 where the identity follows the localization property of the Robbin--Salamon index and \cite[Lemma 2.3.2]{Zhu}. This proves  \eqref{eq:pathlagpertubindex}.

We now apply \eqref{eq:pathlagpertubindex} to the perturbed linearized flow.   For $q \in \tilde{\Sigma}_c,$    denote by $R_{n,q}(t), \; 0 \leq t \leq 1,$ the linearized flow along $t \mapsto \phi_{F_n}^t(q),$  expressed in $\mathfrak{T}.$ Then   $R_{n,q}(t) \to R(t)$ uniformly for $ q \in \tilde{\Sigma}_c$ and $ t\in [0,1].$  Since $R_{n,q}(0)=R(0)=I_4,$  for each $q,$ both paths $t \mapsto {\rm Gr}(R_{n,q}(t))$ and $t \mapsto {\rm Gr}(R(t))$ start at $\Delta.$ By the formula for $\phi_{\tilde{F}}^1$ at the beginning of this section, the map $R(1)$ rotates the $(x_1, y_1)$-plane by $-\pi/2$ and acts as the identity on the $(x_2, y_2)$-plane. Hence, $\dim \left( {\rm Gr}( R(1))\cap \Delta\right) = \dim \ker ( R(1) - I_4) = 2.$ Then applying \eqref{eq:pathlagpertubindex} with $W=\Delta$ and using  \eqref{eq:gradelta} we obtain 
\[
\lvert \mu_{\rm RS}({\rm Gr}(R_{n,q}), \Delta) \rvert \leq 2 ,
\]
 uniformly in $q \in \tilde{\Sigma}_c$ and for all sufficiently large $n.$

For each sufficiently large $n,$  let $z = (y,\tau)$ be a   critical point of $\mathcal{A}_H^{F_n}.$  Define $x_z,$ $\psi_z$ and $\hat{\psi}_z$ as above.  On the second half of the path, the chain rule gives
\[
\hat{\psi}_z(t) = R_{n,q}( \kappa(t))\psi_z(1),\quad \frac{1}{2}\leq t \leq 1, \quad q = y ( \tfrac{1}{2}).
\]
By \eqref{eq:fadslfsih}, additivity and the same   index comparison used to derive \eqref{eq:boundeforgraindin},  we have
\begin{align*}
    \lvert \lvert z \rvert - \mu_{\rm RS}(\psi_z V, V) \rvert & \leq  \lvert \lvert  z \rvert  - \mu_{\rm RS}(\hat{\psi}_z V, V) \rvert  + \lvert \mu_{\rm RS}(\hat{\psi}_z V, V)- \mu_{\rm RS}( {\psi}_z V, V)\rvert\\
    & \leq 1 + \lvert \mu_{\rm RS}( R_{n,q}(t)\psi_z (1)V, V) \rvert \\
    & = 1 + \lvert \mu_{\rm RS}( {\rm Gr}(R_{n,q}), \psi_z(1)V \times V) \rvert \\
    & \leq 1 + \lvert \mu_{\rm RS}({\rm Gr}(R_{n,q}), \Delta) \rvert + 6 \\
    & \leq 9.
\end{align*}
As in the Morse case, $\tau \leq0$ implies $\mu_{\rm RS}(\psi_z V ,V) \leq0.$ Moreover, as $\tau \to 0,$ the paths $\psi_z V$ converge to the constant path $V,$ uniformly for initial points on the Legendrian circle $\tilde{\Sigma}_c \cap L_{\rm col}^E.$ Note that the constant path has index zero. Applying \eqref{eq:pathlagpertubindex} to $\Lambda(t) \equiv V,$ $\Lambda'(t) = \psi_z(t)V$ and $W = V$ gives
\[
\lvert \mu_{\rm RS}(\psi_z V, V) \rvert \leq \frac{1}{2}\dim V =1
\]
for all sufficiently small $\lvert \tau \rvert.$
Since $k_1>10,$ every generator $z$ of degree $k_j$   has $\tau \geq \delta $ for some $\delta>0$ which is independent of $n.$ They also satisfy
\[
\mu_{\rm RS}(\psi_z V, V) \leq k_j + 9 \leq k_{N} + 9 < \Theta.
\]
Hence, $x_z([0,1]) \cap \mathcal{B}_{\Theta} = \emptyset$ by Proposition \ref{prop:iundexesti}.

By \eqref{eq:RFHsuper}, for each $j=1, \ldots, N,$ we   choose a nonzero equivariant RFH class $\mathfrak{a}_j$~of degree $k_j$ for $F_1$ in the component $[{\rm pt.}].$ For each $n \geq 1,$ consider  the homotopy $\mathfrak{H}_{n,s}:=(1-s)F_1 + s F_n, \;  0 \leq s \leq 1.$ This homotopy is   $\sigma$-invariant, and its support~is contained in a compact set, independent of $n$ and $s.$ Equivariant continuation along $\mathfrak{H}_{n,s}$ induces a degree-preserving isomorphism $\mathfrak{I}_n$ from the  RFH for $F_1$ to the RFH for $F_n$ as in Section \ref{sec:basiscsRFH}. Thus  $\mathfrak{a}_{j,n}:=\mathfrak{I}_{ n}(\mathfrak{a}_{j }) $
is a nonzero class of degree $k_j$ for $F_n.$

We now apply the spectral invariant argument from \cite[Section 4.3]{frauenfelder2019a} to these classes, as in \cite[Corollary 4.2]{kang2025consecutive}.
Since $\mathcal{A}_H^{F_n}$ is Morse, the spectral value of $\mathfrak{a}_{j,n}$ is realized by a critical point
\[
z_{j,n} = (y_{j,n}, \tau_{j,n}), \qquad \lvert z_{j,n} \rvert = k_j.
\]
The continuation estimate for $\mathfrak{H}_{n,s}$ gives 
\[
\mathcal{A}_H^{F_n}(z_{j,n}) \leq \mathcal{A}_H^{F_1}(z_{j,1}) +\int_0^1 \| F_n(t, \cdot) - F_1 (t, \cdot)\|_{\infty} dt \leq C
\]
for some constant $C,$ independent of $j$ and $n.$ Indeed, there are only finitely many classes $\mathfrak{a}_j,$ and $F_n  \to F$ in $C^\infty$ with common compact support. 

We use the above bound to obtain an upper bound on the multipliers. Write $z_{j,n}=(y, \tau)$ and $F_{n,t} = F_n(t, \cdot).$ Since $H(y(t))=0$ for $0 \leq t \leq \frac{1}{2}$ and $\lambda|_{ L_{\rm col}^E}=0,$  
\[
 \mathcal{A}_H^{F_n}(z_{j,n}) = \tau  \int_0^{1/2} \beta(t)\lambda(X_H)( y (t)) dt + \int_{1/2}^1 \dot{\kappa}(t) ( \lambda(X_{F_{n, \kappa(t)}}) - F_{n,\kappa(t)}) ( y (t))dt.
 \]
 Since $\beta \geq 0$ and has integral one on $[0,\frac{1}{2}],$ we have   
\[
\tau  \int_0^{1/2} \beta(t)\lambda(X_H)( y (t)) dt \geq m \tau ,
\]
where $m:= \min_{\tilde{\Sigma}_c} \lambda(X_H) >0.$ Indeed, $m$ is positive since $\tilde{\Sigma}_c$ is of contact type and compact. Similarly, since $\dot{\kappa}\geq 0$ and has integral one on $[\frac{1}{2},1],$ 
\[
 \int_{1/2}^1 \dot{\kappa}(t) ( \lambda(X_{F_{n, \kappa(t)}}) - F_{n,\kappa(t)}) ( y (t))dt\geq  - \| \lambda(X_{F_{n,\kappa(t)}})  - F_{n,\kappa(t)} \|_{\infty}.
\]
Since $F_n \to F,$ we have $M_0 := \sup_n \| \lambda(X_{F_{n,\kappa(t)}}) - F_{n,\kappa(t)} \|_{\infty} <\infty.$
Combining the above estimates, we obtain
\[
0 < \delta \leq \tau_{j,n} \leq \frac{C +M_0}{m} 
\]
for every $j$ and all sufficiently large $n.$
  The initial points $y_{j,n}(0)$ lie in the Legendrian circle $\tilde{\Sigma}_c \cap L_{\rm col}^E.$ Thus, after passing to a subsequence, we may assume that 
\[
z_{j,n} \to z_j = (y_j, \tau_j ) \in {\rm Crit}(\mathcal{A}_H^{F}), \quad \tau_j \geq \delta>0, \quad j = 1, \ldots, N.
\]

\begin{lemma}\label{lem:separatbiondegree}
Let $z_n, z_n' \in {\rm Crit}(\mathcal{A}_H^{F_n})$ be two sequences converging to limits in the same $\sigma$-orbit. Then   for $n$ sufficiently large, we have
\[
\lvert \lvert z_n \rvert - \lvert z_n' \rvert \rvert \leq 4.
\]

\end{lemma}
\begin{proof}
Applying $\sigma$ to one sequence, if necessary, preserves its grading and makes both limits equal to some $z.$ The full linearized paths $\hat{\psi}_{z_n}V$ and $\hat{\psi}_{z_n'}V$ start at $V$ and converge uniformly to $\hat{\psi}_z V.$ Applying \eqref{eq:pathlagpertubindex} to $z_n$  then gives
\[
\lvert \mu_{\rm RS}(\hat{\psi}_{z_n}V, V) - \mu_{\rm RS}(\hat{\psi}_{z }V, V) \rvert \leq \tfrac{1}{2} \dim ( \hat{\psi}_z (1)V  \cap V ) \leq 1.
\]
A similar inequality also holds for $z_n'.$  Thus the two full indices differ by at most two. Then this together with   \eqref{eq:fadslfsih} completes the proof.
\end{proof}

By Lemma \ref{lem:separatbiondegree} and the condition $k_{j+1}-k_j>4,$ the limits $z_1, \ldots, z_N$ lie in pairwise distinct $\sigma$-orbits. For each $j,$ the Hamiltonian $H$ segments of $z_{j,n}$ do not intersect $\mathcal{B}_\Theta$ for all sufficiently large $n.$ Since $\mathcal{B}_\Theta$ is open, the limit segment also does not intersect it. To reach the modified region, a segment starting at $x_1=0$ must cross either $\{x_1=x_*(c)\}$ or $\{ x_1 = -x_*(c)\}.$ However, by \eqref{eq:crossinghappen}, every such intersection lies in $\mathcal{B}_{\Theta}.$ Hence, the limit segments follow $K_c.$ They start on $L_{\rm col}^E$ and end on $L_{\rho_a}^E$ by \eqref{eq:compareKandiltdeK}. These $N$ Hamiltonian chords are distinct modulo $\sigma.$  
We now argue and obtain  the desired collision orbits as in the Morse case. The same argument applies to the Moon. This completes the proof.
\hfill $\square$



\appendix

\section{Proof of Lemma \ref{lem:transtionregion}}

By the choice of $r_*,$ \eqref{eq:atcritialenergy} holds on   $\overline{\mathcal{B}}_{r_*}(z_*(c_2)),$    and  the Taylor expansion in \eqref{eq:Taylorofkc}  and the estimates for the remainder term are valid  uniformly for $c \in I_*$ and $\lvert \zeta \rvert \leq r_*.$

  Assume by contradiction that the upper bound of $\lvert \zeta \rvert$ in the statement is false. Then we find a sequence $ \rho_n  \to 0 ,$      $ c_n  \in (c_2, c_2+\rho_n^3)$    and $ \zeta_n \in   \mathcal{S}_{c_n,\rho_n}$ satisfying $\lvert \zeta_n \rvert / \rho_n \to \infty.$
After passing to a subsequence,  we may assume that $\zeta_n \to \hat{\zeta}  $ with $\lvert \hat{\zeta}\rvert \leq r_*.$   Since $-2\rho_n \leq \alpha_{c_n}(\zeta_n) \leq -\rho_n,$ we have    $\alpha_{c_2}(\hat{\zeta})=0.$  Moreover,  $z_*(c_2) + \hat{\zeta} \in K_{c_2}^{-1}(0).$ Then by \eqref{eq:atcritialenergy} we have     $\hat{\zeta}=0,$ implying that $\zeta_n \to 0.$ We now normalize once more by setting   $\bar{\zeta}_n = \zeta_n / \lvert \zeta_n \rvert,$ so that $\lvert \bar{\zeta}_n\rvert =1$ for all $n.$  After passing to a subsequence, we may assume that $\bar{\zeta}_n \to \bar{\zeta}. $ Note that $\lvert \bar{\zeta}\rvert=1.$ Dividing  
\[
0=K_{c_n}( z_*(c_n)+\zeta_n) = K_{c_n}(z_*(c_n))+Q_{c_n}(\zeta_n)+\mathcal{R}_{c_n}(\zeta_n)
\]
by $\lvert \zeta_n\rvert^2$ we  obtain
 \[
 \frac{K_{c_n}( z_*(c_n))}{\lvert \zeta_n\rvert^2 }+Q_{c_n}(\bar{\zeta}_n)+ \frac{\mathcal{R}_{c_n}(\zeta_n)}{\lvert \zeta_n\rvert^2}=0.
 \]
Since $K_c(z_*(c))$ is  a smooth function of $c,$ (see the proof of Lemma \ref{lem:orderofvalues}) and $K_{c_2}(z_*(c_2))=0,$ we find  $\lvert K_{c_n}( z_*(c_n))\rvert \leq C\lvert c_n - c_2\rvert \leq C \rho_n^3$
for some $C>0, $ and thus
\[
\bigg\lvert \frac{K_{c_n}( z_*(c_n))}{\lvert \zeta_n\rvert^2 }\bigg\rvert \leq C \rho_n  \left( \frac{\rho_n  }{\lvert \zeta_n \rvert} \right)^2 \to 0.
\]
Moreover, since $\mathcal{R}_c$ is of order at least three uniformly for $c$ close to $c_2,$ we have
\[
\bigg\lvert \frac{\mathcal{R}_{c_n}(\zeta_n)}{\lvert \zeta_n\rvert^2} \bigg\rvert \leq D \lvert \zeta_n \rvert \to 0. 
\]
Consequently, we have $Q_{c_2}(\bar{\zeta})=0.$ 
Moreover,
\[
\bigg\lvert \frac{\alpha_{c_n}(\zeta_n)}{ \lvert \zeta_n\rvert}\bigg\rvert \leq \frac{2\rho_n}{\lvert \zeta_ n\rvert} \to 0,
\]
and hence  $\alpha_{c_2}(\bar{\zeta})=0.$
It follows  that $\bar{\zeta}=0,$ contradicting $\lvert \bar{\zeta}\rvert=1.$ Thus, there exist  $\rho_0>0$ sufficiently small and $M_0>0$ such that   $\lvert \zeta \rvert \leq M_0 \rho$  
for every $0<\rho<\rho_0 $ and every $c \in (c_2, c_2+\rho^3).$

For the lower bound,   after decreasing $\rho_0>0$ if necessary, we may assume that $0<\rho_0<\varepsilon_*^{1/3},$ so that $c_2 <c<c_2+\rho^3$ implies $ c \in I_*.$
Set $B:= \max_{c \in I_*} \lvert \mu_c \rvert < \infty.$ Therefore, we find
\[
\lvert \alpha_c(\zeta)\rvert  \leq \lvert (\mu_c-\frac{1}{4} )\xi_1 \rvert + \lvert \eta_2\rvert \leq  ( \lvert \mu_c-\frac{1}{4} \rvert +1)\lvert \zeta \rvert \leq (B+2)\lvert \zeta \rvert.
\]
Since $\lvert \alpha_c(\zeta)\rvert \geq \rho$ for all $\zeta \in \mathcal{S}_{c,\rho},$ we conclude that $\lvert \zeta \rvert \geq m_0\rho$, where $ m_0:=\frac{1}{B+2}.$

To finish the proof, it remains to establish the estimate for $\xi_1.$ Assume by contradiction that no such $\kappa$ exists. Then 
there are sequences $  \rho_n \to 0,$   $c_n \in (c_2, c_2+\rho_n^3)$   and   $ \zeta_n \in \mathcal{S}_{c_n, \rho_n}$ such that  
\begin{equation}\label{eq:appendix1last}
  \frac{\xi_{1,n}}{ \rho_n}  \geq - \frac{1}{n},
\end{equation}
where $\zeta_n=(\xi_{1,n}, \xi_{2,n}, \eta_{1,n}, \eta_{2,n})  .$ 
By the bounds already proved above, we have
\[
m_0 \rho_n \leq \lvert \zeta_n \rvert  \leq M_0 \rho_n.
\]
Set $\hat{\zeta}_n = \zeta_n / \rho_n ,$ so that $m_0 \leq \lvert \hat{\zeta}_n \rvert \leq M_0$ for all $n.$ After passing to a subsequence, we may assume that $\hat{\zeta}_n \to \hat{\zeta}=(\hat{\xi}_1, \hat{\xi} _2, \hat{\eta}_1, \hat{\eta}_2). $
 Note that $0< m_0 \leq \lvert \hat{\zeta} \rvert \leq  M_0.$
 Arguing as above, we find $Q_{c_2}(\hat{\zeta})=0.$ On the other hand,
 \[
 -2\rho_n \leq \alpha_{c_n} (\zeta_n) \leq -\rho_n \;\;\Rightarrow \;\; -2 \leq \alpha_{c_n}( \hat{\zeta}_n) \leq -1.
 \]
and hence, passing to the limit,  we obtain $-2 \leq \alpha_{c_2}(\hat{\zeta})\leq-1.$ In particular
 \begin{equation}\label{eq:alphsc2hat}
 \alpha_{c_2}(\hat{\zeta})<0.
 \end{equation}
 By Remark \ref{rmk:separation},  
\[
{\rm sign}( \alpha_{c_2}(\hat{\zeta}) ) = {\rm sign} (\hat{\xi}_1).
\]
This, together with \eqref{eq:alphsc2hat}, shows that $\hat{\xi}_1<0.$ However,  assumption   \eqref{eq:appendix1last} gives a contradiction
 \[
 \hat{\xi}_1 =   \lim_{n \to \infty} \frac{\xi_{1,n}}{\rho_n} \geq 0.
 \]
Therefore, after decreasing $\rho_0>0$ if necessary,  we find $\kappa>0$ such that   $-\xi_1 / \rho \geq \kappa$ for every $0<\rho<\rho_0$ and every $c \in (c_2, c_2+\rho^3).$  This finishes the proof.


\section{Proof of Lemma \ref{lem:cappingfunction}}

 Let $\varrho$ and $c$ be as in the statement.   Since  $\varrho_0^3 <  \varepsilon_*,$ we have $c \in I_*.$  
 Since $\lambda_1(c) = V_{c,1}''(x_*(c))<0$ for every $c \in I_*,$ we choose $a_0>0,$ independent of $\rho,$ such that  $\lambda_1(c) \leq -2a_0$ for every $c \in I_*. $ 
 Then since $x_*(c)$ is a nondegenerate local maximum of $V_{c,1} $ and depends smoothly on $c,$  taking   $\varrho_0>0$  smaller if necessary, we have
 \[
 V''_{c,1}(x_*(c)+s) \leq -a_0, \qquad c \in I_*, \; \; 0 \leq s \leq 4 \varrho_0 .
 \]
Set $C_0:=\max \{    \lvert V_{c,1}''(x_*(c)+s) \rvert : c \in I_*, \; 0 \leq s \leq 4        \}<\infty.$ Since $\varrho_0<1$ and $V_{c,1}'( x_*(c))=0,$ 
\begin{equation*}\label{eq:a0exist}
a_0 s \leq - V_{c,1}'( x_*(c)+s) = - \int_0^s V_{c,1}''(x_*(c) + u)du \leq C_0s, 
\end{equation*}
for every $c \in I_*$ and $0 \leq s \leq 4 \varrho_0.$

We now take $\mathfrak{c}_0 \in (0, a_0/4),$ independent of $c$ and $\rho,$  and construct the desired   smooth function $\psi$   as follows. We first define  $\psi \equiv0$ on $(-\infty, 2\varrho].$ Then for   $s \in [\varrho, 2\varrho],$ we  find  
\begin{align*}
(s-3\varrho) (V_{c,1}'(x_*(c)+s)+\psi'(s) ) &=  (s-3\varrho)V_{c,1}'(x_*(c)+s)\\
&\geq a_0 s(3\varrho -s) \\
&\geq \frac{a_0}{2} (3\varrho-s)^2 \\
&> \mathfrak{c}_0(s-3\varrho)^2.
\end{align*}

We now define $\psi$ on $(2\varrho, 4\varrho).$  Fix $\tau >0   $ sufficiently small, so that
\[
\tau < \min \{ \frac{1}{2}, \frac{9a_0}{4 ( 3 C_0 + 4 \mathfrak{c}_0)  } \} .
\]
Pick any smooth cut-off function $\ell \colon [2\varrho, 3\varrho] \to [0,1]$ such that $\ell \equiv0$ on $[2\varrho, (3-2\tau) \varrho]$ and $\ell \equiv1$ on $[(3-\tau)\varrho, 3\varrho].$ We then set
\[
\psi'(s) = \ell(s)\left( -V_{c,1}'(x_*(c)+s)+2\mathfrak{c}_0(s-3\varrho)\right), \quad 2 \varrho \leq s \leq 3 \varrho
\]
so that $\psi'\equiv0$ on $[2\varrho, (3-2\tau) \varrho]$  and $\psi'(s) =   -V_{c,1}'(x_*(c)+s)+2\mathfrak{c}_0(s-3\varrho)$ for every $(3-\tau) \varrho \leq s \leq 3 \varrho.$
In particular, $\psi$ is smooth on $[2\varrho, 3\varrho].$ We   find
\begin{align*}
   V_{c,1}'(x_*(c)+s)+\psi'(s)     =   (1-\ell(s))V_{c,1}'(x_*(c)+s)+ 2\ell(s) \mathfrak{c}_0(s-3\varrho).
    \end{align*}
Since $s-3\varrho\leq 0$ and  $V_{c,1}'(x_*(c)+s) \leq -a_0 s < -\mathfrak{c}_0s \leq \mathfrak{c}_0(s-3\varrho) ,$ we  have
\begin{align*} 
 V_{c,1}'( x_*(c)+s )+\psi'(s) \leq (1-\ell(s)) \mathfrak{c}_0 (s-3\varrho) + \ell(s) \mathfrak{c}_0 (s-3\varrho) =\mathfrak{c}_0(s-3\varrho)
\end{align*}
from which the second property of $\psi$ holds on $[2\varrho, 3\varrho].$ Moreover, we observe that   $\psi' \geq 0$ on $[2\varrho, 3\varrho].$ Indeed, for $s \in [2 \varrho, 3 \varrho],$ we have 
\[
-V_{c,1}'(x_*(c)+s) +2\mathfrak{c}_0 (s-3\varrho)  \geq a_0 s -2\mathfrak{c}_0 \varrho \geq a_0 s - \frac{1}{2}a_0 s = \frac{1}{2}a_0 s>0.  
\]
It remains to define $\psi$ on $[3\varrho, 4\varrho).$ To this end, we pick any smooth function $h \colon [3\varrho, 4\varrho) \to [0,+\infty)$ such that $h \equiv0$ near $3\varrho$ and 
\begin{equation}\label{eq:hprop}
\int_{3\varrho}^s h(u)du \to +\infty \quad \text{ as } \; s \to 4\varrho^- 
\end{equation}
and then define
\[
\psi'(s)=-V_{c,1}'(x_*(c)+s)+2\mathfrak{c}_0 (s-3\varrho) +h(s), \quad s \in [3\varrho, 4\varrho). 
\]
Note that $\psi$ is now a non-decreasing smooth function on $(-\infty,4\varrho).$ Since $\psi \equiv 0$ on $(-\infty, 2\varrho],$  it is non-negative.  Moreover, because of \eqref{eq:hprop}, the third property of $\psi$ is established. Also, we find
\begin{align*}
    V_{c,1}'(x_*(c)+s) +\psi'(s) = 2\mathfrak{c}_0(s-3\varrho) + h(s) \geq \mathfrak{c}_0 (s-3\varrho)
\end{align*}
from which the second property of $\psi$ is established.

Recall from \eqref{eq:positiefVatx} that
\[
-\frac{c}{4} - V_{c,1}( x_*(c)+3\varrho)>0, \qquad c_2 < c \leq c_2+\varrho^3.
\]
To obtain the last property of $\psi$, we now take $\psi(3\varrho)$ sufficiently small as follows.  By construction we have
\[
\psi(3 \varrho)= \int_{2\varrho}^{3\varrho}\psi'(s)ds = \int_{(3-2\tau)\varrho}^{3\varrho} \psi'(s)ds.
\]
For $s \in [(3-2\tau)\varrho, 3\varrho]$  the definition of $C_0$ gives
\[
\lvert V_{c,1}'(x_*(c)+s)\rvert  \leq C_0 s  \leq 3 C_0 \varrho,
\]
and hence
\[
\lvert \psi'(s)\rvert \leq \lvert V_{c,1}'(x_*(c)+s) \rvert + \lvert 2\mathfrak{c}_0(s-3\varrho) \rvert  \leq C_1 \varrho, 
\]
where $C_1:=3 C_0  + 4\mathfrak{c}_0>0$ is independent of $\rho, c$ and $\varrho.$ Thus we obtain
\[
0 < \psi(3\varrho) \leq \int_{(3-2\tau)\varrho}^{3\varrho} C_1 \varrho ds = 2 C_1 \tau \varrho^2.
\]
On the other hand, recall from Remark \ref{remark:valuev1xstar} that $B_*(c)=V_{c,1}(x_*(c))+\frac{c}{4}$ satisfies $B_*(c_2)=0$ and $B_*(c)<0$ for $c \in (c_2, c_2+\varepsilon_1].$ Hence,
\begin{align*}
    -\frac{c}{4} - V_{c,1}(x_*(c)+3\varrho) &= - B_*(c) + V_{c,1}(x_*(c)) - V_{c,1}(x_*(c)+3\varrho)\\
    &\geq V_{c,1}(x_*(c)) - V_{c,1}(x_*(c)+3\varrho)\\
    &= -\int_0^{3\varrho} V_{c,1}'(x_*(c)+s) ds\\
    &\geq \int_0^{3\varrho} a_0s ds \\
    &= \frac{9}{2}a_0 \varrho^2.
\end{align*}
Thus, by the choice of $\tau,$ we have
\[
0 < \psi(3\varrho) \leq 2C_1 \tau \varrho^2 < \frac{9}{2}a_0 \varrho^2  \leq  - \frac{c}{4} - V_{c,1}( x_*(c)+3\varrho).
\]
Thus,
\[
d_{c,\varrho}:= - \frac{c}{4} - V_{c,1} (x_*(c) +3\varrho)-\psi (3\varrho)>0.
\]
For the upper bound, set $C_B:=\max_{c \in I_*}\lvert B_*'(c)\rvert,$ which is independent of $\rho,c$ and $\varrho.$ Since  $B_*$ is smooth,   the mean value theorem yields  
\[
0 < -B_*(c) \leq C_B ( c-c_2) \leq C_B \varrho^3.
\]
Since $-V_{c,1}'(x_*(c)+s) \leq C_0 s,$  we have 
\[
0< V_{c,1}(x_*(c)) - V_{c,1}(x_*(c)+3\varrho) \leq \int_0^{3\varrho} C_0 s ds =  \frac{9}{2}C_0 \varrho^2.
\]
Therefore, we find, after shrinking $\varrho_0$ so that $\varrho_0<1,$ that 
\begin{align*}
d_{c,\varrho}  &<- \frac{c}{4} - V_{c,1}(x_*(c)+3\varrho)  \\
&= -B_*(c) + V_{c,1}(x_*(c)) - V_{c,1}(x_*(c)+3\varrho)\\
&\leq \mathfrak{c}_1 \varrho^2,  
\end{align*}
where the constant $\mathfrak{c}_1,$ independent of $\rho,c$ and $\varrho,$ is defined by $\mathfrak{c}_1:= C_B + 9C_0/2>0.$ 
The last property thus follows. The desired properties of the derivatives follow from property ii) and the construction.   This completes the proof.


\section{Proof of Lemma \ref{lem:transtionregion2}}

Since $c_2<c \leq c_2+\varrho^3 <c_2 +\varrho_0^3 <c_2+\varepsilon_*, $ we have $c \in I_*.$ 
 Recall from \eqref{eq:Taylorofkc} and \eqref{eq:qudarictopart} that by Taylor's theorem  we have
\[  
K_c( z_*+\zeta) =  B_*(c)+  Q_c(\zeta)+\mathcal{R}_c(\zeta) ,
\]
where  the quadratic part $Q_c$ of $K_c$ at $z_*(c) $ is given as   
\begin{align}
  \nonumber  2Q_c(\zeta)  &=    ( \eta_1 + \frac{1}{4}\xi_2)^2 + ( \eta_2 +  \frac{1}{4}\delta_c\xi_1)^2  +   \lambda_1(c) \xi_1^2  + \lambda_2(c) \xi_2^2 \\
  \label{eq:bcquadratic}   &= ( \eta_1 + \frac{1}{4}\xi_2)^2 + \mathfrak{b}_c(\zeta)^2    +   \lambda_1(c) \xi_1^2  + \lambda_2(c) \xi_2^2 .  
\end{align}
Since $\tilde{K}_c=K_c=0$ along $\tilde{\mathfrak{S}}_{c,\varrho, \mathfrak{n}},$ using the Taylor expansion we find that
\[
\lvert Q_c(\zeta)\rvert \leq \lvert B_*(c) \rvert + C_{\rm Tay} \lvert \zeta \rvert^3, \quad \zeta \in \tilde{\mathfrak{S}}_{c,\varrho, \mathfrak{n}}.
\]
Moreover,   using the mean value theorem,   $\lvert B_*(c)\rvert \leq \mathcal{C}(c-c_2) \leq \mathcal{C} \varrho^3 $ for some  $\mathcal{C}>0,$ independent of $\mathfrak{n},c$ and $\varrho.$ Consequently, by \eqref{eq:boundeforzetaxc} we have
\[
 \lvert \zeta \rvert^2 \leq \mathcal{C}_Q ( \mathcal{C} \varrho^3 + C_{\rm Tay} \lvert \zeta \rvert^3 +\xi_1^2), \quad \zeta \in \tilde{\mathfrak{S}}_{c,\varrho, \mathfrak{n}}.
\]
Moreover, by the choice of $r_*>0$ as in \eqref{eq:SDFK} we have 
\[
\mathcal{C}_Q C_{\rm Tay} \lvert \zeta \rvert^3 \leq \mathcal{C}_Q C_{\rm Tay} r_* \lvert \zeta \rvert^2 \leq \frac{1}{2}\lvert \zeta \rvert^2, 
\]
and thus
\[
\lvert \zeta \rvert^2 \leq 2 \mathcal{C}_Q ( \mathcal{C} \varrho^3 + \xi_1^2), \quad \zeta \in \tilde{\mathfrak{S}}_{c,\varrho, \mathfrak{n}}.
\]

We now claim that $\xi_1>0$ for  $\zeta \in \tilde{\mathfrak{S}}_{c,\varrho, \mathfrak{n}}. $ Assume by contradiction that $\xi_1 \leq 0. $ Recall from the proof of Lemma \ref{lem:newlame1} that $\alpha_c(\zeta)=\mathfrak{a}_c \xi_1 -\mathfrak{b}_c(\zeta).$
Since   $\mathfrak{a}_c - \sqrt{-\lambda_1(c)} \geq 2 d_0$ for $c_2<c<c_2+\rho_1^3,$ see \eqref{eq:cdsd02}, 
using \eqref{eq:bcquadratic} we have
\[
\lvert \mathfrak{b}_c(\zeta) \rvert \leq \sqrt{-\lambda_1(c)}\lvert \xi_1 \rvert + \sqrt{ 2\lvert Q_c(\zeta)\rvert}  ,
\]
and hence
\[
\alpha_c(\zeta)  \leq -2 d_0 \lvert \xi_1 \rvert +\sqrt{2\lvert Q_c(\zeta)\rvert}.
\]
On the other hand,   we have
\[
\lvert Q_c(\zeta)\rvert \leq \lvert B_*(c) \rvert + \lvert \mathcal{R}_c(\zeta)\rvert \leq \mathcal{C}\varrho^3 + C_{\rm Tay}\lvert \zeta \rvert^3  ,
\]
and hence 
\begin{align*}
\sqrt{ \lvert Q_c(\zeta)\rvert} &\leq \sqrt{\mathcal{C}}\varrho^{3/2}+\sqrt{C_{\rm Tay}}\lvert \zeta \rvert^{3/2}\\
&\leq \sqrt{\mathcal{C}}\varrho^{3/2}+\sqrt{C_{\rm Tay}}( 2 \mathcal{C}_Q ( \mathcal{C} \varrho^3 + \xi_1^2))^{3/4} \\
&\leq \sqrt{\mathcal{C}}\varrho^{3/2}+2 \sqrt{C_{\rm Tay}}  \mathcal{C}_Q^{3/4} ( \mathcal{C}^{3/4} \varrho^{9/4} + \lvert \xi_1 \rvert^{3/2}) .  
\end{align*}
Again by  \eqref{eq:SDFK}, we have
\[
2 \sqrt{2C_{\rm Tay}}  \mathcal{C}_Q^{3/4}   \lvert \xi_1 \rvert^{3/2} \leq 2 \sqrt{2C_{\rm Tay}}  \mathcal{C}_Q^{3/4}  \sqrt{r_*} \lvert \xi_1 \rvert  \leq  d_0 \lvert \xi_1 \rvert.
\]
Therefore,
\[
\sqrt{2  \lvert Q_c(\zeta)\rvert}  \leq  \sqrt{2}  ( \sqrt{\mathcal{C}}+2 \sqrt{C_{\rm Tay}}  \mathcal{C}_Q^{3/4}   \mathcal{C}^{3/4} )\varrho^{3/2} +d_0  \lvert \xi_1 \rvert   , 
\]
and hence
\[
\alpha_c(\zeta) \leq -  d_0 \lvert \xi_1 \rvert + \tilde{\mathcal{C}} \varrho^{3/2} \leq \tilde{\mathcal{C}}\varrho^{3/2} ,   
\]
where $\tilde{\mathcal{C}}:=\sqrt{2}  ( \sqrt{\mathcal{C}}+2 \sqrt{C_{\rm Tay}}  \mathcal{C}_Q^{3/4}   \mathcal{C}^{3/4} ).$
Given $\mathfrak{n},$ we  decrease $\varrho_0,$ so that $\tilde{\mathcal{C}}\sqrt{\varrho_0}\mathfrak{n}<1,$ so that $\alpha_c(\zeta) <  {\varrho}/{\mathfrak{n}},$
a contradiction to $\alpha_c \geq \varrho / \mathfrak{n}$ on $\tilde{\mathfrak{S}}_{c,\varrho,\mathfrak{n}}.$ This proves the claim.

By the above claim and the definition of $  \tilde{\mathfrak{S}}_{c,\varrho, \mathfrak{n}},$ we have $0<\xi_1 < \varrho.$ Since  $\varrho_0<1,$ we have  $\varrho_0^3 \leq \varrho_0^2,$ and hence $
\mathcal{C} \varrho^3 +\xi_1^2 \leq \mathcal{C}\varrho^2 + \varrho^2 = (\mathcal{C}+1)\varrho^2.$
The desired upper bound for $\lvert \zeta \rvert $ such that $\zeta \in \tilde{\mathfrak{S}}_{c,\varrho,\mathfrak{n}} $ now holds with $\mathfrak{E}_0:=\sqrt{2\mathcal{C}_Q(\mathcal{C}+1)}.$

To get the stronger bound   for $\zeta \in{\mathfrak{S}}_{c,\varrho,\mathfrak{n}},$ we observe that by \eqref{eq:SDFK} we have 
\begin{equation*}\label{eq:SDFK2}
C_{\rm Tay}\lvert \zeta \rvert^3 \leq C_{\rm Tay} r_* \lvert \zeta \rvert^2 \leq \frac{m}{4}\lvert \zeta\rvert^2.
\end{equation*}
Consequently,
\begin{align*}
\frac{m}{2}\lvert \zeta\rvert^2 &\leq \lvert Q_c( \zeta)\rvert + \alpha_c(\zeta)^2 \\ &\leq \lvert B_*(c)\rvert +C_{\rm Tay}\lvert \zeta\rvert^3 +  \alpha_c(\zeta)^2 \\ &\leq \mathcal{C}\varrho^3   +\frac{m}{4}\lvert \zeta\rvert^2 + \frac{4\varrho^2}{\mathfrak{n}^2},
\end{align*}
and thus 
\[
\frac{m}{4}\lvert \zeta \rvert^2 \leq  ( \mathcal{C}\varrho  + \frac{4 }{\mathfrak{n}^2}) \varrho^2.
\]
Shrinking $\varrho_0>0$ if necessary, we may assume that $\mathcal{C}\mathfrak{n}^2\varrho_0 \leq 1 .$ Then we have
\[
 \lvert \zeta \rvert^2 \leq \frac{4}{m}( \mathcal{C}\varrho_0  + \frac{4 }{\mathfrak{n}^2}) \varrho^2 \leq \mathfrak{R}_0^2\frac{\varrho^2}{\mathfrak{n}^2}, \quad \; \mathfrak{R}_0:=\sqrt{\frac{20}{m  }}.
\]
 This provides the desired upper bound for $\zeta \in  {\mathfrak{S}}_{c,\varrho,\mathfrak{n}}.$

 It remains to obtain the lower    bound for $\xi_1.$ Since $\xi_1>0,$ using the previous computations we have
 \begin{align*}
 \alpha_c(\zeta) &\leq \mathfrak{a}_c \xi_1 +\lvert \mathfrak{b}_c (\zeta)\rvert \\
 &\leq ( \mathfrak{a}_c+ \sqrt{-\lambda_1(c)})\xi_1  + \sqrt{2 \lvert Q_c (\zeta)\rvert}\\
 & \leq ( \mathfrak{a}_c+ \sqrt{-\lambda_1(c)})\xi_1  + \sqrt{2( \mathcal{C}\varrho^3 + C_{\rm Tay}\lvert \zeta \rvert^3)} \\
  & \leq ( \mathfrak{a}_c+ \sqrt{-\lambda_1(c)})\xi_1  + \sqrt{2( \mathcal{C} + C_{\rm Tay} \mathfrak{E}_0^3  )}\varrho^{3/2}. 
 \end{align*}
 By decreasing $\varrho_0>0$ further if necessary, we assume that 
 \[
 \sqrt{2( \mathcal{C} + C_{\rm Tay} \mathfrak{E}_0^3  ) \varrho_0} \leq \frac{1}{2 \mathfrak{n}},
 \]
so that
\[
\alpha_c(\zeta) \leq ( \mathfrak{a}_c+ \sqrt{-\lambda_1(c)})\xi_1  + \frac{1}{2\mathfrak{n}}\varrho. 
\]
Put
\[
M:=\max \{\mathfrak{a}_c+\sqrt{-\lambda_1(c)} \mid c \in I_*\} \in (0,\infty).
\]
We then have
\[
\frac{\varrho}{\mathfrak{n}}\leq \alpha_c(\zeta) \leq (\mathfrak{a}_c+\sqrt{-\lambda_1(c)} )\xi_1 + \frac{\varrho}{2\mathfrak{n}}\leq M \xi_1 + \frac{\varrho}{2\mathfrak{n}},
\]
    implying that
    \[
    \xi_1 \geq \mathfrak{r}_0 \frac{\varrho}{\mathfrak{n}}, \qquad \mathfrak{r}_0:= \frac{1}{2M}>0.
    \]
    This completes the proof.


\section{Proof of Proposition \ref{prop:weinstein}}\label{appendproof}

We apply the following special case of the   cancellation theorem due to Cieliebak and Eliashberg \cite[Proposition 12.22]{CE12book}.

\begin{theorem}\label{thm:cancle}
  Let $(W, \omega, X, \phi)$ be a compact four-dimensional Weinstein cobordism. Suppose that the Morse function $\phi$   has exactly two critical points, a  minimum $q$ and a saddle $p$ of index one, with a unique trajectory of $X$ from $q$ to $p,$  up to time translation.   For any neighborhood $U$ of $ {\rm cl}_W(W^s_X(p)),$ 
  there is a Weinstein homotopy $(\omega, X_t, \phi_t)_{t \in [0,1]}$ with    $(X_0,\phi_0)=(X,\phi),$ having the following properties:
  \begin{enumerate}

\item  For every $t \in [0,1],$ we have   $(X_t,\phi_t)=(X,\phi)$ near $\partial W$ and on $W \setminus U.$

\item   For every $t \in [0,1],$ all critical points of $\phi_t$ lie in ${\rm cl}_W(W^s_X(p)).$

\item  At some $t_0 \in (0,1),$ the   two critical points degenerate into a single birth-death point. For $t>t_0,$ the function $\phi_t$ has no critical points. 

  \end{enumerate}
 In particular,  $d\phi_1(X_1)>0$ on $W.$

\end{theorem}

\begin{remark}
In the original theorem, there is a transversality condition on the stable manifold of $p$ and the unstable manifold of $q.$ In our case above, this is automatic. Indeed, since $q$ is a source of the forward flow of $X,$ $W_X^u(q)$ is open in $W$ and hence intersects $W^s_X(p)$ transversely along the connecting trajectory.
\end{remark}

We shall construct a Weinstein cobordism $W_c$ on which $\Lambda$ acts freely and show that the quotient $W_c/\Lambda$, again a Weinstein cobordism, contains precisely one minimum and one saddle of index one that are connected by a unique trajectory. We then apply the above theorem to this quotient.

Let $ \mathfrak{D}_c$ be the compact domain bounded by $\tilde{\Sigma}_c,$ where we have $\tilde{K}_c \leq 0.$ By Theorem \ref{thm:perturbation}, we have $Y(\tilde{K}_c)>0$ on $\tilde{\Sigma}_c = \partial \mathfrak{D}_c.$

\medskip

\noindent
{\bf Step 1.} We prove   $Y(\tilde{K}_c)>0$ on $\mathfrak{D}_c$ in the four coordinate neighborhoods, away from the critical points. 

\medskip

 Let $q_{c,\varrho}$ denote the minimum introduced by the capping in   Section \ref{sec:capping}. Its coordinates are $(3\varrho, 0,0,\Delta_{c,\varrho})$ as in Section \ref{sec:secondinterpo}. 
 By symmetry, it suffices to show 
 \[
 Y(\tilde{K}_c)>0 \quad \text{ on }  \left(\mathfrak{D}_c \cap  \mathcal{B}_{r_*}(z_*(c))\right)  \setminus \{   z_*(c) ,  q_{c,\varrho} \}.
 \]

 We first extend the bounds of Lemmas \ref{lem:transtionregion} and \ref{lem:transtionregion2} to the sublevel set. As before, we work with the local coordinates $\zeta = (\xi_1, \xi_2, \eta_1, \eta_2)$ centered at $z_*.$  Recall from the proof of Lemma \ref{lem:newlame1} that $\alpha_c = \mathfrak{a}_c \xi_1 - \mathfrak{b}_c(\zeta).$   
Let $Q_c$ be the quadratic part of the Taylor expansion of $K_c$ at $z_*(c),$ see \eqref{eq:qudarictopart}. 
It is easy to check that  $Q_c|_{ \{ \xi_1=0\}}$ is positive-definite. We now assume that $\alpha_c=0,$  so that $\mathfrak{a}_c \xi_1 = \mathfrak{b}_c.$   Then we have
\[
2Q_c = ( \eta_1 + \frac{1}{4}\xi_2)^2 + \left(   (\mu_c+\frac{\delta_c-1}{4})^2 + \lambda_1(c) \right) \xi_1^2 +\lambda_2(c)\xi_2^2.
\]
Inequality \eqref{eq:cdsd02} shows that this is also positive-definite. 
Therefore, there exist uniform constants $C_0, \kappa >0$ such that
\begin{equation}\label{eq:oposdououg}
Q_c + C_0 \alpha_c^2  \geq \kappa \lvert \zeta \rvert^2, \qquad Q_c + C_0 \xi_1^2 \geq \kappa \lvert \zeta \rvert^2.
\end{equation}
Since $\tilde{K}_c \leq 0 $ and the capping function $\psi = \psi_{c,\varrho}$ is nonnegative, by the Taylor expansion and  \eqref{eq:oposdououg}, after decreasing $r_*>0$ if necessary, we have
\begin{equation}\label{eq:compactbss}
  \lvert \zeta \rvert^2 \leq C_1 ( -B_*(c) + \alpha_c^2(\zeta) ), \qquad \lvert \zeta \rvert^2 \leq C_1 ( -B_*(c) + \xi_1^2)   
\end{equation}
on $\mathfrak{D}_c \cap \mathcal{B}_{r_*}(z_*(c)),$ where $C_1:=\max\{ 2/ \kappa, 2C_0/\kappa\}.$ 
Since $\mathfrak{b}_c^2+\lambda_1(c)\xi_1^2 \leq 2Q_c \leq -2B_*(c)-2\mathcal{R}_c,$ we have
\begin{align*}
\lvert \alpha_c - \mathfrak{a}_c \xi_1 \rvert   &= \lvert \mathfrak{b}_c\rvert  \\
&\leq \sqrt{-\lambda_1(c)}\lvert \xi_1 \rvert +  \sqrt{ - 2B_*(c)} + \sqrt{2C_{\rm Tay}}\lvert \zeta \rvert^{3/2}  \\
&\leq \sqrt{-\lambda_1(c)}\lvert \xi_1 \rvert +  \sqrt{ - 2B_*(c)} + \sqrt{2r_*C_{\rm Tay}}  \lvert \zeta \rvert  \\
&\leq  \sqrt{-\lambda_1(c)}\lvert \xi_1 \rvert +  \sqrt{ - 2B_*(c)} + \sqrt{2r_*C_1C_{\rm Tay}}   ( \sqrt{-B_*(c)} +\lvert \xi_1 \rvert ),
\end{align*}
where $C_{\rm Tay}>0$ is as in \eqref{eq:CTay} and in the third inequality we have used \eqref{eq:compactbss}.
 By \eqref{eq:cdsd02},   after shrinking $r_*>0$ further if necessary, we obtain 
 \[
 \lvert \alpha_c - \mathfrak{a}_c \xi_1 \rvert \leq ( \mathfrak{a}_c - d_0   )\lvert \xi_1 \rvert + C_2 \sqrt{-B_*(c)},
 \]
where $C_2:=\sqrt{2} + \sqrt{2r_*C_1C_{\rm Tay}}.$ If $\alpha_c \xi_1 \leq0,$ then 
  $\lvert \alpha_c - \mathfrak{a}_c \xi_1 \rvert = \lvert \alpha_c \rvert + \mathfrak{a}_c \lvert \xi_1 \rvert,$ and hence $\lvert \alpha_c \rvert \leq -d_0 \lvert \xi_1 \rvert + C_2 \sqrt{-B_*(c) }  < 2C_2 \sqrt{-B_*(c)}.$ This implies that   $\alpha_c$ and $\xi_1$ have the same sign, whenever $\lvert \alpha_c \rvert > 2C_2 \sqrt{-B_*(c)}.$  To compare their absolute values, we estimate using the triangle inequality that 
\begin{align*}
d_0 \lvert \xi_1 \rvert- \frac{1}{2}\lvert \alpha_c \rvert \leq \lvert \alpha_c\rvert   \leq (2\mathfrak{a}_c - d_0 ) \lvert \xi_1 \rvert +\frac{1}{2}\lvert \alpha_c\rvert  .
\end{align*}
Set $A = \mathfrak{a}_{c_2}+1.$ By continuity, $\mathfrak{a}_c \leq A$ for all $c>c_2$ sufficiently close to $c_2.$ We then obtain
\begin{equation}\label{eq:compacrabldfolsu}
\frac{2d_0}{3}\lvert \xi_1 \rvert \leq \lvert \alpha_c \rvert \leq ( 4A - 2d_0) \lvert \xi_1 \rvert,
\end{equation}
whenever $\lvert \alpha_c \rvert > 2C_2 \sqrt{-B_*(c)}.$
Note that the coefficients are uniform in $c.$

Choose   $\mathfrak{n}>5/(2d_0)$ sufficiently large and then choose $\rho$ and $\varrho$ sufficiently small as in Section \ref{sec:4} with $\varrho<\rho.$  Since $B_*(c) \to 0 $ as $c \to c_2^+, $ see Remark \ref{remark:valuev1xstar},  we may choose  $0<\varepsilon\leq \varrho^3$ so small that   $2C_2 \sqrt{-B_*(c)} <  \varrho/ \mathfrak{n} $ for $c \in (c_2, c_2+\varepsilon).$ On $\mathfrak{D}_c \cap \mathcal{B}_{r_*}(z_*(c)),$ we have   $\tilde{K}_c = K_c $ whenever $\alpha_c \leq 2\varrho/\mathfrak{n}.$ Indeed, if $\xi_1 \geq 0,$ then  $d_0 \xi_1 \leq \alpha_c + C_2 \sqrt{-B_*(c)} < {5\varrho}/({2\mathfrak{n}})<d_0\varrho,$ 
and thus, $\xi_1 < \varrho$ in this region. Hence, the capping vanishes.

Recall that  the first transition region is given by $-2\rho \leq \alpha_c \leq -\rho,$ where $\tilde{Y}_{\rm loc}$ interpolates  between $\hat{Y}$ and $X_c,$ while the second transition region is given by $\varrho / \mathfrak{n} \leq \alpha_c \leq 2 \varrho / \mathfrak{n},$ where $\bar{Y}_{\rm loc}$ interpolates between $X_c$ and $\mathcal{Z}_{c,\varrho}.$ See Figure \ref{fig:secondinterpolation}. Since $\varrho / \mathfrak{n} <\rho,$ both strips satisfy $\lvert \alpha_c \rvert  \geq \varrho/ \mathfrak{n} > 2C_2 \sqrt{-B_*(c)}.$ Hence,  by the above discussion,  $\alpha_c$ and $\xi_1$ have the same sign. So,    in the first strip, $\xi_1<0,$ and  by \eqref{eq:compacrabldfolsu} 
\[
C_3 \rho  \leq C_3  \lvert \alpha_c\rvert  \leq \lvert \xi_1 \rvert \leq \lvert \zeta \rvert
\]
with  $C_3:= (4A - 2d_0)^{-1} .$ 
Also, by the first inequality in \eqref{eq:compactbss} we have 
\[
\lvert \zeta \rvert^2 \leq C_1(-B_*(c) + \alpha_c^2) \leq C_1 (  (4C_2^2)^{-1}\rho^2   + 4\rho^2).
\]
Consequently, with    $C_4:=\sqrt{C_1 (   (4 C_2^2)^{-1} + 4 )},$
we have 
\[
\xi_1 \leq - C_3 \rho, \quad C_3 \rho \leq \lvert \zeta \rvert \leq C_4 \rho \quad \text{ on } \; \mathfrak{D}_c \cap \{ -2\rho \leq \alpha_c \leq -\rho\}.
\]
Then arguing as in  Lemma \ref{lem:firstinterpolatoinsuccess} we have 
\[
\tilde{Y}_{\rm loc}(\tilde{K}_c) \geq \rho^2 ( a_1 - b_1 \rho)>0 
\]
for some uniform constants $a_1,b_1>0$ and all sufficiently small $\rho>0.$ Similarly, in the second strip, we have
\[
\frac{ C_3 \varrho}{\mathfrak{n}} \leq \xi_1 \leq \lvert \zeta \rvert \leq \frac{C_4 \varrho}{\mathfrak{n}}   \quad \text{ on } \; \mathfrak{D}_c \cap \{ \frac{\varrho}{\mathfrak{n}} \leq \alpha_c \leq \frac{2\varrho}{\mathfrak{n}}\} .
\]
Arguing as in Lemma \ref{lem:latinterpolation} we find
\[
\bar{Y}_{\rm loc}(\tilde{K}_c) \geq \frac{\varrho^2}{\mathfrak{n}^2}(a_2 - \frac{b_2}{\mathfrak{n}})>0
\]
for some uniform constants $a_2, b_2>0,$ provided $\mathfrak{n}$ is sufficiently large and $\varrho$ is then chosen sufficiently small.

We now consider the remaining regions. On $  \{ \alpha_c \leq -2\rho\},$ where $\hat{Y}=Y_1, $ see Remark \ref{remark:slightlyaboveL2}, we also have $\lvert \alpha_c \rvert \geq 2\rho >2 C_2 \sqrt{-B_*(c)}.$ The preceding discussion shows that $\xi_1<0$ and $\lvert \alpha_c \rvert \leq  C_3^{-1} \lvert \xi_1 \rvert. $
Then the first inequality in 
\eqref{eq:compactbss} gives
\[
\lvert \zeta \rvert^2 \leq C_1 \left(  (4C_2^2)^{-1}+1 \right)\alpha_c^2 \leq C_1 C_3^{-2} \left( (4C_2^2)^{-1}+1\right)\lvert \xi_1 \rvert^2.
\]
We now recall from \eqref{eq:TTT} that since $T(c_2)<0$ and $T$ is continuous, there is $\sigma>0$ such that $T(c)<-\sigma<0$ uniformly for $c \in (c_2, c_2+\varepsilon),$ after taking smaller $\varepsilon>0$ if necessary. Since $\xi_1<0,$ equation \eqref{eq:firstesimatedsdfs} now gives
\[
Y_1(K_c) \geq T(c)\xi_1 - (C'+C'')\lvert \zeta\rvert^2 \geq \lvert \xi_1 \rvert ( \sigma - C_* \lvert \xi_1 \rvert) ,
\]
where $C_* := (C'+C'')C_1 C_3^{-2}(  (4C_2^2)^{-1}+1)>0.$ After decreasing $r_*>0$ if necessary, we may assume that $\sigma >C_* r_*,$ so that $Y_1(K_c)>0.$

Between the two  strips, where  $\bar{Y}_{\rm loc} = X_c,$  Lemma \ref{lem:Weinstein} gives $X_c(K_c) \geq \mathfrak{w}\lvert \zeta \rvert^2,$ which is positive   away from $z_*(c).$ 

On $  \{ \alpha_c>2\varrho / \mathfrak{n}\},$   the local Liouville vector field is $  \mathcal{Z}_{c,\varrho}.$  Assume first that $\xi_1 < \varrho.$ Then  $\alpha_c \geq 2\varrho/\mathfrak{n}>2C_2\sqrt{-B_*(c)}.$ By the preceding estimate,   both  $\xi_1$ and $\alpha_c$ are positive. We also have $\xi_1 \geq C_3 \alpha_c \geq 2 C_3  {\varrho}/{\mathfrak{n}}.$ The second inequality in \eqref{eq:compactbss} gives 
\[
\lvert \zeta\rvert^2 \leq C_1 \left( 1 + (4C_2^2 \mathfrak{n}^2)^{-1}\right) \varrho^2 \leq C_1 \left( 1 + (4C_2  ^2)^{-1} \right) \varrho^2. 
\]
As in the proof of Lemma \ref{lem:zoalsosmallstrip}, the lower bound on $\xi_1$ controls the term $(\xi_1 - 3\varrho)V_{c,1}'(x_*(c)+\xi_1),$ and the upper bound on $\lvert \zeta \rvert$ controls the error term $\xi_2^6. $ This therefore gives $\mathcal{Z}_{c,\varrho}(\tilde{K}_c)  >0$ for sufficiently small $\varrho>0.$ We now suppose that $\xi_1 \geq \varrho.$ Since $\xi_1 < 4 \varrho$ by construction of the capping, the second inequality of \eqref{eq:compactbss} gives $\lvert \zeta \rvert \leq C_5 \varrho$ on this part of $\mathfrak{D}_c,$ so    the estimate in the proof of Proposition \ref{prop:radislve} gives $\mathcal{Z}_{c,\varrho}(\tilde{K}_c)>0 $ away from the minimum $q_{c,\varrho}.$

Finally, let  $\mathfrak{f}$ be the cutoff function of Section \ref{sec:globalextn}. On $\{ -5\rho/2 \leq \alpha_c \leq 2\varrho/\mathfrak{n}\},$ we have $\lvert \alpha_c \rvert \leq 5\rho/2, $ 
so the first inequality of \eqref{eq:compactbss} gives
\[
\lvert \zeta \rvert^2 \leq C_1 \left( \frac{1}{4C_2^2} + \frac{25}{4} \right) \rho^2.
\]
On $\{ \alpha_c \geq 2\varrho / \mathfrak{n}\},$ the preceding discussion shows that $\lvert \zeta \rvert < C \varrho$ for some $C>0.$ Therefore, for $\rho$ sufficiently small, 
\[
\mathfrak{D}_c^0:=\mathfrak{D}_c \cap \mathcal{B}_{r_*}(z_*(c)) \cap \{ \alpha_c \geq -\frac{5\rho}{2} \} \subset \mathcal{B}_{r_*/2}(z_*(c)).
\]
Taking a larger neighborhood $\mathcal{U}_{\mathcal{W}}$ of $\mathcal{W}_{c,\rho,\varrho}$ in $\mathcal{B}_{r_*}(z_*(c))$ if necessary,   choose the cutoff $\mathfrak{f}$ of Section \ref{sec:globalextn} to be one on a neighborhood of $\bar{\mathcal{B}}_{r_*/2}(z_*(c)),$ and hence on a neighborhood of $\mathfrak{D}_c^0.$ By definition of $\mathcal{M}_{c,\rho,\varrho},$ we have $\alpha_c \leq -9\rho/4 <-2\rho$ on $\mathcal{M}_{c,\rho,\varrho}\cap{\rm supp}\mathfrak{f}.$ Hence, $\mathcal{F}_c$ vanishes on a neighborhood of this intersection, so $\bar{Y}_{\rm loc} = \hat{Y}$ there. 
On   $(\mathfrak{D}_c \cap \mathcal{B}_{r_*}(z_*(c)) ) \setminus \mathfrak{D}_c^0,$ we have $\alpha_c < -5\rho/2 <-2\rho, $ so  $\mathcal{F}_c$ vanishes on a neighborhood of each such point. Since $\mathfrak{f}=1$ near $\mathfrak{D}_c^0,$ it follows that  $Y = \bar{Y}_{\rm loc}$ on $\mathfrak{D}_c \cap \mathcal{B}_{r_*}(z_*(c)).$ 
This completes the proof of Step 1.

 \medskip

 \noindent
 {\bf Step 2.} We construct  the Weinstein cobordism. 

 \medskip

In each coordinate neighborhood, the modification of $\hat{Y}$ is supported in $\{\alpha_c \geq -2\rho\},$ and the proof of Step 1 gives the bound $\lvert \zeta \rvert \leq C \rho$ on $\mathfrak{D}_c \cap \{\alpha_c \geq -2\rho\}$ for some $C>0.$ Since $\tilde{K}_c \geq K_c \geq B_*(c) - C' \lvert \zeta \rvert^2$ for some $C'>0$ and $ 2C_2 \sqrt{-B_*(c)} <\varrho / \mathfrak{n} <\rho$ by the choices in Step 1, we have
$\tilde{K}_c \geq - ( ({4C_2^2})^{-1}  + C^2 C' ) \rho^2 $ throughout  $\mathfrak{D}_c \cap \{\alpha_c \geq -2\rho\}$ in each coordinate neighborhood.    By property iv) of Lemma \ref{lem:cappingfunction}, we have $0<d_{c,\varrho}\leq \mathfrak{c}_1\varrho^2.$ Fix $\mathfrak{C}> \max \{  (4C_2^2)^{-1}+C^2C', \mathfrak{c}_1\}$ 
and set $h:=-\mathfrak{C}\rho^2.$ Then $\tilde{K}_c>h$ throughout this set and, in particular, on the support of $Y - \hat{Y}$ in $\mathfrak{D}_c.$ 
Moreover, since $\varrho <\rho$ and the modified potential strictly decreases from $x_*(c)$ to $x_*(c)+3\varrho,$ we have
\begin{equation}\label{eq:hdcompare}
 h   < \tilde{K}_c (q_{c,\varrho})=-d_{c,\varrho} < \tilde{K}_c(z_*(c)) = B_*(c)<0 .
\end{equation}

We now choose $\rho$ sufficiently small, so that $-c_2/4 - 1/2 < h.$ Then for every $c \in (c_2, c_2+\varepsilon),$ the critical values $-1/2$ at $(x_1,x_2)=(0,0), (0,\pi)$ and $-c/4 -1/2$ at $(x_1,x_2)=(0,\pi/2), (0,3\pi/2)$ satisfy $-1/2 < -c/4 -1/2 <h.$

Recall the circle  $\mathcal{C} =\{x_1 = y_1 = y_2 =0\} $ from Remark \ref{rem:z2compute}. Then $\max_{\mathcal{C}} \tilde{K}_c = -c/4 - 1/2 <h.$ Since $h$ lies between the consecutive critical values $-c/4-1/2$ and $-d_{c,\varrho},$ it is a regular value of $\tilde{K}_c.$ Set
\[
\mathfrak{D}_c^- := \mathfrak{D}_c \cap \{ \tilde{K}_c \leq h \}.
\]
Each component of $\mathfrak{D}_c^-$ is compact and contains a point where $\tilde{K}_c$ attains its minimum on that component. Since $h$ is a regular value, this minimum is attained in the interior. 
The only local minima of $\tilde{K}_c$ below $h$ are $(0,0,0,0)$ and $(0,\pi,0,0).$ Both lie on the circle $\mathcal{C} \subset \mathfrak{D}_c^-$ and hence belong to the same component. Since every component contains a minimum, this implies that $\mathfrak{D}_c^-$ is connected.   Set  
\[
W_c := \mathfrak{D}_c \setminus {\rm int}(\mathfrak{D}_c^-).
\]
Since $\tilde{K}_c>h$ on the support of $Y-\hat{Y}$ in $\mathfrak{D}_c,$ we have $Y = \hat{Y}$ in a neighborhood of $\mathfrak{D}_c^-.$ 
The proof of  Lemma \ref{lem:locallemman} also applies to $\mathfrak{D}_c$ since $K_c \leq \tilde{K}_c \leq 0$ there and hence $p_1^2 +p_2^2 \leq -2V_c.$ Hence, by symmetry, we may choose $\delta_0>0$ such that, for every $c>c_2$ sufficiently close to $c_2$
\[
\mathfrak{D}_c \cap \{ \lvert \lvert x_1 \rvert - x_*(c) \rvert <\delta_0 \} \subset \cup_{\lambda \in \Lambda} \lambda( \mathcal{B}_{r_*}(z_*(c))).
\]
We then take $4\varrho<\delta_0.$
Now  Lemma \ref{lem:additionalestimateforhaty} gives $\hat{Y}(K_c) \geq \sigma_0$ on the center part $\mathcal{K}_c.$ Since $h \to 0$ as $\rho \to 0,$ compactness and continuity give $\hat{Y}(K_c) \geq \sigma_0/2$ whenever $\lvert \xi_1 \rvert \leq x_*(c)-\delta_0$ and $h < K_c \leq 0,$ after decreasing $\rho$ and $\varepsilon$ if necessary. On $W_c  \setminus \cup_{\lambda \in \Lambda} \lambda (\mathcal{B}_{r_*}( z_*(c))),$ we have $Y = \hat{Y}$ and $\tilde{K}_c=K_c,$ and hence
\[
Y(\tilde{K}_c) \geq \frac{\sigma_0}{2} >0 \qquad \text{ on } W_c  \setminus \cup_{\lambda \in \Lambda} \lambda (\mathcal{B}_{r_*}( z_*(c))).
\] 
Together with Step 1, this gives  $Y(\tilde{K}_c)>0$ on $W_c$ away from the four saddle points $\lambda(z_*(c))$ and the four minima $\lambda(q_{c,\varrho}),$ $\lambda \in \Lambda.$  At $z_*(c)$ the quadratic form $Q_c$ from  \eqref{eq:qudarictopart} has index one.
In the coordinates $(x_1, x_2,P_1,P_2),$ the Hessian of $\tilde{K}_c$ at $q_{c,\varrho}$ is diagonal. Its $x_1$-component is positive by  Lemma \ref{lem:cappingfunction},  and the $x_2$-component is $V_{c,2}''(0)>0.$ The kinetic energy $\frac{1}{2}(P_1^2+P_2^2)$ gives a diagonal entry of $1$ in each momentum coordinate.  Thus, the Hessian is positive-definite at $q_{c,\varrho}.$  By symmetry, the same holds for the other three minima.   Hence, $\tilde{K}_c|_{W_c}$ is Morse. For each $\lambda \in \Lambda,$ we have $Y=\lambda_*X_c$ near $\lambda(z_*(c))$ and $Y=\lambda_*\mathcal{Z}_{c,\varrho}$ near $\lambda(q_{c,\varrho}).$  
Then
Lemma \ref{lem:Weinstein}  and the quadratic bound in the proof of Proposition \ref{prop:radislve} show that  $Y$ is a gradient-like vector field for $\tilde{K}_c.$ Consequently, $(W_c, \omega, Y, \tilde{K}_c)$ is a Weinstein cobordism with $\partial^-W_c = \partial \mathfrak{D}_c^-$ and $\partial^+ W_c = \tilde{\Sigma}_c .$ The action of $\Lambda$ on $W_c$ is free, since every fixed point of a non-identity element lies on the circle $\mathcal{C} \subset {\rm int}(\mathfrak{D}_c^-).$
  Therefore, the quotient $ W_c/\Lambda$  is a Weinstein cobordism having exactly two critical points, the minimum $[q_{c,\varrho}]$ and the saddle $[z_*(c)]$ of index one. 
  This completes the proof of Step 2.

 \medskip

 \noindent
 {\bf Step 3.} We establish the uniqueness of the connecting trajectory and show that ${\rm cl}_{W_c} ( W_Y^s(z_*(c)))  \subset \mathcal{B}_{r_*}(z_*(c)).$

 \medskip


On the plane $\Pi := \{\xi_2 = \eta_1=0\},$ we have $E_c= \Omega_{c,\varrho}=0,$ see \eqref{eq:firstinterlE} and \eqref{eq:secondomnega}, so  
\[
\bar{Y}_{\rm loc} = (1-\ell_c)\hat{Y} + (\ell_c - g_{c,\varrho,\mathfrak{n}})X_c + g_{c,\varrho,\mathfrak{n}}\mathcal{Z}_{c,\varrho}.
\]
Since each of the three vector fields in the sum is tangent to $\Pi,$ we conclude that $\bar{Y}_{\rm loc}$ is tangent to $\Pi.$ Near $z_*(c),$ the Liouville vector field $\bar{Y}_{\rm loc} = X_c $ is linear. Its restriction to $\Pi$ has eigenvalues $a<0<1-b,$ with the $\eta_2$-axis as the expanding direction. Therefore, the two   branches of the local stable manifold of $z_*(c)$ lie in $\Pi,$ one in $\{\xi_1>0\}$ and the other in $\{\xi_1<0\}.$  
Along either branch, $\tilde{K}_c$ strictly increases to $B_*(c).$ On $\{ \xi_1=0\},$ we have $B_*(c) \leq \tilde{K}_c.$  Thus neither branch meets $\{ \xi_1=0\}.$ In other words, one lies entirely in $\{\xi_1>0\}$ and the other in $\{\xi_1<0\}.$

We first consider   the branch with $\xi_1>0.$   We claim that it converges to $q_{c,\varrho}$ in backward time. Since $0< \xi_1 < 4 \varrho$ (and we already have $2C_2\sqrt{-B_*(c)}<\varrho/\mathfrak{n}$ in Step 1), by the second inequality of \eqref{eq:compactbss} we have
\[
\lvert \zeta \rvert^2 \leq C_1 \left(  ({4C_2^2 \mathfrak{n}^2})^{-1}+16\right)\varrho^2 <  C_1 \left( ({4C_2^2 })^{-1}+16\right)\varrho^2.
\]
Then take smaller $\varrho>0$ if necessary, so that $\sqrt{ C_1 (   ({4C_2^2 })^{-1}+16)} \varrho <r_*. $  By the above inequality, the chosen branch cannot reach $\partial \mathcal{B}_{r_*}(z_*(c))$ while it remains in $W_c.$  Since $\xi_1>0,$ the modified potential attains its minimum at $q_{c,\varrho},$ so $\tilde{K}_c \geq -d_{c,\varrho}>h$ by \eqref{eq:hdcompare}.  Thus, the branch cannot reach the lower boundary either, so the backward trajectory stays in a compact subset of $W_c \cap \mathcal{B}_{r_*}(z_*(c)).$ Since $\tilde{K}_c$ strictly increases along nonconstant trajectories, its backward limit set consists of critical points with value less than $B_*(c).$  The only such point in $\mathcal{B}_{r_*}(z_*(c))$ is $q_{c,\varrho},$ proving the claim.

We now follow the  branch with $\xi_1<0$ in backward time. We claim that it reaches  $\partial^- W_c = \partial \mathfrak{D}_c^-$ in finite time without leaving $\mathcal{B}_{r_*}(z_*(c)).$ While $\alpha_c \geq -2\rho,$ the bound  $\lvert \zeta \rvert \leq C \rho$ from Step 1 prevents the branch from leaving $\mathcal{B}_{r_*}(z_*(c))$. By Step 2, the branch cannot reach the lower boundary $\{ \tilde{K}_c =h \}$ while $\alpha \geq -2\rho,$ since Step 2 gives   $\tilde{K}_c>h$ on this part of the branch.  Thus it must reach $\alpha = -2\rho$ in finite backward time. Indeed, otherwise, its backward trajectory from any point sufficiently close to the saddle remains in a compact subset of $\mathcal{B}_{r_*}(z_*(c)).$ Then its backward limit set contains a critical point with $\xi_1 \leq 0$ and value strictly less than $B_*(c),$ since $\tilde{K}_c $ decreases in backward time from a value already below $B_*(c).$ However, the only critical point in $\mathcal{B}_{r_*}(z_*(c))$ below $B_*(c)$ is $q_{c,\varrho},$ where $\xi_1 = 3\varrho>0.$

Set $t=0$ at the first intersection point with $\alpha = -2\rho.$ At $z_*(c)$  we have    $Y_1(\alpha_c)=\frac{1}{2}\mathfrak{a}_c \tanh(2x_*(c))>0$, so we may take $r_*>0$ small enough, so that  $Y_1(\alpha_c)>0$ on $\mathcal{B}_{r_*}(z_*(c)).$ Since $Y=Y_1$ on $\{ \alpha_c \leq -2\rho\}$ in this ball, the backward trajectory enters $\{\alpha_c \leq -2\rho\}$ and follows $Y_1$ there, with $\alpha_c$ decreasing as long as the trajectory remains in the coordinate ball.  
Along this trajectory in $\Pi,$ we have $\dot{x}_1  = \frac{1}{2}\tanh(2x_1)$ and $ \dot{p}_2 = p_2.$ 
Since $\lvert \zeta(0)\rvert \leq C \rho$ and $p_2= \eta_2 + A_1(x_*(c)+\xi_1 ) - A_1(x_*(c)),$   the mean value theorem gives $\lvert p_2(0)\rvert \leq C\rho$ possibly with different $C>0.$  Thus $\lvert p_2(t)\rvert  = e^t \lvert p_2(0)\rvert \leq C \rho$ for $t \leq 0$ on this trajectory. Since $x_2 = p_1=0,$ we have   $K_c = \frac{1}{2}p_2^2 + V_{c,1}(x_1)+V_{c,2}(0). $

Since $V'_{c,1}(x_*(c))=0$ and $V_{c,1}''(x_*(c))<0$ uniformly for $c \in (c_2, c_2+\varepsilon),$ by Taylor's theorem there is $\tilde{C}>0,$ independent of $c,$ such that
\[
V_{c,1}(x_1) \leq V_{c,1}(x_*(c)) - \tilde{C}( x_*(c)-x_1)^2
\]
for $x_1$ sufficiently close to $x_*(c).$ Since $\lvert p_2 \rvert \leq C\rho$ and $B_*(c)=V_{c,1}(x_*(c)) + V_{c,2}(0),$  
\[
K_c \leq B_*(c) + \frac{1}{2} C ^2\rho^2 -\tilde{C}(x_*(c)-x_1)^2.
\]
Choose $R>0$ so large that $\tilde{C}R^2 > \frac{1}{2} C ^2 + \mathfrak{C},$ where $\mathfrak{C}$ is as in Step 2. At $t=0,$ we have
$ x_*(c)-x_1(0) = \lvert \xi_1(0)\rvert \leq \lvert \zeta(0)\rvert \leq C \rho < R \rho.$ If the branch reaches $x_1 = x_*(c)-R\rho$ inside the coordinate ball, then the above inequalities give
\begin{equation}\label{eq:sdligd}
K_c \leq B_*(c) + (\frac{1}{2} C ^2   -\tilde{C}R^2) \rho^2 < - \mathfrak{C} \rho^2 = h 
\end{equation}
 for $\rho$ sufficiently small. 
For $x_*(c)-R\rho \leq x_1 \leq x_*(c),$ we have 
$\dot{x}_1 = \frac{1}{2}\tanh(2x_1) \geq  \frac{1}{2}\tanh( 2(x_*(c)-R\rho))>0,$ where the last inequality holds for $\rho$ sufficiently small. Thus, $x_1$ decreases in backward time and the branch can stay in this interval  only for a finite time while following $Y_1.$ On the interval  $x_*(c)-R\rho \leq x_1 \leq x_*(c),$ we have $\lvert \xi_1 \rvert = x_*(c)-x_1 \leq R \rho.$ For points of the branch in $W_c \cap \mathcal{B}_{r_*}(z_*(c))$ satisfying this bound,   \eqref{eq:compactbss} gives $\lvert \zeta \rvert^2 \leq C_1 ( (4C_2^2)^{-1}+R^2)\rho^2<r_*^2/4 $ for sufficiently small $\rho.$ 
This prevents the branch from leaving the coordinate ball.  Since $\xi_1<0$ along this branch, $\tilde{K}_c=K_c.$ Thus $K_c \geq h$ as long as the branch remains in $W_c.$ If the branch reaches $x_1 = x_*(c)-R\rho$ while still in $W_c,$ then \eqref{eq:sdligd} gives $\tilde{K}_c = K_c <h,$ a contradiction. Therefore, the branch meets the lower boundary in finite backward time at a point with $x_1 > x_*(c)-R\rho.$

Consequently, the branch  in $\{\xi_1>0\}$ is the unique trajectory from $q_{c,\varrho}$ to $z_*(c).$ Moreover, the closure of the stable manifold of $z_*(c)$ in $W_c$ lies in $\mathcal{B}_{r_*}(z_*(c)).$ By symmetry the same holds in the other three coordinate neighborhoods. These four connecting trajectories descend to a single trajectory in $W_c/\Lambda.$  Every connecting trajectory in the quotient lifts to one ending at $z_*(c),$ so the trajectory from $[q_{c,\varrho}]$ to $[z_*(c)]$ is unique up to time translation. 
This completes the proof of Step 3.

 \medskip

 \noindent
 {\bf Step 4.} We prove the proposition.

 \medskip

By Steps 2 and 3, Theorem \ref{thm:cancle} applies to $W_c / \Lambda$ to the critical points $[q_{c,\varrho}]$ and $[z_*(c)].$ By Step 3, the closure of the stable manifold of $[z_*(c)]$ lies in the image of $W_c \cap \mathcal{B}_{r_*}(z_*(c))$ under the quotient map $W_c \to W_c /\Lambda. $ Thus we can choose the homotopy, so that its lift to $W_c $  is supported in the four coordinate neighborhoods.  Note that the lift is $\Lambda$-invariant, fixes $\omega,$ and agrees with the original Weinstein structure $(\omega, Y, \tilde{K}_c)$ near both boundary components for every $t \in [0,1].$ At $t=1,$ we denote the lifted Liouville vector field and Morse function by $Y'$ and $\tilde{K}_c'.$ Then $d\tilde{K}_c'(Y')>0$ on $W_c.$ Since $Y'=Y$ near $\partial W_c,$ we extend $Y'$ by $Y$ outside $W_c.$ The resulting vector field is then Liouville with respect to $\omega$ and agrees with $Y$ near the upper boundary $\tilde{\Sigma}_c.$ 
Since the four coordinate neighborhoods are disjoint from $L_{\rm col}^E,$ we also have $Y'=Y$ near $L_{\rm col}^E.$

By compactness, we set $m:=\min_{W_c} d\tilde{K}_c'(Y')>0.$ Then a trajectory of the flow of $Y'$ which stays in $W_c$ for time $\tau$  satisfies
\begin{equation}\label{eq:timeesiat}
m \tau \leq \max_{W_c} \tilde{K}_c' - \min_{W_c} \tilde{K}_c'.
\end{equation}
Thus, every trajectory starting in ${\rm int}(W_c)$ reaches $\partial W_c$ in finite time in both forward and backward time. Since $Y'$ points into $W_c$ along the lower boundary $\partial \mathfrak{D}_c^-,$ the trajectory reaches $\partial \mathfrak{D}_c^- $ in backward time. Similarly, since $Y'$ points out of $W_c$ along the upper boundary $\tilde{\Sigma}_c,$ the trajectory reaches $\tilde{\Sigma}_c $ in forward time.

Recall the  reference energy  $d \in (c_1, c_2)$ fixed in Section \ref{secL:72}.  Let $\Sigma_d$ be the compact component of $K_d^{-1}(0) $ and let $\mathfrak{D}_d$ be the compact domain bounded by $\Sigma_d.$ Since the four coordinate neighborhoods supporting the modification are disjoint from $\Sigma_d,$ we have $Y'=Y=\hat{Y}$ near $\Sigma_d.$ 
Since $\min_{x_2} V_d(\pm x_*(c_2), x_2) = ((c_2-d)/4)\sinh^2(x_*(c_2))>0,$ we have $V_d>0$ on $\{ x_1= \pm x_*(c)\}$ for $c$ sufficiently close to $c_2.$  Then $\mathfrak{D}_d \subset \{ K_d \leq 0\}$ intersects $\{x_1=0\}$ and hence lies in $\{ -x_*(c) < x_1 < x_*(c)\}.$ In $\mathfrak{D}_d,$ the capping function vanishes, and 
\[
\tilde{K}_c = K_c=K_d - \frac{c-d}{4}( \cosh^2 x_1 - \cos^2 x_2) \leq K_d \leq 0,
\]
 implying that $\mathfrak{D}_d \subset \mathfrak{D}_c. $ A trajectory starting on $\tilde{\Sigma}_c$ enters $\mathfrak{D}_c^-$ in backward time and remains there. Since $Y'=Y =\hat{Y}$ there,   by Remark \ref{rem:z2compute} it converges backward to a point on the circle $\{ x_1 = y_1= y_2=0\},$ which lies in ${\rm int}(\mathfrak{D}_d). $ Therefore, every trajectory starting on $\tilde{\Sigma}_c$ meets $\Sigma_d$ in backward time, possibly at time zero.
 
 Conversely, let $z \in \Sigma_d$ and follow the trajectory of $Y'$ starting at $z$ in forward time. If $z \in \tilde{\Sigma}_c,$ the hitting time is zero. If $z \in W_c \setminus \tilde{\Sigma}_c,$ then \eqref{eq:timeesiat} shows that the trajectory reaches $\tilde{\Sigma}_c$ in finite time, using inward transversality if $z \in \partial \mathfrak{D}_c^-.$ Since $\mathfrak{D}_d \subset \mathfrak{D}_c,$ it remains to consider the case $z \in {\rm int}(\mathfrak{D}_c^-).$ Since $\{x_1 = y_1 = y_2=0\} \subset { \rm int}(\mathfrak{D}_d),$ the boundary point $z$ is not on this circle. Since $Y' = \hat{Y}$ on $\mathfrak{D}_c^-,$ Remark \ref{rem:z2compute} implies that the trajectory reaches $\partial \mathfrak{D}_c^-$ in finite forward time. Then it enters $W_c$  and reaches $\tilde{\Sigma}_c$ in finite time again by \eqref{eq:timeesiat}.
Since $Y'$ points strictly outward along both hypersurfaces $\Sigma_d$ and $\tilde{\Sigma}_c,$ each intersection is unique. 
This completes the proof of Step 4 and hence of the proposition.

\begin{remark}\label{rmk:functioanlnunchange}

The above modification from $Y$ to $Y'$ preserves the action functional. Indeed, for  $\Sigma=\Sigma_d$ or $\Sigma=\tilde{\Sigma}_c\},$ the two vector fields agree near $\Sigma, $ so we may use the same defining Hamiltonian $H$ before and after the modification. Set $\lambda= \iota_Y \omega$ and $\lambda' = \iota_{Y'}\omega.$ The closed one-form $\lambda' - \lambda$ is compactly supported in the four coordinate balls,  and therefore equal to $dF$ for  a function $F$ supported in their union. Since these balls are disjoint from $L_{\rm col}^E,$ both $F$ and $dF = \lambda'-\lambda$ vanish in a neighborhood of $L_{\rm col}^E.$  Hence, we may use the same function $\ell,$ with $d\ell = \lambda|_{L_{\rm col}^E},$ for the endpoint term $\ell(x(0))-\ell(x(1))$ in both action functionals. Thus, for every path $x$ with endpoints on $L_{\rm col}^E$ and every $\tau \in \R$ we have
\[
\mathcal{A}_H ( x, \tau;\lambda') - \mathcal{A}_H(x,\tau;\lambda) = \int_0^1 x^*dF = F(x(1)) - F(x(0))=0.
\]
The same calculation applies to the perturbed functional, since its Hamiltonian terms are independent of the primitive.

\end{remark}




\bibliographystyle{abbrv}
\bibliography{mybibfile}

\end{document}